\documentclass[11pt]{article}
\usepackage{amsfonts}
\usepackage{bbm}
\usepackage{mathrsfs}
\usepackage{amscd}
\usepackage{color}
\usepackage{amsmath,amsfonts,amssymb,amscd}
\usepackage{indentfirst,graphics,epsfig,psfrag}
\input{epsf}
\usepackage{ifpdf}
\usepackage{enumerate}
\usepackage{appendix}
\usepackage{enumerate}

\usepackage{lineno}

\makeatletter \@addtoreset{figure}{section} \makeatother
\makeatletter
\long\def\@makecaption#1#2{%
   \vskip 10\p@
   \setbox\@tempboxa\hbox{{#1}\ \ #2}%
   \ifdim \wd\@tempboxa >\hsize

       {#1}\ \ #2\par
   \else
       \hbox to\hsize{\hfil\box\@tempboxa\hfil}%
   \fi}
\makeatother

\newtheorem{thm}{Theorem}[section]
\newtheorem{cor}[thm]{Corollary}
\newtheorem{lem}[thm]{Lemma}
\newtheorem{rem}[thm]{Remark}
\newtheorem{con}[thm]{Conjecture}

\newtheorem{obs}[thm]{Observation}
\newtheorem{pro}[thm]{Proposition}
\newtheorem{prob}[thm]{Problem}

\begin{document}
\title{\textbf{A survey on the\\ generalized connectivity of graphs}\footnote{Supported by NSFC Nos.11371205 and 11531011.}}
\author{
\small  Xueliang Li, \ Yaping Mao\\[8pt]
\small Center for Combinatorics and LPMC-TJKLC\\
\small Nankai University, Tianjin 300071, P.R. China\\
\small E-mails: lxl@nankai.edu.cn; maoyaping@ymail.com
 }
\date{}
\maketitle
\begin{abstract}
The generalized $k$-connectivity $\kappa_k(G)$ of a graph $G$ was
introduced by Hager before 1985. As its a natural counterpart, we
introduced the concept of generalized edge-connectivity
$\lambda_k(G)$, recently. In this paper we summarize the known
results on the generalized connectivity and generalized
edge-connectivity. After an introductory section, the
paper is then divided into nine sections: the generalized (edge-)connectivity of some graph
classes, algorithms and computational complexity, sharp bounds of
$\kappa_k(G)$ and $\lambda_k(G)$, graphs with large generalized
(edge-)connectivity, Nordhaus-Gaddum-type results, graph operations,
extremal problems, and some results for random graphs and
multigraphs. It also contains some conjectures and open problems
for further studies. \\[2mm]
\end{abstract}

\section{Introduction}

In this introductory section, we will give both theoretical and
practical motivation for introducing the concept of generalized
(edge-)connectivity of graphs. Some useful definitions on graph theory are
also given. It is divided into the following five subsections.

\subsection{Connectivity and its generalizations}

Connectivity is one of the most basic concepts of graph-theoretic
subjects, both in a combinatorial sense and an algorithmic sense. As
we know, the classical connectivity has two equivalent definitions.
The \emph{connectivity} of $G$, written $\kappa(G)$, is the minimum
order of a vertex set $S\subseteq V(G)$ such that $G\setminus S$ is
disconnected or has only one vertex. We call this definition the
``cut" version definition of connectivity. A well-known theorem of
Whitney provides an equivalent definition of connectivity, which can
be called the ``path" version definition of connectivity. For any
two distinct vertices $x$ and $y$ in $G$, the \emph{local
connectivity} $\kappa_{G}(x,y)$ is the maximum number of internally
disjoint paths connecting $x$ and $y$. Then
$\kappa(G)=\min\{\kappa_{G}(x,y)\,|\,x,y\in V(G), \ x\neq y\}$ is
defined to be the \emph{connectivity} of $G$. In contrast to this
parameter, $\overline{\kappa}(G)=\max\{\kappa_{G}(x,y)\,|\,x,y\in
V(G), \ x\neq y\}$, first introduced by Bollob\'{a}s (see
\cite{Bollobas1} for example), is called the \emph{maximum local
connectivity} of $G$. As we have seen, the connectivity and maximum
local connectivity are two extremes of the local connectivity of a
graph. An invariant lying between these two extremes is the {\it
average connectivity} $\widehat{\kappa}(G)$ of a graph, which is
defined to be $\widehat{\kappa}(G)=\sum_{x,y\in
V(G)}\kappa_{G}(x,y)/{n \choose 2}$; see \cite{Beineke2}.

Similarly, the classical edge-connectivity also has two equivalent definitions.
The \emph{edge-connectivity} of $G$, written $\lambda(G)$, is the
minimum size of an edge set $M\subseteq E(G)$ such that $G\setminus
M$ is disconnected. We call this definition
the ``cut" version definition of edge-connectivity. Whitney also
provided an equivalent definition of edge-connectivity, which can be
called the ``path" version definition. For any two distinct vertices
$x$ and $y$ in $G$, the \emph{local edge-connectivity}
$\lambda_{G}(x,y)$ is the maximum number of edge-disjoint
paths connecting $x$ and $y$. Then
$\lambda(G)=\min\{\lambda_{G}(x,y)\,|\,x,y\in V(G), \ x\neq y\}$,
$\overline{\lambda}(G)=\max\{\lambda_{G}(x,y)\,|\,x,y\in V(G), \
x\neq y\}$ and $\widehat{\lambda}(G)=\sum_{x,y\in
V(G)}\kappa_{G}(x,y)/{n \choose 2}$ are the {\it edge-connectivity,
maximum local edge-connectivity} and {\it average
edge-connectivity}, respectively. For connectivity and
edge-connectivity, Oellermann gave a survey paper on this subject;
see \cite{Oellermann2}.

\subsubsection{$k$-Connectivity and $k$-edge-connectivity}

Although there are many elegant and powerful results on connectivity
in Graph Theory, the classical connectivity and edge-connectivity
cannot be satisfied considerably in practical uses. So
people tried to generalize these concepts. For the ``cut" version
definition of connectivity, one can see that the above minimum
vertex set does not regard to the number of components of
$G\setminus S$. Two graphs with the same connectivity may have
different degrees of vulnerability in the sense that the deletion of
a vertex cut-set of minimum cardinality from one graph may produce a
graph with considerably more components than in the case of the
other graph. For example, the star $K_{1,n}$ and the path $P_{n+1}\
(n\geq 3)$ are both trees of order $n+1$ and therefore connectivity
$1$, but the deletion of a cut-vertex from $K_{1,n}$ produces a
graph with $n$ components, while the deletion of a cut-vertex from
$P_{n+1}$ produces only two components. The above statement suggests
a generalization of the connectivity of a graph. In 1984, Chartrand
et al. \cite{Chartrand1} generalized the ``cut" version definition
of connectivity. For an integer $k \ (k\geq 2)$ and a graph $G$ of
order $n \ (n\geq k)$, the \emph{$k$-connectivity} $\kappa'_k(G)$ is
the smallest number of vertices whose removal from $G$ produces a
graph with at least $k$ components or a graph with fewer than $k$
vertices. Thus, for $k=2$, $\kappa'_2(G)=\kappa(G)$. For more
details about the $k$-connectivity, we refer to \cite{Chartrand1,
Day, Oellermann2, Oellermann3}.

If two graphs have the same edge-connectivity, then the removal of
an edge set of minimum cardinality from either graph produces
exactly two components. On the other hand, disconnecting these
graphs into three components may require the removal of considerably
more edges in the one case than the other. Take for example, if
$H_1$ is obtained from two copies of complete graph $K_n \ (n\geq
2)$ by joining two vertices (one from each copy of $K_n$) by an edge
and $H_2$ is a path of order $2n$, then both graphs have order $2n$
and edge-connectivity $1$. However, $n$ edges need to be removed
from $H_1$ but only two edges from $H_2$ to produce a graph with three
components. This observation suggests a generalization of the ``cut"
version definition of classical edge-connectivity. For an integer $k
\ (k\geq 2)$ and a graph $G$ of order $n \ (n\geq k)$, the
\emph{$k$-edge-connectivity} $\lambda'_k(G)$ is the smallest number
of edges whose removal from $G$ produces a
graph with at least $k$ components. Thus, for $k=2$,
$\lambda'_2(G)=\lambda(G)$. The $k$-edge-connectivity was initially
introduced by Boesch and Chen \cite{Boesch} and subsequently studied
by Goldsmith in \cite{Goldsmith1, Goldsmith2} and Goldsmith et al.
\cite{Goldsmith3}. In all these papers, the computational difficulty
of finding $\lambda'_k(G)$ for $k\geq 3$ leads to the development
of heuristics and bounds for approximating this parameter. For more
details on $k$-edge-connectivity, we refer to \cite{Beineke1,
Oellermann1}.

\subsubsection{Generalized $k$-connectivity and generalized $k$-edge-connectivity}

The generalized connectivity of a graph $G$, introduced by Hager, is
a natural generalization of the ``path" version definition of
connectivity. For a graph $G=(V,E)$ and a set $S\subseteq V(G)$ of
at least two vertices, \emph{an $S$-Steiner tree} or \emph{a Steiner
tree connecting $S$} (or simply, \emph{an $S$-tree}) is a such
subgraph $T=(V',E')$ of $G$ that is a tree with $S\subseteq V'$.
Note that when $|S|=2$ a minimal Steiner tree
connecting $S$ is just a path connecting the two vertices of $S$.
Two Steiner trees $T$ and $T'$ connecting $S$ are said to be
\emph{internally disjoint} if $E(T)\cap E(T')=\varnothing$ and
$V(T)\cap V(T')=S$. For $S\subseteq V(G)$ and $|S|\geq 2$, the
\emph{generalized local connectivity} $\kappa_G(S)$ is the maximum
number of internally disjoint Steiner trees connecting $S$ in $G$,
that is, we search for the maximum cardinality of edge-disjoint
Steiner trees which contain $S$ and are vertex-disjoint with the
exception of the vertices in $S$. For an integer $k$ with $2\leq
k\leq n$, the \emph{generalized $k$-connectivity} (or
\emph{$k$-tree-connectivity}) is defined as
$\kappa_k(G)=\min\{\kappa_G(S)\,|\,S\subseteq V(G),|S|=k\}$, that
is, $\kappa_k(G)$ is the minimum value of $\kappa_G(S)$ when $S$
runs over all $k$-subsets of $V(G)$. Clearly, when $|S|=2$,
$\kappa_2(G)$ is just the connectivity $\kappa(G)$ of $G$, that is,
$\kappa_2(G)=\kappa(G)$, which is the reason why one addresses
$\kappa_k(G)$ as the generalized connectivity of $G$. By convention,
for a connected graph $G$ with less than $k$ vertices, we set
$\kappa_k(G)=1$, and $\kappa_k(G)=0$ when $G$ is disconnected.
Note that the generalized
$k$-connectivity and the $k$-connectivity of a graph are indeed
different. Take for example, the graph $G_0$ obtained from a
triangle with vertex set $\{v_1,v_2,v_3\}$ by adding three new
vertices $u_1,u_2,u_3$ and joining $v_i$ to $u_i$ by an edge for $1
\leq i\leq 3$. Then $\kappa_3(G_0)=1$ but $\kappa'_3(G_0)=2$. We
knew this concept in \cite{Chartrand2} for the first time. There the
authors obtained the exact value of the generalized $k$-connectivity
of complete graphs. Recently, from \cite{Hager, Hager2}, we know
that the concept was introduced actually by Hager in his another
paper, but we do not know whether his this paper has been published,
yet. For results on the generalized connectivity (or
tree-connectivity), we refer to \cite{Chartrand2, CLLM, CLL, Gao,
GLS, LLM, LLM2, LLMS, LLMY, LLSun, SLi, LLL1, LLL2, LL, LLShi, LLZ,
WLi, LM, LLZ, LM3, LM4, LMS, LMW2, LYZ, LZ, Okamoto, SL, SZ}.

The following Table 1 shows how the generalization proceeds.

\begin{tabular}{|c|c|c|}
\hline $$ & $Classical~connectivity$& $Generalized~
connectivity$\\
\cline{1-3}
$Vetex~subset $ & $S=\{x,y\}\subseteq V(G) \ (|S|=2)$ & $S\subseteq V(G) \ (|S|\geq 2))$\\
\cline{1-3} $Set~of~Steiner~trees$ & $\left\{
\begin{array}{ll}
\mathscr{P}_{x,y}=\{P_1,P_2,\cdots,P_{\ell}\}\\
\{x,y\}\subseteq V(P_i),\\
E(P_i)\cap
E(P_j)=\varnothing\\
V(P_i)\cap V(P_j)=\{x,y\}\\
\end{array}
\right.$ & $\left\{
\begin{array}{ll}
\mathscr{T}_{S}:T_1,T_2,\cdots,T_{\ell}\\
S\subseteq V(T_i),\\
E(T_i)\cap E(T_j)=\varnothing,\\
V(T_i)\cap
V(T_j)=S\\
\end{array}
\right.$\\
\cline{1-3}
$Local~parameter$ & $\kappa(x,y)=\max|\mathscr{P}_{x,y}|$ & $\kappa(S)=\max|\mathscr{T}_{S}|$\\
\cline{1-3} $Global~parameter$ & $\kappa(G)=\underset{x,y\in
V(G)}{\min} \kappa(x,y)$ & $\kappa_k(G)=\underset{S\subseteq V(G),
|S|=k}{\min} \kappa(S)$\\
\cline{1-3}
\end{tabular}
\begin{center}
{Table~1.~Classical~connectivity and generalized connectivity}
\end{center}

As a natural counterpart of the generalized connectivity, we
introduced the concept of generalized edge-connectivity in
\cite{LMS}. For $S\subseteq V(G)$ and $|S|\geq 2$, the
\emph{generalized local edge-connectivity} $\lambda(S)$ is the
maximum number of edge-disjoint Steiner trees connecting $S$ in $G$.
For an integer $k$ with $2\leq k\leq n$, the \emph{generalized
$k$-edge-connectivity} $\lambda_k(G)$ of $G$ is then defined as
$\lambda_k(G)= min\{\lambda(S)\,|\,S\subseteq V(G) \ and \ |S|=k\}$.
It is also clear that when $|S|=2$, $\lambda_2(G)$ is just the
standard edge-connectivity $\lambda(G)$ of $G$, that is,
$\lambda_2(G)=\lambda(G)$, which is the reason why we address
$\lambda_k(G)$ as the generalized edge-connectivity of $G$. Also set
$\lambda_k(G)=0$ when $G$ is disconnected. Results on the
generalized edge-connectivity can be found in \cite{LM, LM4, LMS,
LMW2, LYZ}.

The following Table 2 shows how the generalization of the
edge-version definition proceeds.

\begin{tabular}{|c|c|c|}
\hline $$ & $Edge$-$connectivity$& $Generalized~
edge$-$connectivity$\\
\cline{1-3}
$Vetex~subset $ & $S=\{x,y\}\subseteq V(G) \ (|S|=2)$ & $S\subseteq V(G) \ (|S|\geq 2))$\\
\cline{1-3} $Set~of~Steiner~trees$ & $\left\{
\begin{array}{ll}
\mathscr{P}_{x,y}=\{P_1,P_2,\cdots,P_{\ell}\}\\
\{x,y\}\subseteq V(P_i),\\
E(P_i)\cap
E(P_j)=\varnothing\\
\end{array}
\right.$ & $\left\{
\begin{array}{ll}
\mathscr{T}_{S}:T_1,T_2,\cdots,T_{\ell}\\
S\subseteq V(T_i),\\
E(T_i)\cap E(T_j)=\varnothing,\\
\end{array}
\right.$\\
\cline{1-3}
$Local~parameter$ & $\lambda(x,y)=\max|\mathscr{P}_{x,y}|$ & $\lambda(S)=\max|\mathscr{T}_{S}|$\\
\cline{1-3} $Global~parameter$ & $\lambda(G)=\underset{x,y\in
V(G)}{\min} \lambda(x,y)$ & $\lambda_k(G)=\underset{S\subseteq V(G),
|S|=k}{\min} \lambda(S)$\\
\cline{1-3}
\end{tabular}
\begin{center}
{Table~2.~Classical~edge-connectivity and generalized
edge-connectivity}
\end{center}

\begin{rem} The difference between the ``path" version generalized connectivity $\kappa_k(G)$ and
the ``cut" version $k$-connectivity $\kappa'_k(G)$ was discussed very clearly by Sun and Li in \cite{SL}, where
they got sharp lower and upper bounds for the difference $\kappa'_k(G)-\kappa_k(G)$, and investigated the problem
that under what conditions for a graph $G$ one has $\kappa'_k(G)=\kappa_k(G)$.
\end{rem}

\subsubsection{Mader's generalization}

In fact, Mader \cite{Mader3} studied an extension of Menger's
theorem to independent sets of three or more vertices. We know that
from Menger's theorem that if $S=\{u,v\}$ is a set of two
independent vertices in a graph $G$, then the maximum number of
internally disjoint $u$-$v$ paths in $G$ equals the minimum number
of vertices that separate $u$ and $v$. For a set
$S=\{u_1,u_2,\cdots,u_k\}$ of $k \ (k\geq 2)$ vertices in a graph
$G$, an \emph{$S$-path} is defined as a path between a pair of
vertices of $S$ that contains no other vertices of $S$. Two
$S$-paths $P_1$ and $P_2$ are said to be \emph{internally disjoint}
if they are vertex-disjoint except for the vertices of $S$. If $S$
is a set of independent vertices of a graph $G$, then a vertex set
$U\subseteq V(G)$ with $U\cap S=\varnothing$ is said to
\emph{totally separate $S$} if every two vertices of $S$ belong to
different components of $G\setminus U$. Let $S$ be a set of at least
three independent vertices in a graph $G$. Let $\mu(G)$ denote the
maximum number of internally disjoint $S$-paths and $\mu'(G)$ the
minimum number of vertices that totally separate $S$. A natural
extension of Menger' s theorem may well be suggested, namely: If $S$
is a set of independent vertices of a graph $G$ and $|S|\geq 3$,
then $\mu(S)=\mu'(S)$. However, the statement is not true in
general. Take the above graph $G_0$ for example. For
$S=\{u_1,u_2,u_3\}$, $\mu(S)=1$ but $\mu'(S)=2$. Mader proved that
$\mu(S)\geq \frac{1}{2}\mu'(S)$. Moreover, the bound is sharp.
Lov\'{a}sz conjectured an edge analogue of this result and Mader
proved this conjecture and established its sharpness. For more
details, we refer to \cite{Mader3, Mader4, Oellermann1}.

\subsubsection{Pendant tree-connectivity and path-connectivity}

Except for the concept of tree-connectivity, Hager also introduced
another tree-connectivity parameter, called the {\it pendant
tree-connectivity} of a graph in \cite{Hager}. For the
tree-connectivity (or generalized connectivity),
we only search for edge-disjoint trees which include $S$ and are
vertex-disjoint with the exception of the vertices in $S$. But
pendant tree-connectivity further requests the degree of each vertex
of $S$ in a Steiner tree connecting $S$ is equal to one. Note that
it is a specialization of the generalized connectivity (or
tree-connectivity), but it is a generalization of the classical
connectivity. The detailed definitions are stated as follows. For an
$S$-Steiner tree, if the degree of each vertex in $S$ is equal to
one, then this tree is called a \emph{pendant $S$-Steiner tree}. Two
pendant $S$-Steiner trees $T$ and $T'$ are said to be
\emph{internally disjoint} if $E(T)\cap E(T')=\varnothing$ and
$V(T)\cap V(T')=S$. For $S\subseteq V(G)$ and $|S|\geq 2$, the
\emph{local pendant tree-connectivity} $\tau_G(S)$ is the maximum
number of internally disjoint pendant $S$-Steiner trees in $G$. For
an integer $k$ with $2\leq k\leq n$, the \emph{pendant
$k$-tree-connectivity} is defined as
$\tau_k(G)=\min\{\tau_G(S)\,|\,S\subseteq V(G),|S|=k\}$. Set
$\tau_k(G)=0$ when $G$ is disconnected. It is clear that
$$
\left\{
\begin{array}{ll}
\tau_k(G)=\kappa_k(G),&for~k=1,2;\\
\tau_k(G)\leq \kappa_k(G),&for~k\geq 3.
\end{array}
\right.
$$

Dirac \cite{Dirac} showed that in a $(k-1)$-connected graph there is
a path through each $k$ vertices. Related problems were inquired in
\cite{Wilson}. In \cite{Hager2}, Hager revised this statement to the
question of how many internally disjoint paths $P_i$ with the
exception of a given set $S$ of $k$ vertices exist such that
$S\subseteq V(P_i)$. Another concept of connectivity, the path-connectivity,
of a graph $G$ was also introduced by Hager in \cite{Hager2}, which
is a specialization of the generalized connectivity and
is also a generalization of the ``path" version definition of the classical connectivity.
For a graph $G=(V,E)$ and a set $S\subseteq V(G)$ of at least two vertices, \emph{a path
connecting $S$} (or simply, \emph{an $S$-path}) is a subgraph
$P=(V',E')$ of $G$ that is a path with $S\subseteq V'$. Note that an
$S$-path is also a tree connecting $S$. Two $S$-paths $P$ and $P'$
are said to be \emph{internally disjoint} if $E(P)\cap
E(P')=\varnothing$ and $V(P)\cap V(P')=S$. For $S\subseteq V(G)$ and
$|S|\geq 2$, the \emph{local path-connectivity} $\pi_G(S)$ is the
maximum number of internally disjoint $S$-paths in $G$, that is, we
search for the maximum cardinality of edge-disjoint paths which
contain $S$ and are vertex-disjoint with the exception of the
vertices in $S$. For an integer $k$ with $2\leq k\leq n$, the
\emph{$k$-path-connectivity} of a graph $G$ on $n$ vertices is defined as
$\pi_k(G)=\min\{\pi_G(S)\,|\,S\subseteq V(G),|S|=k\}$, that is,
$\pi_k(G)$ is the minimum value of $\pi_G(S)$ when $S$ runs over all
$k$-subsets of $V(G)$. Clearly, we have
$$
\left\{
\begin{array}{ll}
\pi_k(G)=\delta(G),&for~k=1;\\
\pi_k(G)=\kappa(G),&for~k=2;\\
\pi_k(G)\leq \kappa_k(G),&for~k\geq 3.
\end{array}
\right.
$$

The relations between the pendant tree-connectivity,
generalized connectivity and path-connectivity are shown in the
following Table 3.

\begin{tabular}{|c|c|c|}
\hline $Pendant~tree$-$connectivity$ & $Generalized~
connectivity$& $Path$-$connectivity$\\
\cline{1-3}
$S\subseteq V(G) \ (|S|\geq 2))$ & $S\subseteq V(G) \ (|S|\geq 2))$ & $S\subseteq V(G) \ (|S|\geq 2))$\\
\cline{1-3} $\left\{
\begin{array}{ll}
\mathscr{T}_{S}:T_1,T_2,\cdots,T_{\ell}\\
S\subseteq V(T_i),\\
d_{T_i}(v)=1~for~every~v\in S\\
E(T_i)\cap E(T_j)=\varnothing,\\
\end{array}
\right.$ & $\left\{
\begin{array}{ll}
\mathscr{T}_{S}:T_1,T_2,\cdots,T_{\ell}\\
S\subseteq V(T_i),\\
E(T_i)\cap E(T_j)=\varnothing,\\
\end{array}
\right.$ & $\left\{
\begin{array}{ll}
\mathscr{P}_{S}:P_1,P_2,\cdots,P_{\ell}\\
S\subseteq V(P_i),\\
E(P_i)\cap E(P_j)=\varnothing,\\
\end{array}
\right.$\\
\cline{1-3}
$\tau(S)=\max|\mathscr{T}_{S}|$ & $\kappa(S)=\max|\mathscr{T}_{S}|$ & $\pi(S)=\max|\mathscr{P}_{S}|$\\
\cline{1-3} $\tau_k(G)=\underset{S\subseteq V(G), |S|=k}{\min}
\tau(S)$ & $\kappa_k(G)=\underset{S\subseteq V(G), |S|=k}{\min}
\kappa(S)$ & $\pi_k(G)=\underset{S\subseteq V(G),
|S|=k}{\min} \pi(S)$\\
\cline{1-3}
\end{tabular}
\begin{center}
{Table~3.~Three kinds of tree-connectivities}
\end{center}

\begin{rem} There are many other kinds of generalizations
of the classical connectivity and edge-connectivity, such as the
restricted (edge-)connectivity in \cite{Esf} and
super edge-connectivity in \cite{LCXu}. However, our intention of
this survey is to only focus on the generalized (edge-)connectivity.
In very rare occasions, we will mention some results on the most
closely related concepts: pendant tree-connectivity,
path-connectivity and $k$-connectivity, in order to show the
differences among them.
\end{rem}

\subsection{Generalized connectivity and Steiner tree packing
problem}

The generalized edge-connectivity is related to two important
problems. For a given graph $G$ and $S\subseteq V(G)$, the problem
of finding a set of maximum number of edge-disjoint Steiner trees
connecting $S$ in $G$ is called the \emph{Steiner tree packing
problem}. The difference between the Steiner tree packing problem
and the generalized edge-connectivity is as follows: The Steiner
tree packing problem studies local properties of graphs since $S$ is
given beforehand, but the generalized edge-connectivity focuses on
global properties of graphs since it first needs to compute the
maximum number $\lambda(S)$ of edge-disjoint trees connecting $S$
and then $S$ runs over all $k$-subsets of $V(G)$ to get the minimum
value of $\lambda(S)$.

The problem for $S=V(G)$ is called the \emph{spanning tree packing
problem}. Note that spanning tree packing problem is a
specialization of Steiner tree packing problem (For $k=n$, each
Steiner tree connecting $S$ is a spanning tree of $G$). For any
graph $G$ of order $n$, the \emph{spanning tree packing number} or
\emph{$STP$ number}, is the maximum number of edge-disjoint spanning
trees contained in $G$. From the definitions of $\kappa_k(G)$ and
$\lambda_k(G)$, $\kappa_n(G)=\lambda_n(G)$ is exactly the spanning
tree packing number of $G$ (For $k=n$, both internally disjoint
Steiner trees connecting $S$ and edge-disjoint Steiner trees
connecting $S$ are edge-disjoint spanning trees). For the spanning
tree packing number, we refer to \cite{OY, Palmer}. Observe that
spanning tree packing problem is a special case of both the
generalized $k$-connectivity and the generalized
$k$-edge-connectivity. This problem has two practical applications.
One is to enhance the ability of fault tolerance \cite{Fragopoulou,
Itai}. Consider a source node $u$ that wants to broadcast a message
on a network with $\ell$ edge-disjoint spanning trees. The node $u$
copies $\ell$ messages to different spanning trees. If there are no
more than $\ell-1$ fault edges, all the other nodes can receive the
message. The other application is to develop efficient collective
communication algorithms in distributed memory parallel computers
\cite{Barden, Libeskind, Wang}. If the above source node has a large
number of data to transmit, we can let every edge-disjoint spanning
tree be responsible for only $1/\ell$ data to increase the
throughput. For any graph $G$, the maximum number of edge-disjoint
spanning trees in $G$ can be found in polynomial time; see
(\cite{Schrijver}, page 879). Actually, Roskind and Tarjan
\cite{Roskind} proposed an $O(m^2)$ time algorithm for finding
the maximum number of edge-disjoint spanning trees in an arbitrary
graph, where $m$ is the number of edges in the graph.

\subsection{Theoretical and application backgrounds of generalized connectivity}

In addition to being natural combinatorial measures, the generalized
connectivity and generalized edge-connectivity can be motivated by
their interesting interpretation in practice as well as theoretical
consideration.

From a theoretical perspective, both extremes of this problem are
fundamental theorems in combinatorics. One extreme of the problem is
when we have two terminals. In this case internally (edge-)disjoint
trees are just internally (edge-)disjoint paths between the two
terminals, and so the problem becomes the well-known Menger theorem.
The other extreme is when all the vertices are terminals. In this
case internally disjoint trees and edge-disjoint trees are just
edge-disjoint spanning trees of the graph, and so the problem becomes the
classical Nash-Williams-Tutte theorem.

\begin{thm}(Nash-Williams {\upshape \cite{Nash}}, Tutte {\upshape \cite{Tutte}})\label{th1-1}
A multigraph $G$ contains a system of $k$ edge-disjoint spanning
trees if and only if
$$
\|G/\mathscr{P}\|\geq k(|\mathscr{P}|-1)
$$
holds for every partition $\mathscr{P}$ of $V(G)$, where
$\|G/\mathscr{P}\|$ denotes the number of edges in $G$ between
distinct blocks of $\mathscr{P}$.
\end{thm}

The next theorem is due to Nash-Williams.

\begin{thm}{\upshape \cite{Nash2}}\label{th1-2}
Let $G$ be a graph. Then the edge set of $G$ can be covered by $t$
forests if and only if, for every nonempty subset $S$ of vertices of
$G$, $|E_G[S]|\leq t(|S|-1)$.
\end{thm}

The following corollary can be easily derived from Theorem
\ref{th1-1}.

\begin{cor}\label{cor1-3}
Every $2\ell$-edge-connected graph contains a system of $\ell$
edge-disjoint spanning trees.
\end{cor}

Kriesell \cite{Kriesell1} conjectured that this corollary can be
generalized for Steiner trees.

\begin{con}(Kriesell {\upshape \cite{Kriesell1}})\label{con1-3}
If a set $S$ of vertices of $G$ is $2k$-edge-connected (see later in Section
$1.5$ for the definition), then there is a set of $k$ edge-disjoint
Steiner trees in $G$.
\end{con}

Motivated by this conjecture, the Steiner Tree Packing Problem has
obtained wide attention and many results have been worked out, see
\cite{Kriesell1, Kriesell2, West, Jain, Lau}. In \cite{LMS} we set
up the relationship between the Steiner tree packing problem and the
generalized edge-connectivity.

The generalized edge-connectivity and the Steiner tree packing
problem have applications in $VLSI$ circuit design, see
\cite{Grotschel1, Grotschel2, Sherwani}. In this application, a
Steiner tree is needed to share an electronic signal by a set of
terminal nodes. Steiner tree is also used in computer communication
networks (see \cite{Du}) and optical wireless communication networks
(see \cite{Cheng}). Another application, which is our primary focus,
arises in the Internet Domain. Imagine that a given graph $G$
represents a network. We choose arbitrary $k$ vertices as nodes.
Suppose one of the nodes in $G$ is a \emph{broadcaster}, and all
other nodes are either \emph{users} or \emph{routers} (also called
\emph{switches}). The broadcaster wants to broadcast as many streams
of movies as possible, so that the users have the maximum number of
choices. Each stream of movie is broadcasted via a tree connecting
all the users and the broadcaster. So, in essence we need to find
the maximum number Steiner trees connecting all the users and the
broadcaster, namely, we want to get $\lambda (S)$, where $S$ is the
set of the $k$ nodes. Clearly, it is a Steiner tree packing problem.
Furthermore, if we want to know whether for any $k$ nodes the
network $G$ has above properties, then we need to compute
$\lambda_k(G)=\min\{\lambda (S)\}$ in order to prescribe the
reliability and the security of the network.

\subsection{Strength and generalized connectivity}

The \emph{strength} of a graph $G$ is defined as
$$
\eta(G)=\underset{X\subseteq E(G)}{\min} \
\frac{|X|}{\omega(G-X)-\omega(G)},
$$
where the minimum is taken over whenever the denominator is non-zero
and $\omega(G)$ denotes the number of components of $G$. From
Nash-Williams-Tutte theorem, a multigraph $G$ contains a system of
$k$ edge-disjoint spanning trees if and only if for any $X\subseteq
E(G)$, $|X|\geq k(\omega(G-X)-1)$. One can see that the concept of
the strength of a graph may be derived from the Nash-Williams-Tutte
theorem for connected graphs. By Nash-Williams-Tutte theorem,
$\kappa_n(G)=\lambda_n(G)=\lfloor \eta(G)\rfloor$ for a simple
connected graph $G$. For more details, we refer to \cite{Catlin,
Gusfield, Welsh}. In addition, the generalized (edge)-connectivity
and the strength of a graph can be used to measure the reliability
and the security of a network, see \cite{Cunningham, Matula}.

Similar to the strength of a graph, another interesting concept
involving the vertex set is the toughness of a graph. A graph $G$ is
\emph{$t$-tough} if $|S|\geq t\omega(G-S)$ for every subset $S$ of
the vertex set $V(G)$ with $\omega(G-S)>1$. The \emph{toughness} of
$G$, denoted by $\tau(G)$, is the maximum value of $t$ for which $G$
is $t$-tough (taking $\tau(K_n) = \infty$ for all $n>1$). Hence if
$G$ is not complete, then $\tau(G)=
\min\{\frac{|S|}{\omega(G-S)}\}$, where the minimum is taken over
all cut sets of vertices in $G$. Bauer, Broersma and Schmeichel had
a survey on this subject, see \cite{Bauer}.

\subsection{Notation and terminology}

All graphs considered in this paper are undirected, finite and
simple. We refer to book \cite{bondy} for graph theoretical notation
and terminology not described here. For a graph $G$, let $V(G)$,
$E(G)$, $e(G)$, $L(G)$, $\overline{G}$ and
$\alpha(G)$ denote the set of vertices, the set of
edges, the size or number of edges, the line graph, the complement and
the independence (or stable) number of $G$, respectively. As
usual, the {\it union}, denoted by $G\cup H$, of two graphs $G$ and
$H$ is the graph with vertex set $V(G)\cup V(H)$ and edge set
$E(G)\cup E(H)$. The disjoint union of $k$ copies of a same graph
$G$ is denoted by $k G$. The \emph{join} $G_1\vee G_2$ of $G_1$ and
$G_2$ is obtained from $G_1\cup G_2$ by joining each vertex of $G_1$
to every vertex of $G_2$. For $S\subseteq V(G)$, we denote by
$G\setminus S$ the subgraph obtained by deleting the vertices of $S$ together
with the edges incident with them from $G$. If $S=\{v\}$, we simply
write $G\setminus v$ for $G\setminus \{v\}$. If $S$ is a subset of
vertices of a graph $G$, the subgraph of $G$ induced by $S$ is
denoted by $G[S]$. If $M$ is an edge subset of $G$, then $G\setminus
M$ denote the subgraph by deleting the edges of $M$. The subgraph of
$G$ induced by $M$ is denoted by $G[M]$. If $M=\{e\}$, we simply
write $G\setminus e$ or $G-e$ for $G\setminus \{e\}$. We denote by
$E_G[X,Y]$ the set of edges of $G$ with one end in $X$ and the other
end in $Y$. If $X=\{x\}$, we simply write $E_G[x,Y]$ for
$E_G[\{x\},Y]$.

For two distinct vertices $x,y$ in $G$, let $\lambda(x,y)$ denote
the local edge-connectivity of $x$ and $y$. A subset $S\subseteq
V(G)$ is called \emph{$t$-edge-connected}, if $\lambda(x,y)\geq t$
for all $x\neq y$ in $S$. A $k$-connected graph $G$ is
\emph{minimally $k$-connected} if the graph $G-e$ is not
$k$-connected for any edge of $G$. A graph $G$ is \emph{$k$-regular}
if $d(v)=k$ for every $v\in V(G)$. A $3$-regular graph is called
\emph{cubic}. For $X=\{x_1,x_2,\cdots,x_k\}$ and
$Y=\{y_1,y_2,\cdots,y_k\}$, an \emph{$XY$-linkage} is defined as a
set of $k$ vertex-disjoint paths $x_iP_iy_i$ for every $i$ with
$1\leq i\leq k$. The \emph{Linkage Problem} is the problem of
deciding whether there exists an $XY$-linkage for given sets $X$ and
$Y$.

The \emph{Cartesian product} (also called the square product) of two
graphs $G$ and $H$, written as $G\square H$, is the graph with
vertex set $V(G)\times V(H)$, in which two vertices $(u,v)$ and
$(u',v')$ are adjacent if and only if $u=u'$ and $(v,v')\in E(H)$,
or $v=v'$ and $(u,u')\in E(H)$. Clearly, the Cartesian product is
commutative, that is, $G\Box H\cong H\Box G$. The lexicographic
product of two graphs $G$ and $H$, written as $G\circ H$, is defined
as follows: $V(G\circ H)=V(G)\times V(H)$, and two distinct vertices
$(u,v)$ and $(u',v')$ of $G\circ H$ are adjacent if and only if
either $(u,u')\in E(G)$ or $u=u'$ and $(v,v')\in E(H)$. Note that
unlike the Cartesian product, the lexicographic product is a
non-commutative product since $G\circ H$ is usually not isomorphic
to $H\circ G$.

A \emph{decision problem} is a question whose answer is either
¡°yes¡± or ¡°no¡±. Such a problem belongs to the class $\mathcal
{P}$ if there is a polynomial-time algorithm that solves any
instance of the problem in polynomial time. It belongs to the class
$\mathcal {N}\mathcal {P}$ if, given any instance of the problem
whose answer is ¡°yes¡±, there is a certificate validating this
fact, which can be checked in polynomial time; such a certificate is
said to be \emph{succinct}. It is immediate from these definitions
that $\mathcal {P}\subseteq \mathcal {N}\mathcal {P}$, inasmuch as a
polynomial-time algorithm constitutes, in itself, a succinct
certificate. A \emph{polynomial reduction} of a problem $P$ to a
problem $Q$ is a pair of polynomial-time algorithms, one of which
transforms each instance $I$ of $P$ to an instance $J$ of $Q$, and
the other of which transforms a solution for the instance $J$ to a
solution for the instance $I$. If such a reduction exists, we say
that $P$ is \emph{polynomially reducible} to $Q$, and write
$P\preceq Q$. A problem $P$ in $\mathcal {N}\mathcal {P}$ is
\emph{$\mathcal {N}\mathcal {P}$-complete} if $X\preceq P$ for
every problem $X$ in $\mathcal {N}\mathcal {P}$. The following two problems
are well-known $\mathcal {N}\mathcal {P}$-complete problems.

\noindent \textbf{$3$-DIMENSIONAL MATCHING($3$-DM):} Given three
sets $U$, $V$ and $W$ of equal cardinality, and a subset $T$ of
$U\times V\times W$, decide whether there is a subset $M$ of $T$
with $|M|=|U|$ such that whenever $(u,v,w)$ and $(u',v',w')$ are
distinct triples in $M$, then $u\neq u'$, $v\neq v'$, and $w\neq
w'$?

\noindent \textbf{BOOLEAN $3$-SATISFIABILITY ($3$-SAT):} Given a
boolean formula $\phi$ in conjunctive normal form with three
literals per clause, decide whether $\phi$ is satisfiable ?

\section{Results for some graph classes}

The following two observations are easily seen.
\begin{obs}\label{obs2-1}
If $G$ is a connected graph, then $\kappa_k(G)\leq \lambda_k(G)\leq
\delta(G)$.
\end{obs}
\begin{obs}\label{obs2-2}
If $H$ is a spanning subgraph of $G$, then $\kappa_k(H)\leq
\kappa_k(G)$ and $\lambda_k(H)\leq \lambda_k(G)$.
\end{obs}

\subsection{Results for complete graphs}

Chartrand, Okamoto and Zhang in \cite{Chartrand2} proved that if $G$
is the complete $3$-partite graph $K_{3,4,5}$, then $\kappa_3(G)=6$.
They also got the exact value of the generalized $k$-connectivity
for complete graph $K_n$.

\begin{thm}{\upshape \cite{Chartrand2}}\label{th2-3}
For every two integers $n$ and $k$ with $2\leq k\leq n$,
$$
\kappa_k(K_n)=n-\lceil k/2\rceil.
$$
\end{thm}

In \cite{LMS}, Li, Mao and Sun obtained the explicit value for
$\lambda_k(K_n)$. One may not expect that it is the same as
$\kappa_k(K_n)$.

\begin{thm} {\upshape \cite{LMS}}\label{th2-4}
For every two integers $n$ and $k$ with $2\leq k\leq n$,
$$
\lambda_k(K_n)=n-\lceil k/2\rceil.
$$
\end{thm}

From Theorems \ref{th2-3} and \ref{th2-4}, we get that
$\lambda_k(G)=\kappa_k(G)$ for a complete graph $G$. However, this
is a very special case. Actually, $\lambda_k(G)-\kappa_k(G)$ could
be very large. For example, let $G$ be a graph obtained from two
copies of the complete graph $K_n$ by identifying one vertex in each
of them. For $k\leq n$, $\lambda_k(G)=n-\lceil\frac{k}{2}\rceil$,
but $\kappa_k(G)=1$.

For the pendant $k$-tree-connectivity and $k$-path-connectivity of
the complete graph $K_n$, Hager obtained
that $\tau_k(K_n)=n-k$ in \cite{Hager}, and $\pi_k(K_n)=
\left\lfloor\frac{2n+k^2-3k}{2(k-1)}\right\rfloor$ in \cite{Hager2},
respectively. One can see that they are very different from the
generalized $k$-connectivity. Hager also gave the exact values of
the pendant $k$-tree-connectivity and $k$-path-connectivity for other
special graphs, such as complete bipartite graphs. For more results,
also see Mao \cite{Mao2}.

\subsection{Results for complete multipartite graphs}

Okamoto and Zhang \cite{Okamoto} investigated the generalized
$k$-connectivity of a regular complete bipartite graph $K_{a,a}$.
Naturally, one may ask whether we can compute the value of
generalized $k$-connectivity of a complete bipartite graph
$K_{a,b}$, or even a complete multipartite graph. For $k=n$, Peng,
Chen and Koh \cite{Peng}, and Peng and Tay \cite{Peng2} later,
obtained the $STP$ number of a complete multipartite graph.

\begin{thm}{\upshape \cite{Peng, Peng2}}\label{th2-5}
For a complete multipartite graph $G$, the $STP$ number of $G$ is
$$
\Big\lfloor\frac{e(G)}{|V(G)|-1}\Big\rfloor.
$$
\end{thm}

The above result means that for a complete multipartite graph $G$,
$\lambda_n(G)=\kappa_n(G)=\Big\lfloor\frac{e(G)}{|V(G)|-1}\Big\rfloor.$
Recently, Li, Li and Li \cite{LLL1, LLL2, WLi} devoted to solving
this problem for a general $k$. Restricting to simple graphs, they
rediscovered the result of Theorem \ref{th2-5} for complete
bipartite graphs and complete equipartition $3$-partite graphs. But,
it is worth to point out that their proof method, called the
\emph{List Method}, is more constructive, different from that of
Peng et al., and can exactly give all the
$\Big\lfloor\frac{e(G)}{|V(G)|-1}\Big\rfloor$ edge-disjoint spanning
trees.

Actually, Li, Li and Li used their \emph{List Method} and obtained
the value of generalized $k$-connectivity of all complete bipartite
graphs for $2\leq k\leq n$.

\begin{thm}{\upshape \cite{LLL1}}\label{th2-6}
Given any three positive integers $a,b,k$ such that $a\leq b$ and
$2\leq k\leq a+b$, let $K_{a,b}$ denote a complete bipartite graph
with a bipartition of sizes $a$ and $b$, respectively. Then we have
the following results:

if $k>b-a+2$ and $a-b+k$ is odd, then

$$\kappa_k(K_{a,b})=\frac{a+b-k+1}{2}+\Big\lfloor\frac{(a-b+k-1)(b-a+k-1)}
{4(k-1)}\Big\rfloor,$$

if $k>b-a+2$ and $a-b+k$ is even, then

$$\kappa_k(K_{a,b})=\frac{a+b-k}{2}+\Big\lfloor\frac{(a-b+k)(b-a+k)}
{4(k-1)}\Big\rfloor$$

and if $k\leq b-a+2$, then

$$\kappa_k(K_{a,b})=a$$
\end{thm}

It is not easy to obtain the exact value of generalized
$k$-connectivity of a complete multipartite graph. So they focused
on the complete equipartition $3$-partite graph and got the
following result.

\begin{thm}{\upshape \cite{LLL2}}\label{th2-7}
Given any positive integer $b\geq 2$, let $K_{b}^3$ denote a
complete equipartition $3$-partite graph in which every part
contains exactly $b$ vertices. Then we have the following results:
$$
\kappa_k(K_b^3)=\left\{
\begin{array}{ll}
\Big\lfloor\frac{\lceil k^2/3
\rceil+k^2-2kb}{2(k-1)}\Big\rfloor+3b-k, &if~k\geq
\frac{3b}{2};\\[7pt]
\Big\lfloor\frac{3bk+3b-k+1}{2k+1}\Big\rfloor, &if~\frac{3b}{4}<k<\frac{3b}{2}~and~k=1 \ (mod~3);\\[7pt]
\Big\lfloor\frac{3bk+6b-2k+1}{2k+2}\Big\rfloor, &if~b\leq k<\frac{3b}{2}~and~k=2 \ (mod~3);\\[7pt]
\big\lfloor\frac{3b}{2}\big\rfloor, &if~k<\frac{3b}{2}~and~k=0 \ (mod~3);\\[7pt]
\big\lfloor\frac{3b+1}{2}\big\rfloor, &otherwise.
\end{array}
\right.$$
\end{thm}

Let $U=\{u_1,u_2,\cdots,u_a\}$ and $V=\{v_1,v_2,\cdots,v_{b}\}$ be
the two parts of a complete bipartite graph $K_{a,b}$. Set
$S_i=\{u_1,u_2,\cdots,u_x,v_1,v_2,\cdots,v_{k-x}\}$ for $0\leq x\leq
k$. If $k>b-a+2$ and $a-b+k$ is odd, then
$\kappa_k(K_{a,b})=\kappa(S_{\frac{a-b+k-1}{2}})$, in the part $X$
there are $\frac{a-b+k-1}{2}$ vertices not in $S$, and in the part
$Y$ there are $\frac{a-b+k-1}{2}$ vertices not in $S$. The number of
vertices in each part but not in $S$ is almost the same. And if
$k>b-a+2$ and $a-b+k$ is even, then
$\kappa_k(K_{a,b})=\kappa(S_{\frac{a-b+k}{2}})$, in the part $X$
there are $\frac{a-b+k}{2}$ vertices not in $S$, and in the part $Y$
there are $\frac{a-b+k}{2}$ vertices not in $S$. The number of
vertices in each part but not in $S$ is the same.

Similarly, let $U=\{u_1,u_2,\cdots,u_b\}$,
$V=\{v_1,v_2,\cdots,v_{b}\}$ and $W=\{w_1,w_2,\cdots,w_{b}\}$ be the
three parts of a complete equipartition $3$-partite graph $K_b^3$.
Set
$S_{x,y,z}=\{u_1,u_2,\cdots,u_x,v_1,v_2,\\
\cdots,v_{y},v_1,v_2,\cdots,v_{z}\}$ for $0\leq x,y,z\leq k$ with
$x+y+z=k$. If $k=0 \ (mod~3)$, then
$\kappa_k(K_b^3)=\kappa(S_{\frac{k}{3},\frac{k}{3},\frac{k}{3}})$,
in the part $U$ there are $b-\frac{k}{3}$ vertices not in $S$, in
the part $V$ there are $b-\frac{k}{3}$ vertices not in $S$, and in
the part $W$ there are $b-\frac{k}{3}$ vertices not in $S$. The
number of vertices in each part but not in $S$ is the same. If $k=1
\ (mod~3)$, then
$\kappa_k(K_b^3)=\kappa(S_{\frac{k+2}{3},\frac{k-1}{3},\frac{k-1}{3}})$.
And if $k=2 \ (mod~3)$, then
$\kappa_k(K_b^3)=\kappa(S_{\frac{k+1}{3},\frac{k+1}{3},\frac{k-2}{3}})$.
In both cases, the number of vertices in each part but not in $S$ is
almost the same.

So, W. Li proposed the following two conjectures in her Ph.D. thesis
\cite{WLi}.

\begin{con}{\upshape \cite{WLi}}\label{con2-8}
For a complete equipartition $a$-partite graph $G$ with partition
$(X_1,X_2,\cdots,$ $X_a)$ and integer $k=ab+c$, where $b,c$ are
integers and $0\leq c\leq a-1$, we have $\kappa_k(G)=\kappa(S)$,
where $S$ is a $k$-subset of $V(G)$ such that $|S\cap
X_1|=\cdots=|S\cap X_c|=b+1$ and $|S\cap X_{c+1}|=\cdots=|S\cap
X_a|=b$.
\end{con}

\begin{con}{\upshape \cite{WLi}}\label{con2-9}
For a complete multipartite graph $G$, we have
$\kappa_k(G)=\kappa(S)$, where $S$ is a $k$-subset of $V(G)$ such
that the number of vertices in each part but not in $S$ are almost
the same.
\end{con}

\subsection{Results for Cayley graphs}

Let $X$ be a finite Abelian group, and its operation be called addition,
denoted by $+$. Let $A$ be a subset of $X\setminus \{0\}$ such that
$a\in A$ implies $-a\in A$, where $0$ is the identity element of $X$.
The \emph{Cayley graph} $Cay(X, A)$ is defined to have vertex set $X$
such that there is an edge between $x$ and $y$ if and only if $x-y\in A$.
A \emph{circulant graph} is a Cayley graph on a cyclic group. Observe
that $Cay(X,A)$ is connected if and only if $A$ is a generating set
of $X$. Cayley graphs are important objects of study in
algebraic graph theory; see \cite{Babai, Biggs}. In fact,
many mathematicians and computer scientists recommend Cayley graphs
as models for interconnection networks because they exhibit many
properties that ensure high performance; see \cite{AkersK,
Heydemann, Zhou}. In fact, a number of networks of both theoretical
and practical importance, including hypercubes, butterflies,
cube-connected cycles, star graphs and their generalizations, are
Cayley graphs. For the results pertaining to Cayley graphs as models
for interconnection networks, we refer to the survey papers
\cite{Heydemann, Lakshmivarahan}.

Because of the importance of Cayley graphs in network design and
the significance of reliability of networks, Sun and Zhou \cite{SZ}
studied the generalized connectivity of Cayley graphs.

\begin{thm}{\upshape \cite{SZ}}\label{th2-10}
Let $G$ be a cubic connected Cayley graph on an Abelian group with
order $n\geq 8$. Then
$$
\kappa_k(G)=\lambda_k(G)=\left\{
\begin{array}{ll}
2, &if~3\leq k\leq 6;\\[0.2cm]
1, &if~7\leq k\leq n.
\end{array}
\right.
$$
\end{thm}

\begin{thm}\label{th2-11}
Let $G$ be a connected Cayley graph of degree $4$ on an Abelian
group with order $n\geq 3$. Then
$$
\kappa_k(G)=\left\{
\begin{array}{ll}
3, &if~k=3;\\[0.2cm]
1~or~2, &if~8\leq k\leq n-2;\\[0.2cm]
2, &if~k=n-1,n.
\end{array}
\right.
$$
and
$$
\lambda_k(G)=\left\{
\begin{array}{ll}
3, &if~k=3;\\[0.2cm]
2~or~3, &if~4\leq k\leq 7;\\[0.2cm]
2, &if~8\leq k\leq n.
\end{array}
\right.
$$
\end{thm}

\section{Algorithm and complexity}

As it is well-known that, for any graph $G$, we have polynomial-time
algorithms to get the classical connectivity $\kappa(G)$ and the
edge-connectivity $\lambda(G)$, a natural question is whether there
is a polynomial-time algorithm to get the new parameters $\kappa_k(G)$
and $\lambda_k(G)$.

\subsection{Results for $\kappa_k$}

For a graph $G$, by the definition of $\kappa_3(G)$, it is natural
to study $\kappa(S)$ first, where $S$ is a $3$-subset of $V(G)$. A
question is then raised: for any fixed positive integer $\ell$,
given a $3$-subset $S$ of $V(G)$, is there a polynomial-time
algorithm to determine whether $\kappa(S)\geq \ell$ ? Li, Li and Zhou \cite{LLZ} gave a
positive answer by converting the problem into the \emph{$k$-Linkage
Problem} \cite{Robertson}. From this together with
$\kappa_3(G)=\min\{\kappa(S)\}$, the following theorem can be easily
obtained.

\begin{thm}{\upshape \cite{LLZ}}\label{th3-1}
Given a fixed positive integer $\ell \ (\ell\geq 2)$, for any graph
$G$ the problem of deciding whether $\kappa_3(G)\geq \ell$ can be
solved by a polynomial-time algorithm.
\end{thm}

The following two corollaries are immediate from the relation
$\kappa_3\leq \kappa\leq \delta$.

\begin{cor}{\upshape \cite{LLZ}}\label{cor3-2}
Given a fixed positive integer $\kappa$, for any graph $G$ with
connectivity $\kappa$, the problem of deciding $\kappa_3(G)$ can be
solved by a polynomial-time algorithm.
\end{cor}

\begin{cor}{\upshape \cite{LLZ}}\label{cor3-3}
Given a fixed positive integer $\delta$, for any graph $G$ with
minimum degree $\delta$, the problem of deciding $\kappa_3(G)$ can
be solved by a polynomial-time algorithm.
\end{cor}

Furthermore, for a planar graph they derived the following result.

\begin{pro}{\upshape \cite{LLZ}}\label{pro3-4}
For a planar graph $G$ with connectivity $\kappa(G)$, the problem of
determining $\kappa_3(G)$ has a polynomial-time algorithm and its
complexity is bounded by $O(n^8)$.
\end{pro}

They mentioned that the above complexity is not very good, and so
the problem of finding a more efficient algorithm is interesting.
The complexity of the problem of determining $\kappa_3(G)$ for a
general graph is not known: Can it be solved in polynomial time or
$\mathcal {N}\mathcal {P}$-hard ? Nevertheless, they derived a
polynomial-time algorithm to determine it approximately with a
constant ratio.

\begin{pro}{\upshape \cite{LLZ}}\label{pro3-5}
The problem of determining $\kappa_3(G)$ for any graph $G$ can be
solved by a polynomial-time approximation algorithm with a constant
ratio about $\frac{3}{4}$.
\end{pro}

Later, Li and Li \cite{LL} considered to generalize the result of
Theorem \ref{th3-1} to that for general $k$ and obtained the
following theorem.

\begin{thm}{\upshape \cite{LL}}\label{th3-6}
Given two fixed positive integers $k$ and $\ell$, for any graph $G$
the problem of deciding whether $\kappa_{k}(G)\geq \ell$ can be
solved by a polynomial-time algorithm.
\end{thm}

For $k$ a fixed integer but $\ell$ an arbitrary integer, Li and Li
proposed the following problem.

\begin{prob}{\upshape \cite{LL}}\label{prob3-7}
Given a graph $G$, a $4$-subset $S$ of $V(G)$ and an integer $\ell \
(\ell\geq 2)$, decide whether there are $\ell$ internally disjoint
trees connecting $S$, namely decide whether $\kappa(S)\geq \ell$?
\end{prob}

At first, they proved that Problem \ref{prob3-7} is $\mathcal
{N}\mathcal {P}$-complete by reducing $3$-DM to it. Next, they
showed that for a fixed $k\geq 5$, in Problem \ref{prob3-7}
replacing the $4$-subset of $V(G)$ with a $k$-subset of $V(G)$, the
problem is still $\mathcal {N}\mathcal {P}$-complete, which can be
proved by reducing Problem \ref{prob3-7} to it. Thus, they obtained
the following result.

\begin{pro}{\upshape \cite{LL}}\label{pro3-8}
For any fixed integer $k\geq 4$, given a graph $G$, a $k$-subset $S$
of $V(G)$ and an integer $\ell \ (\ell\geq 2)$, deciding whether
there are $\ell$ internally disjoint trees connecting $S$, namely
deciding whether $\kappa(S)\geq \ell$, is $\mathcal {N}\mathcal
{P}$-complete.
\end{pro}

As shown in above proposition, Li and Li \cite{LL} only showed that
for any fixed integer $k\geq 4$, deciding whether $\kappa(S)\geq
\ell$ is $\mathcal {N}\mathcal {P}$-complete. For $k=3$, the
complexity is yet not known. So, S. Li in her Ph.D. thesis
\cite{SLi} conjectured that it is $\mathcal {N}\mathcal
{P}$-complete.

\begin{con}{\upshape \cite{SLi}}\label{con3-9}
Given a graph $G$ and a $3$-subset $S$ of $V(G)$ and an integer
$\ell \ (\ell\geq 2)$, deciding whether there are $\ell$ internally
disjoint trees connecting $S$, namely deciding whether
$\kappa(S)\geq \ell$, is $\mathcal {N}\mathcal {P}$-complete.
\end{con}

Recently, Chen, Li, Liu and Mao \cite{CLLM} confirmed
the conjecture. In their proof, they employed the following new
$\mathcal {N}\mathcal {P}$-complete problem.

\begin{prob}{\upshape \cite{CLLM}} \label{pro3-10}
Given a tripartite graph $G=(V,E)$ with three partitions
$(\overline{U},\overline{V},\overline{W})$, and
$|\overline{U}|=|\overline{V}|=|\overline{W}|=q$, decide whether
there is a partition of $V$ into $q$ disjoint $3$-sets
$V_1,V_2,\ldots, V_q$ such that every
$V_i=\{v_{i_1},v_{i_2},v_{i_3}\}$ satisfies that $v_{i_1}\in
\overline{U}$, $v_{i_2}\in \overline{V}$, $v_{i_3}\in \overline{W}
$, and $G[V_i]$ is connected ?
\end{prob}

By reducing the $3$-$DM$ to Problem \ref{pro3-10}, they proved that
Problem \ref{pro3-10} is $\mathcal {N}\mathcal {P}$-complete.
Furthermore, they confirmed that Conjecture \ref{con3-9} is true by
reducing Problem \ref{pro3-10} to it.

\begin{pro}{\upshape \cite{CLLM}}\label{pro3-11} Given a graph
$G$, a $3$-subset $S$ of $V(G)$ and an integer $\ell \ (\ell\geq
2)$, the problem of deciding whether $G$ contains $\ell$ internally
disjoint trees connecting $S$ is $\mathcal {N}\mathcal
{P}$-complete.
\end{pro}

From Propositions \ref{pro3-8} and \ref{pro3-11}, we conclude that
if $k \ (k\geq 3)$ is a fixed integer and $\ell \ (\ell\geq 2)$ is
an arbitrary positive integer, the problem of deciding whether
$\kappa(S)\geq \ell$ is $\mathcal {N}\mathcal {P}$-complete. S. Li
in her Ph.D. thesis conjectured that the problem of deciding whether
$\kappa_{k}(G)\geq \ell$ is also $\mathcal {N}\mathcal
{P}$-complete.

\begin{con} {\upshape \cite{SLi}}\label{con3-12}
For a fixed integer $k\geq 3$, given a graph $G$ and an integer
$\ell \ (\ell\geq 2)$, the problem of deciding whether
$\kappa_{k}(G)\geq \ell$, is $\mathcal {N}\mathcal {P}$-complete.
\end{con}

The above conjecture is still open. Li and Li turned to considering
the case that $\ell$ is a fixed integer but $k$ is an arbitrary integer,
and they employed another problem.

\begin{prob}{\upshape \cite{LL}}\label{prob3-14}
Given a graph $G$, a subset $S$ of $V(G)$, decide whether there are
two internally disjoint trees connecting $S$, namely decide whether
$\kappa(S)\geq 2$?
\end{prob}

By reducing the $3$-SAT to Problem \ref{prob3-14}, they also
verified that Problem \ref{prob3-14} is $\mathcal {N}\mathcal
{P}$-complete. Next they showed that for a fixed integer $\ell\geq
3$, similar to Problem \ref{prob3-14} if we want to decide whether
there are $\ell$ internally disjoint trees connecting $S$ rather
than two, the problem is still $\mathcal {N}\mathcal {P}$-complete,
which can be easily proved by reducing Problem \ref{prob3-14} to it.
Then they got the following theorem.

\begin{thm}{\upshape \cite{LL}}\label{th3-15}
For any fixed integer $\ell\geq 2$, given a graph $G$ and a subset
$S$ of $V(G)$, deciding whether there are $\ell$ internally disjoint
trees connecting $S$, namely deciding whether $\kappa(S)\geq \ell$,
is $\mathcal {N}\mathcal {P}$-complete.
\end{thm}

\subsection{Results for $\lambda_k$}

For the computational complexity of the generalized
edge-connectivity $\lambda_k(G)$, Chen, Li, Liu and Mao
in the same paper \cite{CLLM} got the following result.

\begin{thm}{\upshape \cite{CLLM}} \label{th3-16} Given two
fixed positive integers $k$ and $\ell$, for any graph $G$ the
problem of deciding whether $\lambda_k(G)\geq \ell$ can be solved by
a polynomial-time algorithm.
\end{thm}

If $k$ or $\ell$ is/are not fixed, the problem for the
computational complexity of the generalized edge-connectivity
$\lambda_k(G)$ is still not known. To conclude this chapter,
we propose the following conjectures.

\begin{con}
For any fixed integer $k\geq 3$, given a graph $G$, a
$k$-subset $S$ of $V(G)$, and an integer $ \ell \ (2\leq \ell\leq
n-2)$, deciding whether there are $\ell$ edge-disjoint Steiner trees
connecting $S$, namely deciding whether $\lambda(S)\geq \ell$, is
$\mathcal {N}\mathcal {P}$-complete.
\end{con}

\begin{con}
For a fixed integer $k\geq 3$, given a graph $G$ and an
integer $ \ell \ (2\leq \ell\leq n-2)$, the problem of deciding
whether $\lambda_k(G)\geq \ell$ is $\mathcal {N}\mathcal
{P}$-complete.
\end{con}

\begin{con}
For any fixed integer $\ell\geq 2$, given a graph $G$, a
subset $S$ of $V(G)$, deciding whether there are $\ell$
edge-disjoint Steiner trees connecting $S$, namely deciding whether
$\lambda(S)\geq \ell$, is $\mathcal {N}\mathcal {P}$-complete.
\end{con}

\section{Sharp bounds for the generalized connectivity}

From the last section we know that it is almost impossible to get
the exact value of the generalized (edge-)connectivity for a given
arbitrary graph. So people tried to give some nice bounds for it,
especially sharp upper and lower bounds.

\subsection{Bounds for $\kappa_k(G)$ and $\lambda_k(G)$}

From Theorems \ref{th2-3} and \ref{th2-4}, i.e.,
$\kappa_k(K_n)=n-\lceil k/2 \rceil$ and $\lambda_(K_n)=n-\lceil k/2
\rceil$, one can see the following two consequences since any
connected graph $G$ is a subgraph of a complete graph.

\begin{pro}{\upshape \cite{LMS}}\label{pro4-1}
Let $k,n$ be two integers with $2\leq k\leq n$. For a connected
graph $G$ of order $n$, we have $1\leq \kappa_k(G)\leq n-\lceil k/2 \rceil$.
Moreover, the upper and lower bounds are sharp.
\end{pro}

\begin{pro}{\upshape \cite{LMS}}\label{pro4-2}
Let $k,n$ be two integers with $2\leq k\leq n$. For a connected
graph $G$ of order $n$, we have $1\leq \lambda_k(G)\leq n-\lceil k/2
\rceil$. Moreover, the upper and lower bounds are sharp.
\end{pro}
For the above two propositions, one can easily check that the
complete graph $K_n$ attains the upper bound and any tree $T_{n}$
attains the lower bound.

People mainly focus on sharp upper and lower bounds of $\kappa_k(G)$
and $\lambda_k(G)$ in terms of $\kappa$ and $\lambda$, respectively.
Li and Mao \cite{LM} derived a lower bound from Corollary
\ref{cor1-3}.

\begin{pro}{\upshape \cite{LM}}\label{pro4-3}
For a connected graph $G$ of order $n$ and and an integer $k$ with
$3\leq k\leq n$, we have $\lambda_k(G)\geq
\lfloor\frac{1}{2}\lambda(G)\rfloor$. Moreover, the lower bound is
sharp.
\end{pro}

In order to show the sharpness of this lower bound for $k=n$, they
showed that the Harary graph $H_{n,2r}$ attains this bound. For a
general $k \ (3\leq k\leq n)$, one can check that the cycle $C_n$
can attain the lower bound since
$\frac{1}{2}\lambda(C_n)=1=\lambda_k(C_n)$.

It seems difficult to get the sharp lower bound of $\kappa_k(G)$.
So, Li, Li and Zhou focused on the case $k=3$. By their method,
called the \emph{Path-Bundle Transformation} method, they obtained
the following result.

\begin{thm}{\upshape \cite{LLZ}}\label{th4-4}
Let $G$ be a connected graph with $n$ vertices. For every two
integers $s$ and $r$ with $s\geq 0$ and $r\in \{0,1,2,3\}$, if
$\kappa(G)=4s+r$, then $\kappa_3(G)\geq 3s+\lceil\frac{r}{2}\rceil$.
Moreover, the lower bound is sharp. We simply write $\kappa_3(G)\geq
\frac{3\kappa-2}{4}$.
\end{thm}

To show that the lower bound of Theorem \ref{th4-4} is sharp, they
gave the following example.

\noindent \textbf{Example 4.1}. For $\kappa(G)=4k+2i$ with $i=0$ or
$1$, they constructed a graph $G$ as follows (see Figure 4.1 $(a)$
): Let $Q=Y_1\cup Y_2$ be a vertex cut of $G$, where $Q$ is a clique
and $|Y_1|=|Y_2|=2k+i$, $G-Q$ has $2$ components $C_1, C_2$.
$C_1=\{v_3\}$ and $v_3$ is adjacent to every vertex in $Q$;
$C_2=\{v_1\}\cup \{v_2\}\cup X$, $|X|=2k+i$, the subgraph induced by
$X$ is an empty graph, each vertex in $X$ is adjacent to every
vertex in $Q\cup \{v_1,v_2\}$, $v_i$ is adjacent to every vertex
$Y_i$ for $i=1,2$. It can be checked that $\kappa(G)=4k+2i$ and
$\kappa_3(G)=3k+i$, which means that $G$ attains the lower bound.

For $\kappa(G)=4k+2i+1$ with $i=0$ or $1$, they constructed a graph
$G$ as follows (see Figure 4.1 $(b)$ ): Let $Q=Y_1\cup Y_2\cup
\{y_0\}$ be a vertex cut of $G$, where $Q$ is a clique and
$|Y_1|=|Y_2|=2k+i$. $G-Q$ has $2$ components $C_1, C_2$.
$C_1=\{v_3\}$ and $v_3$ is adjacent to every vertex in $Q$;
$C_2=\{v_1\}\cup \{v_2\}\cup X$, $|X|=2k+i$, the subgraph induced by
$X$ is an empty graph, each vertex in $X$ is adjacent to every
vertex in $Q\cup \{v_1,v_2\}$, $v_i$ is adjacent to every vertex
$Y_i$ for $i=1,2$, and both $v_1$ and $v_2$ are adjacent to $y_0$.
It can be checked that $\kappa(G)=4k+2i+1$ and $\kappa_3(G)=3k+i+1$,
which means that $G$ attains the lower bound.

\begin{figure}[!hbpt]
\begin{center}
\includegraphics[scale=0.8]{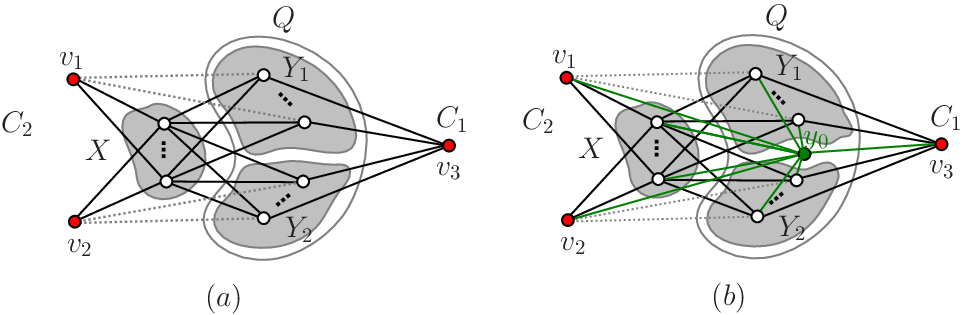}
\end{center}
\begin{center}
\caption{$(a)$ For $\kappa(G)=4k+2i$ with $i=0$, the graph attaining
the lower bound of Theorem \ref{th4-4}. $(b)$ For $\kappa(G)=4k+2i$
with $i=1$, the graph attaining the lower bound of Theorem
\ref{th4-4}.}
\end{center}\label{fig1}
\end{figure}

Kriesell \cite{Kriesell1} obtained a result on the Steiner tree
packing problem: Let $t\geq 1$ be a natural number and $G$ a graph,
and let $\{a,b,c\}\subseteq V(G)$ be
$\lfloor\frac{8t+3}{6}\rfloor$-edge-connected in $G$. Then there
exists a system of $t$ edge-disjoint $\{a,b,c\}$-spanning trees.
Using his result, Li, Mao and Sun derived a sharp lower bound of
$\lambda_3(G)$ and gave graphs attaining the bound. With this lower
bound, they got some results for line graphs (see Section $7$) and
planar graphs.

\begin{pro}{\upshape \cite{LMS}}\label{pro4-5}
Let $G$ be a connected graph with $n$ vertices. For every two
integers $s$ and $r$ with $s\geq 0$ and $r\in \{0,1,2,3\}$, if
$\lambda(G)=4s+r$, then $\lambda_3(G)\geq
3s+\lceil\frac{r}{2}\rceil$. Moreover, the lower bound is sharp. We
simply write $\lambda_3(G)\geq \frac{3\lambda-2}{4}$.
\end{pro}

They gave the following graph class to show that the lower bound is
sharp.

\noindent \textbf{Example 4.2}. For $\lambda=4s$ with $s\geq 1$, let
$P=X_1\cup X_2$ and $Q=Y_1\cup Y_2$ be two cliques with
$|X_1|=|Y_1|=2s$ and $|X_2|=|Y_2|=2s$. Let $x,y$ be adjacent to
every vertex in $P, Q$, respectively, and $z$ be adjacent to every
vertex in $X_1$ and $Y_1$. Finally, they finished the construction
of $G$ by adding a perfect matching between $X_2$ and $Y_2$. It can
be checked that $\lambda=4s$ and $\lambda(S)\geq 3s$. One can also
check that for other three vertices of $G$ the number of
edge-disjoint trees connecting them is not less than $3s$. So,
$\lambda_3(G)=3s$ and the graph $G$ attains the lower bound.

\begin{figure}[!hbpt]
\begin{center}
\includegraphics[scale=0.8]{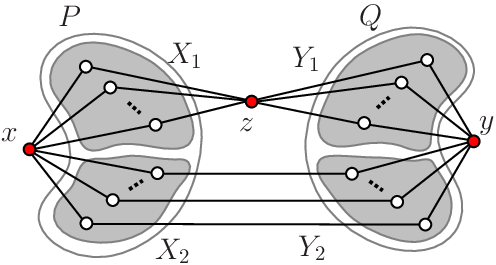}
\end{center}
\begin{center}
\caption{The graph with $\lambda(G)=4s$ and $\lambda_3(G)=3s$ that
attains the lower bound of Proposition \ref{pro4-5}.}
\end{center}\label{fig2}
\end{figure}

For $\lambda=4s+1$, let $|X_1|=|Y_1|=2s+1$ and $|X_2|=|Y_2|=2s$; for
$\lambda=4s+2$, let $|X_1|=|Y_1|=2s+1$ and $|X_2|=|Y_2|=2s+1$; for
$\lambda=4s+3$, let $|X_1|=|Y_1|=2s+2$ and $|X_2|=|Y_2|=2s+1$, where
$s\geq 1$. Similarly, one can check that $\lambda_3(G)=3s+1$ for
$\lambda=4s+1$; $\lambda_3(G)=3s+1$ for $\lambda=4s+2$;
$\lambda_3(G)=3s+2$ for $\lambda=4s+3$.

For the case $s=0$, $G=P_n$ satisfies that
$\lambda(G)=\lambda_3(G)=1$; $G=C_n$ satisfies that $\lambda(G)=2$
and $\lambda_3(G)=1$; $G=H_t$ satisfies that $\lambda(G)=3$ and
$\lambda_3(G)=2$, where $H_t$ denotes the graph obtained from $t$
copies of $K_4$ by identifying a vertex from each of them in the way
shown in Figure $4.3$.

\begin{figure}[!hbpt]
\begin{center}
\includegraphics[scale=0.8]{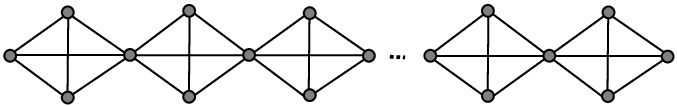}
\end{center}
\begin{center}
\caption{The graph $H_t$ with $\lambda(H_t)=3,\lambda_3(H_t)=2$.}
\end{center}\label{fig3}
\end{figure}

Li, Mao and Sun gave a sharp upper bound of $\lambda_k(G)$.

\begin{pro} {\upshape \cite{LMS}}\label{pro4-6}
For any graph $G$ of order $n$, $\lambda_k(G)\leq \lambda(G)$.
Moreover, the upper bound is sharp.
\end{pro}

But, for $\kappa_k(G)$, Li \cite{SLi} only proved that
$\kappa_k(G)\leq \kappa(G)$ for $3\leq k\leq 6$.

\begin{thm} {\upshape \cite{SLi}}\label{th4-7}
Let $G$ be a connected graph of order $n\geq 6$. Then for $3\leq
k\leq 6$, $\kappa_k(G)\leq \kappa(G)$. Moreover, the upper bound is
always sharp for $3\leq k\leq 6$.
\end{thm}

A natural question is why $\kappa_k(G)\leq \kappa(G)$ is not true
for $k\geq 6$ ? One may want to solve this problem by proving
$\kappa_k(G)\leq \kappa_{k-1}(G)$ for $3\leq k\leq n$, namely,
considering whether $\kappa_k$ is monotonically decreasing in $k$.
Unfortunately, Li found a counterexample $G$ such that
$\kappa_4(G)\geq \kappa_3(G)$. See the graph $G$ shown in Figure
$4.1$ $(a)$ for $i=1$. Li showed that $\kappa(G)=4k+2$ and
$\kappa_{3}(G)=3k+1$. It can be checked that the generalized
$4$-connectivity $\kappa_4(G)=3k+2$, which means that
$\kappa_4(G)\geq \kappa_3(G)$ for the graph $G$.

She also gave a graph $H(k,t)=(K_{\frac{k}{2}}\cup
K_{\frac{k}{2}})\vee K_t$, where $k \geq 6$ and $t \geq 1$, to show
that the monotone property of $\kappa_k$,
namely, $\kappa_n\leq \kappa_{n-1}\leq \cdots \kappa_4\leq
\kappa_3\leq \kappa$, is not true for $2\leq k\leq n$.

\begin{pro} {\upshape \cite{SLi}}\label{th4-8}
For any two integer $k \geq 6$ and $t \geq 1$,
$\kappa_{k+2}(H(k+1,t))\geq \kappa_{k+1}(H(k,t))$.
\end{pro}

However, for cubic graphs the conclusion holds.

\begin{thm} {\upshape \cite{SLi}}\label{th4-9}
If $G$ is a cubic graph, then $\kappa_k(G)\leq \kappa_{k-1}(G)$ for
$3\leq k\leq n$.
\end{thm}

Li and Mao \cite{LM} showed that the monotone property of
$\lambda_k$ is true for $2\leq k\leq n$ although it is not true for
$\kappa_k$.

\begin{pro}{\upshape \cite{LM}}\label{pro4-10}
For two integers $k$ and $n$ with $2\leq k\leq n-1$, and a connected
graph $G$, we have $\lambda_{k+1}(G)\leq \lambda_{k}(G)$.
\end{pro}

From Observation \ref{obs2-1}, we know that $\kappa_k(G)\leq
\lambda_k(G) \leq \delta$. Actually, Li, Mao and Sun \cite{LMS}
showed that the graph $G=K_k\vee (n-k)K_1 \ (n\geq 3k)$ satisfies
that $\kappa_k(G)=\lambda_k(G)=\kappa(G)=\lambda(G)=\delta(G)=k$,
which implies that the upper bounds of Observation \ref{obs2-1},
Proposition \ref{pro4-6} and Theorem \ref{th4-7} are sharp.

Li and Mao \cite{LM4} gave a sufficient condition for
$\lambda_k(G)\leq \delta-1$. Li \cite{SLi} obtained similar results
on the generalized $k$-connectivity.

\begin{pro} {\upshape \cite{LM4}}\label{pro4-11}
Let $G$ be a connected graph of order $n$ with minimum degree
$\delta$. If there are two adjacent vertices of degree $\delta$,
then $\lambda_k(G)\leq \delta-1$ for $3\leq k\leq n$. Moreover, the
upper bound is sharp.
\end{pro}

\begin{pro} {\upshape \cite{SLi}}\label{pro4-12}
Let $G$ be a connected graph of order $n$ with minimum degree
$\delta$. If there are two adjacent vertices of degree $\delta$,
then $\kappa_k(G)\leq \delta-1$ for $3\leq k\leq n$. Moreover, the
upper bound is sharp.
\end{pro}

With the above bounds, we will focus on their applications. From
Theorems \ref{th4-4} and \ref{th4-7}, Li, Li and Zhou derived
sharp bounds for planar graphs.

\begin{thm} {\upshape \cite{LLZ}}\label{th4-13}
If $G$ is a connected planar graph, then $\kappa(G)-1\leq
\kappa_3(G)\leq \kappa(G)$.
\end{thm}

Motivated by constructing graphs to show that the upper and lower
bounds are sharp, they obtained some lemmas. By the well-known
Kuratowski's theorem \cite{bondy}, they verified the following
lemma.

\begin{lem} {\upshape \cite{LLZ}}\label{lem4-14}
For a connected planar graph $G$ with $\kappa_3(G)=k$, there are no
$3$ vertices of degree $k$ in $G$, where $k\geq 3$.
\end{lem}

They also studied the generalized $3$-connectivity of four kinds of
graphs.

\begin{lem} {\upshape \cite{LLZ}}\label{lem4-15}
If $\kappa(G)\geq 3$, then $\kappa_3(G-e)\geq 2$ for any edge $e\in
E(G)$.
\end{lem}

\begin{lem} {\upshape \cite{LLZ}}\label{lem4-16}
If $G$ is a planar minimally $3$-connected graph, then
$\kappa_3(G)=2$.
\end{lem}

\begin{lem} {\upshape \cite{LLZ}}\label{lem4-17}
Let $G$ be a $4$-connected graph and let $H$ be a graph obtained
from $G$ by adding a new vertex $w$ and jointing it to $3$ vertices
of $G$. Then $\kappa_3(H)=\kappa(H)=3$.
\end{lem}

\begin{lem} {\upshape \cite{LLZ}}\label{lem4-18}
If $G$ is a planar minimally $4$-connected graph, then
$\kappa_3(H)=3$.
\end{lem}

If $G$ is a connected planar graph, then $1\leq \kappa(G)\leq 5$ by
Theorem \ref{th4-13}. Then, for each $1\leq \kappa(G)\leq 5$, they
gave some classes of planar graphs attaining the bounds of
$\kappa_3(G)$, respectively.

\noindent \emph{Case 1}: $\kappa(G)=1$. For any graph $G$ with
$\kappa(G)=1$, obviously $\kappa_3(G)\geq 1$ and so
$\kappa_3(G)=\kappa(G)=1$. Therefore, all planar graphs with
connectivity $1$ can attain the upper bound, but can not attain the
lower bound.

\noindent \emph{Case 2}: $\kappa(G)=2$. Let $G$ be a planar graph
with $\kappa(G)=2$ and having two adjacent vertices of degree $2$.
Then by Theorem \ref{th4-13} $\kappa_3(G)\leq 1$ and so
$\kappa_3(G)=1=\kappa(G)-1$. Therefore, this class of graphs attain
the lower bound.

Let $G$ be a planar minimally $3$-connected graph. By the
definition, for any edge $e\in E(G)$, $\kappa_3(G-e)=2$. Then by
Lemma \ref{lem4-15}, it follows that $\kappa_3(G-e)=2$. Therefore,
the $2$-connected planar graph $G-e$ attains the upper bound.

\noindent \emph{Case 3}: $\kappa(G)=3$. For any planar minimally
$3$-connected graph $G$, we know that $\kappa(G)=3$ and by Lemma
\ref{lem4-16}, $\kappa_3(G)=2=\kappa(G)-1$. So this class of graphs
attain the lower bound.

Let $G$ be a planar $4$-connected graph and let $H$ be a graph
obtained from $G$ by adding a new vertex $w$ in the interior of a
face for some planar embedding of $G$ and joining it to $3$ vertices
on the boundary of the face. Then $H$ is still planar and by Lemma
\ref{lem4-17}, one can immediately get that
$\kappa_3(H)=\kappa(H)=3$, which means that $H$ attains the upper
bound.

\begin{figure}[!hbpt]
\begin{center}
\includegraphics[scale=0.7]{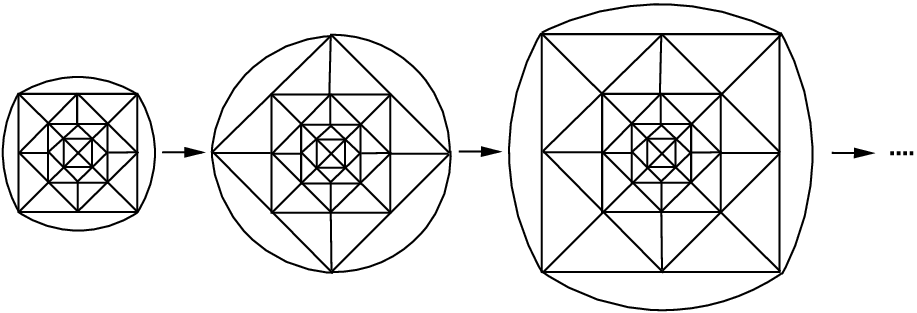}
\end{center}
\begin{center}
\caption{The graphs for the upper bound of Case $4$.}
\end{center}\label{fig4}
\end{figure}

\noindent \emph{Case 4}: $\kappa(G)=4$. For any planar minimally
$4$-connected graph $G$, one knows that $\kappa(G)=4$, and by Lemma
\ref{lem4-18}, $\kappa_3(G)=3=\kappa(G)-1$. So this class of graphs
attain the lower bound.

For every graph in Figure $4.4$, the vertex in the center has degree
$4$ and it can be checked that for any $2$ vertices there always
exist four pairwise internally disjoint paths connecting them, which
means that $\kappa(G)=4$. It can also be checked that for any $3$
vertices there always exist four pairwise internally disjoint trees
connecting them. Combining this with Theorem \ref{th4-13}, one can
get that $\kappa_3(G)=4$. Therefore, the graphs attain the upper
bound. Moreover, we can construct a series of graphs according to
the pattern of Figure $4.4$, which attain the upper bound.

\emph{Case 5}: $\kappa(G)=5$. For any planar graph $G$ with
$\kappa(G)=5$, by Lemma \ref{lem4-14} one can get that
$\kappa_3(G)=4$. So, any planar graph $G$ with connectivity $5$ can
attain the lower bound, but obviously can not attain the upper
bound.

Similarly, the following result is obvious from Propositions
\ref{pro4-5} and \ref{pro4-6}.
\begin{pro}\label{pro4-19}{\upshape \cite{LMS}}
If $G$ be a connected planar graph, then $\lambda(G)-1\leq
\lambda_3(G)\leq \lambda(G)$.
\end{pro}

\section{Graphs with given generalized (edge-)connectivity}

From the last section, we know that $1\leq \kappa_k(G)\leq n-\lceil
k/2 \rceil$ and $1\leq \lambda_k(G)\leq n-\lceil k/2 \rceil$ for a
connected graph $G$. Li, Mao and Sun \cite{LMS} considered to
characterize graphs attaining the upper bounds, namely, graphs with
$\kappa_k(G)=n-\lceil k/2 \rceil$ or $\lambda_k(G)=n-\lceil k/2
\rceil$. Since a complete graph $K_n$ possesses the maximum
generalized (edge-)connectivity, they wanted to find out the
critical value of the number of edges, denoted by $\mathbbm{t}$,
such that the generalized (edge-)connectivity of the resulting
graph will keep being $n-\lceil k/2 \rceil$ by deleting
$\mathbbm{t}$ edges from a complete graph $K_n$ but will not keep
being $n-\lceil k/2 \rceil$ by deleting $\mathbbm{t}+1$ edges. By
further investigation, they conjectured that $\mathbbm{t}$ may be
$0$ for $k$ even and $\frac{k-1}{2}$ for $k$ odd.

First, they noticed that for arbitrary $S\subseteq V(G)$ there are
two types of edge-disjoint trees connecting $S$: A tree of Type
$I$ is a tree whose edges all belong to $E(G[S])$; a tree of
Type $II$ is a tree containing at least one edge of
$E_G[S,\bar{S}]$. We denote the set of the edge-disjoint trees of
Type $I$ and Type $II$ by $\mathscr{T}_1$ and $\mathscr{T}_2$,
respectively. Let $\mathscr{T}=\mathscr{T}_1\cup \mathscr{T}_2$.

\begin{lem}{\upshape \cite{LMS}}\label{lem5-1}
Let $S\subseteq V(G)$, $|S|=k$ and $T$ be a tree connecting $S$. If
$T\in \mathscr{T}_1$, then $T$ uses $k-1$ edges of $E(G[S])\cup
E_G[S,\bar{S}]$; If $T\in \mathscr{T}_2$, then $T$ uses at least $k$
edges of $E(G[S])\cup E_G[S,\bar{S}]$.
\end{lem}

\subsection{Graphs with $\kappa_k(G)=n-\lceil\frac{k}{2}\rceil$ and
$\lambda_k(G)=n-\lceil\frac{k}{2}\rceil$}

They found that $|E(G[S])\cup E_G[S,\bar{S}]|$ is fixed once the
graph $G$ is given whatever there exist how many trees of Type $I$
and how many trees of Type $II$. From Lemma \ref{lem5-1}, each tree
will use certain number of edges in $E(G[S])\cup E_G[S,\bar{S}]$.
Deleting excessive edges from a complete graph $K_n$ will result in
that the remaining edges in $E(G[S])\cup E_G[S,\bar{S}]$ will not
form $n-\lceil k/2 \rceil$ trees. By using such an idea, they proved
that $\lambda_k(G)<n-\lceil\frac{k}{2}\rceil$ for $\mathbbm{t}\geq
1$ ($k$ is even) and $\lambda_k(G)<n-\lceil\frac{k}{2}\rceil$ for
$\mathbbm{t}\geq \frac{k+1}{2}$ ($k$ is odd). Furthermore, from
Observation \ref{obs2-1}, $\kappa_k(G)<n-\lceil\frac{k}{2}\rceil$
for $\mathbbm{t}\geq 1$ ($k$ is even) and
$\kappa_k(G)<n-\lceil\frac{k}{2}\rceil$ for $\mathbbm{t}\geq
\frac{k+1}{2}$ ($k$ is odd).

Next, they only needed to find out $n-\lceil k/2 \rceil$ internally
disjoint trees connecting $S$ in $G$, where $G=K_n$ for $k$ even;
$G=K_n\setminus M$ and $M$ is an edge set such that
$|M|=\frac{k-1}{2}$ for $k$ odd. Obviously, it only needs to
consider the case that $k$ is odd. But the difficulty is that each
edge of $E(G[S])\cup E_G[S,\bar{S}]$ belongs to a tree connecting
$S$ and can not be wasted. Fortunately, Nash-Williams-Tutte theorem
provides a perfect solution. They first derived the following lemma
from Theorem \ref{th1-1}.

\begin{lem}{\upshape \cite{LMS}}\label{lem5-2}
If $n$ is odd and $M$ is an edge set of the complete graph $K_n$
such that $0\leq |M|\leq \frac{n-1}{2}$, then $G=K_n\setminus M$
contains $\frac{n-1}{2}$ edge-disjoint spanning trees.
\end{lem}

They wanted to find out $\frac{k-1}{2}$ edge-disjoint spanning trees
in $G[S]$ (By the definition of internally disjoint trees, these
trees are internally disjoint trees connecting $S$, as required).
Then their basic idea is to seek for some edges in $G[S]$, and let
them together with the edges of $E_G[S,\bar{S}]$ form $n-k$
internally disjoint trees. They proved that there are indeed $n-k$
internally disjoint trees in the premise that $G[S]$ contains
$\frac{k-1}{2}$ edge-disjoint spanning trees. Actually, Lemma
\ref{lem5-2} can guarantee the existence of such $\frac{k-1}{2}$
trees. Then they found out $n-\frac{k-1}{2}$ internally disjoint
trees connecting $S$ and accomplished the proof of the following
theorem.

\begin{thm}{\upshape \cite{LMS}}\label{th5-3}
Let $k,n$ be two integers with $2\leq k\leq n$. Then for a connected
graph $G$ of order $n$, $\kappa_k(G)=n-\lceil\frac{k}{2}\rceil$
if and only if $G=K_n$ for $k$ even; $G=K_n\setminus M$ for $k$ odd,
where $M$ is an edge set such that $0\leq |M|\leq \frac{k-1}{2}$.
\end{thm}

Combining Theorem \ref{th5-3} and Observation \ref{obs2-1}, they
obtained the following theorem for $\lambda_k(G)$.

\begin{thm}{\upshape \cite{LMS}}\label{th5-4}
Let $k,n$ be two integers with $2\leq k\leq n$. Then for a connected
graph $G$ of order $n$, $\lambda_k(G)=n-\lceil\frac{k}{2}\rceil$ if
and only if $G=K_n$ for $k$ even; $G=K_n\setminus M$ for $k$ odd,
where $M$ is an edge set such that $0\leq |M|\leq \frac{k-1}{2}$.
\end{thm}

\subsection{Graphs with $\kappa_k(G)=n-\lceil\frac{k}{2}\rceil-1$ and
$\lambda_k(G)=n-\lceil\frac{k}{2}\rceil-1$}

As a continuation of their investigation, Li and Mao later turned
their attention to characterizing graphs $G$ with
$\kappa_k(G)=n-\lceil\frac{k}{2}\rceil-1$ and
$\lambda_k(G)=n-\lceil\frac{k}{2}\rceil-1$. One may notice that
$\kappa_k(G)=n-\lceil\frac{k}{2}\rceil$ if and only if $G$ itself is
the complete graph $K_n$ for $k$ even. So for $k$ even it is
possible to continue to characterize
$\kappa_k(G)=n-\lceil\frac{k}{2}\rceil-1$ by deleting more edges
from the complete graph $K_n$.

\begin{thm}{\upshape \cite{LMW2}}\label{th5-5}
Let $n$ and $k$ be two integers such that $k$ is even and $4\leq
k\leq n$, and $G$ be a connected graph of order $n$. Then
$\kappa_k(G)=n-\frac{k}{2}-1$ if and only if $G=K_n\setminus M$
where $M$ is an edge set such that $1\leq \Delta(K_n[M])\leq
\frac{k}{2}$ and $1\leq |M|\leq k-1$.
\end{thm}

Different from the proof of Theorem \ref{th5-3}, in order to find
$n-\frac{k}{2}-1$ internally disjoint trees in $G[S]\cup
G[S,\bar{S}]$, they designed a procedure to emphasize seeking for
some edges ``evenly" in $G[S]$, and let them together with the edges
of $E_G[S,\bar{S}]$ form $n-k$ internally disjoint trees
$T_1,T_2,\cdots,T_{n-k}$ with its root $w_1,w_2,\cdots,w_{n-k}\in
\bar{S}$, respectively. Applying this procedure designed by them,
they proved that the remaining edges in $G[S]$ can form
$\frac{k-2}{2}$ spanning trees, which are also $\frac{k-2}{2}$
internally disjoint trees connecting $S$. These trees together with
$T_1,T_2,\cdots,T_{n-k}$ are $n-\frac{k}{2}-1$ internally disjoint
trees connecting $S$, accomplishing the proof of the above theorem.

With the help of Theorem \ref{th5-5} and Observation \ref{obs2-1},
they obtained the following theorem for $\lambda_k(G)$.

\begin{thm}{\upshape \cite{LMW2}}\label{th5-6}
Let $n$ and $k$ be two integers such that $k$ is even and $2\leq
k\leq n$, and $G$ be a connected graph of order $n$. Then
$\lambda_k(G)=n-\frac{k}{2}-1$ if and only if $G=K_n\setminus M$
where $M$ is an edge set satisfying one of the following conditions:

$(1)$ $\Delta(K_n[M])=1$ and $1\leq |M|\leq
\lfloor\frac{n}{2}\rfloor$;

$(2)$ $2\leq \Delta(K_n[M])\leq \frac{k}{2}$ and $1\leq |M|\leq
k-1$.
\end{thm}

By Nash-Williams-Tutte theorem, they luckily characterized the
graphs attaining the upper bound and graphs with
$\kappa_k(G)=n-\lceil\frac{k}{2}\rceil-1$ and
$\kappa_k(G)=n-\lceil\frac{k}{2}\rceil-1$ for $k$ even. But, for $k$
odd, it is not easy to characterize the graphs with
$\kappa_k(G)=n-\lceil\frac{k}{2}\rceil-1$. So, Li, Li, Mao and Sun
considered the case $k=3$, namely, they considered graphs with
$\kappa_3(G)=n-3$.

\begin{thm} {\upshape \cite{LLMS}}\label{th5-7}
Let $G$ be a connected graph of order $n \ (n\geq 3)$. Then
$\kappa_3(G)=n-3$ if and only if $G$ is a graph obtained from the
complete graph $K_n$ by deleting an edge set $M$ such that
$K_n[M]=P_4$ or $K_n[M]=P_3\cup r P_2 \ (r=1,2)$ or $K_n[M]=C_3\cup
r P_2 \ (r=1,2)$ or $K_n[M]=s P_2 \ ( 2\leq s\leq
\lfloor\frac{n}{2}\rfloor)$.
\end{thm}

But, for the edge case, Li and Mao \cite{LM4} showed that the
statement is different.

\begin{thm}{\upshape \cite{LM4}}\label{th5-8}
Let $G$ be a connected graph of order $n$. Then $\lambda_3(G)=n-3$
if and only if $G$ is a graph obtained from the complete graph $K_n$
by deleting an edge set $M$ such that $K_n[M]=r P_2 \ (2\leq r\leq
\lfloor\frac{n}{2}\rfloor)$ or $K_n[M]=P_4\cup sP_2 \ (0\leq s\leq
\lfloor\frac{n-4}{2}\rfloor)$ or $K_n[M]=P_3\cup t P_2 \ (0\leq
t\leq \lfloor\frac{n-3}{2}\rfloor)$ or $K_n[M]=C_3\cup t P_2 \
(0\leq t\leq \lfloor\frac{n-3}{2}\rfloor)$.
\end{thm}

\section{Nordhaus-Gaddum-type results}

Let $\mathcal {G}(n)$ denote the class of simple graphs of order $n$
and $\mathcal {G}(n,m)$ the subclass of $\mathcal {G}(n)$ having
graphs with $n$ vertices and $m$ edges. Give a graph parameter
$f(G)$ and a positive integer $n$, the \emph{Nordhaus-Gaddum
(\textbf{N-G}) Problem} is to determine sharp bounds for: $(1)$
$f(G)+f(\overline{G})$ and $(2)$ $f(G)\cdot f(\overline{G})$, as $G$
ranges over the class $\mathcal {G}(n)$, and characterize the
extremal graphs. The Nordhaus-Gaddum type relations have received
wide attention; see a recent survey paper \cite{Aouchiche} by
Aouchiche and Hansen.

Alavi and Mitchem in \cite{Alavi} investigated Nordhaus-Gaddum-type
results for the classical connectivity and edge-connectivity in $\mathcal
{G}(n)$. Achuthan and Achuthan \cite{Achuthan} considered the same
problem in $\mathcal {G}(n,m)$.

\subsection{Results for graphs in  $\mathcal {G}(n)$}

Li and Mao \cite{LM} investigated the Nordhaus-Gaddum type relations
on the generalized edge-connectivity. At first, they focused on the
graphs in $\mathcal {G}(n)$.

\begin{thm}{\upshape \cite{LM}}\label{th6-1}
Let $G\in \mathcal {G}(n)$, and $k,n$ be two integers with $2\leq
k\leq n$. Then

$(1)$ $1\leq \lambda_k(G)+\lambda_k(\overline{G})\leq n-\lceil k/2
\rceil$;

$(2)$ $0\leq \lambda_k(G)\cdot \lambda_k(\overline{G})\leq
\big[\frac{n-\lceil k/2 \rceil}{2}\big]^2$.

Moreover, the upper and lower bounds are sharp.
\end{thm}

The following observation indicates the graphs attaining the above
lower bound.

\begin{obs}{\upshape \cite{LM}}\label{obs6-2}
$\lambda_k(G)\cdot \lambda_k(\overline{G})=0$ if and only if $G$ or
$\overline{G}$ is disconnected.
\end{obs}

For $n\geq 5$, $\mathcal {G}_n^1$ is a graph class as shown in
Figure 6.1 $(a)$ such that $\lambda(G)=1$ and $d_{G}(v_1)=n-1$ for
$G\in \mathcal {G}_n^1$, where $v_1\in V(G)$; $\mathcal {G}_n^2$ is
a graph class as shown in Figure 6.1 $(b)$ such that $\lambda(G)=2$
and $d_{G}(u_1)=n-1$ for $G\in \mathcal {G}_n^2$, where $u_1\in
V(G)$; $\mathcal {G}_n^3$ is a graph class as shown in Figure 6.1
$(c)$ such that $\lambda(G)=2$ and $d_{G}(v_1)=n-1$ for $G\in
\mathcal {G}_n^3$, where $v_1\in V(G)$; $\mathcal {G}_n^4$ is a
graph class as shown in Figure 6.1 $(d)$ such that $\lambda(G)=2$.

\begin{figure}[!hbpt]
\begin{center}
\includegraphics[scale=0.9]{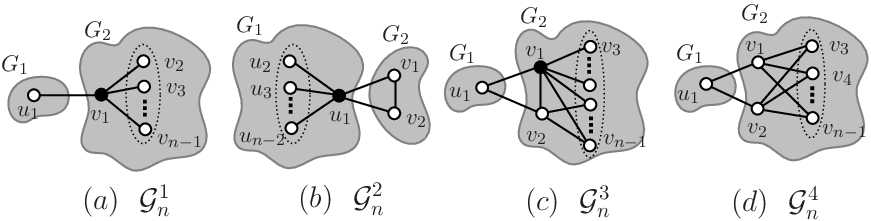}
\end{center}
\begin{center}
\caption{The graph classes $\mathcal {G}_n^i \ (1\leq i\leq 4)$}
\end{center}\label{fig5-1}
\end{figure}

As we know, it is not easy to characterize the graphs with
$\lambda_k(G)=1$, even with $\lambda_3(G)=1$. So, Li and Mao wanted
to add some conditions to attack such a problem. Motivated by such
an idea, they hope to characterize the graphs with
$\lambda_k(G)+\lambda_k(\overline{G})=1$. Actually, the
Norhaus-Gaddum-type problems also need to characterize the extremal
graphs attaining the bounds.

\begin{pro}{\upshape \cite{LM}}\label{pro6-3}
$\lambda_k(G)+\lambda_k(\overline{G})=1$ if and only if $G$
(symmetrically, $\overline{G}$) satisfies one of the following
conditions:

$\bullet$ $G\in \mathcal {G}_n^1$ or $G\in \mathcal {G}_n^2$;

$\bullet$ $G\in \mathcal {G}_n^3$ and there exists a component $G_i$
of $G\setminus v_1$ such that $G_i$ is a tree and $|V(G_i)|<k$;

$\bullet$ $G\in \{P_3,C_3\}$ for $k=n=3$, or $G\in
\{C_4,K_4\setminus e\}$ for $k=n=4$, or $G=K_{3,3}$ for $k=n=6$, or
$G=K_{2,n-2}$ for $k=n-1$ and $n\geq 5$, or $G=C_4$ for $k=n-1=3$,
or $G\in \{K_{2,n-2}^{+},K_{2,n-2}\}$ for $k=n$ and $n\geq 5$ where
$K_{2,n-2}^{+}$ denote the graph obtained from the complete
bipartite graph $K_{2,n-2}$ by adding one edge in the part with
$n-2$ vertices.
\end{pro}

Let us focus on $(1)$ of Theorem \ref{th6-1}. If one of $G$ and
$\overline{G}$ is disconnected, we can characterize the graphs
attaining the upper bound by Theorem \ref{th5-4}.

\begin{pro}{\upshape \cite{LM}}\label{pro6-4}
For any graph $G$ of order $n$, if $G$ is disconnected, then
$\lambda_k(G)+ \lambda_k(\overline{G})=n-\lceil\frac{k}{2}\rceil$ if
and only if $\overline{G}=K_n$ for $k$ even;
$\overline{G}=K_n\setminus M$ for $k$ odd, where $M$ is an edge set
such that $0\leq |M|\leq \frac{k-1}{2}$.
\end{pro}

If both $G$ and $\overline{G}$ are connected, we can obtain a
structural property of the graphs attaining the upper bound although
it seems too difficult to characterize them.

\begin{pro}{\upshape \cite{LM}}\label{pro6-5}
If $\lambda_k(G)+\lambda_k(\overline{G})=n-\lceil\frac{k}{2}\rceil$,
then $\Delta(G)-\delta(G)\leq \lceil\frac{k}{2}\rceil-1$.
\end{pro}

One can see that the graphs with
$\lambda_k(G)+\lambda_k(\overline{G})=n-\lceil\frac{k}{2}\rceil$
must have a uniform degree distribution. By this property, Li and Mao
constructed a graph class to show that the two upper bounds of
Theorem \ref{th6-1} are tight for $k=n$.

\noindent \textbf{Example 6.1.}~~Let $n,r$ be two positive integers
such that $n=4r+1$. From Theorem \ref{th2-5}, we know that
$\kappa_{n}(K_{2r,2r+1})=\lambda_{n}(K_{2r,2r+1})=r$. Let $\mathcal
{E}$ be the set of the edges of these $r$ spanning trees in
$K_{2r,2r+1}$. Then there exist $2r(2r+1)-4r^2=2r$ remaining edges
in $K_{2r,2r+1}$ except for the edges in $\mathcal {E}$. Let $M$ be the
set of these $2r$ edges. Set $G=K_{2r,2r+1}\setminus M$. Then
$\lambda_{n}(G)=r$, $M\subseteq E(\overline{G})$ and $\overline{G}$
is a graph obtained from two cliques $K_{2r}$ and $K_{2r+1}$ by
adding $2r$ edges in $M$ between them, that is, one endpoint of each
edge belongs to $K_{2r}$ and the other endpoint belongs to
$K_{2r+1}$. Note that $E(\overline{G})=E(K_{2r})\cup M\cup
E(K_{2r+1})$. Now we show that $\lambda_{n}(\overline{G})\geq r$. As
we know, $K_{2r}$ contains $r$ Hamiltonian paths, say
$P_{1},P_{2},\cdots,P_{r}$, and so does $K_{2r+1}$, say
$P_{1}',P_{2}',\cdots,P_{r}'$. Pick up $r$ edges from $M$, say
$e_1,e_2,\cdots,e_r$, let $T_i=P_{i}\cup P_{i}'\cup e_i(1\leq i\leq
r)$. Then $T_{1},T_{2},\cdots,T_{r}$ are $r$ spanning trees in
$\overline{G}$, namely, $\lambda_{n}(\overline{G})\geq r$. Since
$|E(\overline{G})|={{2r}\choose{2}}+{{2r+1}\choose{2}}+2r=4r^2+2r$
and each spanning tree uses $4r$ edges, these edges can form at most
$\lfloor\frac{4r^2+2r}{4r}\rfloor=r$ spanning trees, that is,
$\lambda_{n}(\overline{G})\leq r$. So $\lambda_{n}(\overline{G})=r$.
Clearly, $\lambda_{n}(G)+\lambda_{n}(\overline{G})=2r=\frac{n-1}{2}
=n-\lceil\frac{n}{2}\rceil$ and $\lambda_{n}(\overline{G})\cdot
\lambda_{n}(\overline{G})=r^2=\big[\frac{n-\lceil n/2
\rceil}{2}\big]^2$.

Li, Mao and Sun \cite{LMS} were concerned with analogous
inequalities involving the generalized $k$-connectivity for the
graphs in $\mathcal {G}(n)$.

\begin{thm}{\upshape \cite{LMS}}\label{th6-10}
Let $G\in \mathcal {G}(n)$, and $k,n$ be two integers with $2\leq
k\leq n$. Then

$(1)$ $1\leq \kappa_k(G)+\kappa_k(\overline{G})\leq n-\lceil k/2
\rceil$;

$(2)$ $0\leq \kappa_k(G)\cdot \kappa_k(\overline{G})\leq
[\frac{2n-\lceil k \rceil}{4}]^2$.

Moreover, the upper and lower bounds are sharp.
\end{thm}

\subsection{Results for graphs in $\mathcal {G}(n,m)$}

Then, Li and Mao also focused on the graphs in $\mathcal {G}(n,m)$ in \cite{LM}.
Let us begin with another problem, called the maximum connectivity of a
graph. It was pointed out by Harary \cite{Harary} that given the
number of vertices and edges of a graph, the largest connectivity
possible can also be read out of the inequality $\kappa(G)\leq
\lambda(G)\leq \delta(G)$.

\begin{thm}{\upshape \cite{Harary}}\label{th6-6}
For each pair of integers $n,m$ with $0\leq n-1\leq m\leq {{n}\choose{2}}$,

$$\kappa(G)\leq \lambda(G)\leq \Big\lfloor \frac{2m}{n}\Big\rfloor,$$
where the maximum is taken over all graphs $G\in \mathcal
{G}(n,m)$.
\end{thm}

Li and Mao considered the similar problem for the generalized
edge-connectivity.

\begin{cor}{\upshape \cite{LM}}\label{cor6-7}
For any graph $G\in \mathcal {G}(n,m)$ and $3\leq k\leq n$,
$\lambda_k(G)=0$ for $m<n-1$; $\lambda_k(G)\leq
\lfloor\frac{2m}{n}\rfloor$ for $m\geq n-1$.
\end{cor}

Although the above bound of $\lambda_k(G)$ is the same as
$\lambda(G)$, the graphs attaining the upper bound seems to be very
rare. Actually, we can obtain some structural properties of these
graphs.

\begin{pro}{\upshape \cite{LM}}\label{pro6-8}
For any $G\in \mathcal {G}(n,m)$ and $3\leq k\leq n$, if
$\lambda_k(G)=\lfloor\frac{2m}{n}\rfloor$ for $m\geq n-1$, then

$\bullet$ $\frac{2m}{n}$ is not an integer;

$\bullet$ $\delta(G)=\lfloor\frac{2m}{n}\rfloor$;

$\bullet$ for $u,v\in V(G)$ such that
$d_G(u)=d_G(v)=\lfloor\frac{2m}{n}\rfloor$, $uv\notin E(G)$.
\end{pro}

By Theorem \ref{th1-1} and Corollary \ref{cor6-7}, they derived the
following theorem.

\begin{thm}{\upshape \cite{LM}}\label{th6-9}
Let $G\in \mathcal {G}(n,m)$. For $n\geq 6$, we have

$(1)$ $L(n,m)\leq \lambda_k(G)+\lambda_k(\overline{G})\leq M(n,m)$;

$(2)$ $0\leq \lambda_k(G)\cdot \lambda_k(\overline{G})\leq N(n,m)$,

where
$$
L(n,m)=\left\{
\begin{array}{ll}
max\{1,\lfloor\frac{1}{2}(n-2-m)\rfloor\}, &if~\lfloor\frac{n}{3}\rfloor+1\leq m\leq {{n}\choose{2}};\\
min \{n-2m-1, \lfloor \frac{n}{2}-\frac{2m}{n-1}\rfloor\},&if~0\leq
m\leq \lfloor\frac{n}{3}\rfloor.
\end{array}
\right.
$$

$$
M(n,m)=\left\{
\begin{array}{ll}
n-\lceil\frac{k}{2}\rceil, &if~m\geq n-1,\\
&~or~k~is~even~and~m=0,\\
&~or~k~is~odd~and~0\leq m\leq \frac{k-1}{2};\\
n-\lceil\frac{k}{2}\rceil-1,&if~k~is~even~and~1\leq m<
n-1,\\
&~or~k~is~odd~and~\frac{k+1}{2}\leq m<n-1.
\end{array}
\right.
$$

$$
N(n,m)=\left\{
\begin{array}{ll}
0, &if~0\leq m\leq n-2;\\
(\frac{2m}{n}-1)(n-2-\frac{2m}{n}),&if~m\geq n-1~and~2m\equiv
0(mod~n);\\
\lfloor\frac{2m}{n}\rfloor(n-2-\lfloor\frac{2m}{n}\rfloor),&otherwise.
\end{array}
\right.
$$

Moreover, the upper and lower bounds are sharp.
\end{thm}

\section{Results for graph operations}

In this section we will survey the results for line graphs and graph
products.

\subsection{Results for line graphs}

Chartrand and Steeart \cite{Steeart} investigated the relation
between the connectivity and edge-connectivity of a graph and its
line graph. They proved that if $G$ is a connected graph, then $(1)$
$\kappa(L(G))\geq \lambda(G)$ if $\lambda(G)\geq 2$; $(2)$
$\lambda(L(G))\geq 2\lambda(G)-2$; $(3)$ $\kappa(L(L(G)))\geq
2\kappa(G)-2$. With the help of Proposition \ref{pro4-5}, Li, Mao
and Sun also considered the generalized $3$-connectivity and
$3$-edge-connectivity for line graphs.

\begin{pro}{\upshape \cite{LMS}}\label{pro7-1}
If $G$ is a connected graph, then

$(1)$ $\lambda_3(G)\leq \kappa_3(L(G))$.

$(2)$ $\lambda_3(L(G))\geq \frac{3}{2}\lambda_3(G)-2$.

$(3)$ $\kappa_3(L(L(G))\geq \frac{3}{2}\kappa_3(G)-2$.
\end{pro}

First, they proved $(1)$ of this theorem. Next, combining
Proposition \ref{pro4-5} with $(1)$ of Proposition \ref{pro7-1},
they derived $(2)$ and $(3)$ of Proposition \ref{pro7-1}. One can
check that $(1)$ of this proposition is sharp since $G=C_n$ can
attain this bound.

Let $L^0(G)=G$ and $L^1(G)=L(G)$. Then for $r\geq 2$, the $r$-$th$
iterated line graph $L^r(G)$ of $G$ is defined by $L(L^{r-1}(G))$. The next
statement follows immediately from Proposition \ref{pro7-1} and a
routine application of induction.

\begin{cor}{\upshape \cite{LMS}}\label{cor7-2}
$\lambda_3(L^r(G))\geq (\frac{3}{2})^r(\lambda_3(G)-4)+4$, and
$\kappa_3(L^r(G))\geq
(\frac{3}{2})^{\lfloor\frac{r}{2}\rfloor}(\kappa_3(G)-4)+4$.
\end{cor}

\subsection{Results for graph products}

Product networks were proposed based upon the idea of using the
cross product as a tool for ``combining'' two known graphs with
established properties to obtain a new one that inherits properties
from both \cite{DayA}. Recently, there has been an increasing
interest in a class of interconnection networks called Cartesian
product networks; see \cite{Bao, DayA, Ku}. In \cite{Ku}, Ku, Wang
and Hung studied the problem of constructing the maximum number of
edge-disjoint spanning trees in Cartesian product networks, and gave
a sharp lower bound of $\kappa_n(G\Box H)$.

\begin{thm}{\upshape \cite{Ku}}\label{th7-3}
For two connected graphs $G$ and $H$, $\kappa_n(G \Box H)\geq
\kappa_n(G)+\kappa_n(H)-1$. Moreover, the lower bound is sharp.
\end{thm}

But the upper bound of $\kappa_n(G\Box H)$ is still unknown. A
natural question is to study the following problems:

$\bullet$ Give sharp upper and lower bounds of $\kappa_k(G*H)$,
where $*$ is a kind of graph product.

$\bullet$ Give sharp upper and lower bounds of $\lambda_k(G*H)$,
where $*$ is a kind of graph product.

Sabidussi in \cite{Sabidussi} derived a result on the classical connectivity
of Cartesian product graphs: for two connected graphs $G$ and $H$,
$\kappa(G\square H)\geq \kappa(G)+\kappa(H)$. But we mention that it
was incorrectly claimed in (\cite{Hammack}, page 308) that $\kappa(G\Box
H)=\kappa(G)+\kappa(H)$ holds for any connected $G$ and $H$. In
\cite{Spacapan}, \u{S}pacapan proved that $\kappa(G\Box
H)=\min\{\kappa(G)|V(H)|,\kappa(H)|V(G)|,\delta(G)+\delta(H)\}$ for
two nontrivial graphs $G$ and $H$.

\subsubsection{The case $k=3$}

In \cite{LLSun}, Li, Li and Sun investigated the generalized
$3$-connectivity of Cartesian product graphs. Their results could be
seen as a generalization of Sabidussi's result. As usual, in order
to get a general result, they first began with a special case.

\begin{pro}{\upshape \cite{LLSun}}\label{pro7-4}
Let $G$ be a graph and $P_m$ be a path with $m$ edges. The following
assertions hold:

$(1)$ If $\kappa_3(G)=\kappa(G)\geq 1$, then $\kappa_3(G\square
P_m)\geq \kappa_3(G)$. Moreover, the bound is sharp.

$(2)$ If $1\leq \kappa_3(G)<\kappa(G)$, then $\kappa_3(G\square
P_m)\geq \kappa_3(G)+1$. Moreover, the bound is sharp.

\end{pro}

Note that $Q_n\cong P_2\square P_2\square\cdots \square P_2$, where
$Q_n$ is the $n$-hypercube. They got the following corollary.

\begin{cor}{\upshape \cite{LLSun}}\label{cor7-5}
Let $Q_n$ be the $n$-hypercube with $n\geq 2$. Then
$\kappa_3(Q_n)=n-1$.
\end{cor}

\noindent \textbf{Example 4}. Let $H_1$ and $H_2$ be two complete
graphs of order $n$, and let $V(H_1)=\{u_1,u_2,\cdots,$ $u_n\}$,
$V(H_2)=\{v_1,v_2,\cdots,v_n\}$. We now construct a graph $G$ as
follows: $V(G)=V(H_1)\cup V(H_2)\cup \{w\}$ where $w$ is a new
vertex; $E(G)=E(H_1)\cup E(H_2)\cup \{u_iv_j|1\leq i,j\leq n\}\cup
\{wu_i|1\leq i\leq n\}$. It is easy to check that $\kappa_3(G\square
K_2)=\kappa_3(G)=n$.

They showed that the bounds of $(1)$ and $(2)$ in Proposition
\ref{pro7-4} are sharp by Example $4$ and Corollary \ref{cor7-5}.

Next, they studied the generalized $3$-connectivity of the Cartesian
product of a graph $G$ and a tree $T$, which will be used in Theorem
\ref{th7-7}.

\begin{pro}{\upshape \cite{LLSun}}\label{pro7-6}
Let $G$ be a graph and $T$ be a tree. The following assertions hold:

$(1)$ If $\kappa_3(G)=\kappa(G)\geq 1$, then $\kappa_3(G\square
T)\geq \kappa_3(G)$. Moreover, the bound is sharp.

$(2)$ If $1\leq \kappa_3(G)<\kappa(G)$, then $\kappa_3(G\square
T)\geq \kappa_3(G)+1$. Moreover, the bound is sharp.

\end{pro}

The bounds of $(1)$ and $(2)$ in Proposition \ref{pro7-6} are sharp
by Example $4$ and Corollary \ref{cor7-5}.

They mainly investigated the generalized $3$-connectivity of the
Cartesian product of two connected graphs $G$ and $H$. By
decomposing $H$ into some trees connecting $2$ vertices or $3$
vertices, they considered the Cartesian product of a graph $G$ and a
tree $T$ and obtained Theorem \ref{th7-7} by Proposition
\ref{pro7-6}.

\begin{thm}{\upshape \cite{LLSun}}\label{th7-7}
Let $G$ and $H$ be connected graphs such that
$\kappa_3(G)>\kappa_3(H)$. The following assertions hold:

$(1)$ If $\kappa(G)=\kappa_3(G)$, then $\kappa_3(G\square H)\geq
\kappa_3(G)+\kappa_3(H)-1$. Moreover, the bound is sharp.

$(2)$ If $\kappa(G)>\kappa_3(G)$, then $\kappa_3(G\square H)\geq
\kappa_3(G)+\kappa_3(H)$. Moreover, the bound is sharp.
\end{thm}

They also showed that the bounds of $(1)$ and $(2)$ in Theorem
\ref{th7-7} are sharp. Let $K_n$ be a complete graph with $n$
vertices, and $P_m$ be a path with $m$ vertices, where $m\geq 2$.
Since $\kappa_3(P_m)=1$ and $\kappa_3(K_n)=n-2$, it is easy to see
that $\kappa_3(K_n\square P_m)=n-1$. Thus, $K_n\square P_m$ is a
sharp example for $(1)$. For $(2)$, Example $4$ is a sharp one.

Lexicographic product is one of the standard products, which are
studied extensively; see \cite{Hammack}. Recently, some applications
in networks of the lexicographic product were studied; see
\cite{Blasiak, Feng, LXZW}. Yang and Xu \cite{Yang} investigated the
classical connectivity of the lexicographic product of two graphs:
For two graphs $G$ and $H$, if $G$ is non-trivial, non-complete and
connected, then $\kappa(G\circ H)=\kappa(G)|V(H)|$.

Using Fan Lemma (\cite{West}, page 170) and Expansion Lemma
(\cite{West}, page 162), Li and Mao \cite{LM3} obtained the following
lower bound of $\kappa_3(G\circ H)$, which could be seen as an
extension of Yang and Xu' result.

\begin{thm}{\upshape \cite{LM3}}\label{th7-8}
Let $G$ and $H$ be two connected graphs. Then
$$
\kappa_3(G\circ H)\geq \kappa_3(G)|V(H)|.
$$
Moreover, the bound is sharp.
\end{thm}

For a tree $T$ and a connected graph $H$, they showed that
$\kappa_3(T\circ H) =|V(H)|$, which can be seen as an improvement of
Theorem \ref{th7-8}. From Theorem \ref{th7-7}, one may wonder
whether $\kappa_3(T\Box H)=\kappa_3(T)+\kappa_3(H)-1$ for a
connected graph $H$ and a tree $T$ (note that
$\kappa_3(T)=\kappa(T)=1$). For example, let $T=P_3$ and $H=K_4$.
Then $\kappa_3(T)=\kappa(T)=1$ and $\kappa_3(H)=2$. One can check
that $\kappa_3(T\Box H)=3>2=\kappa_3(T)+\kappa_3(H)-1$. So the
equality does not hold for the Cartesian product of a tree and a
connected graph.

For the edge version of the above mentioned problem, Yang and Xu
\cite{Yang} also derived that $\lambda(G\circ H)=\min\{
\lambda(G)|V(H)|^2,\delta(H)+\delta(G)|V(H)|\}$ for a connected
graph $G$ and a non-trivial graphs $H$. Recently, Li, Yue and Zhao
\cite{LYZ} gave a lower bound of $\lambda_3(G \circ H)$.

\begin{thm}{\upshape \cite{LYZ}} \label{th7-10}
Let $G$ and $H$ be a connected graph. Then
$$
\lambda_3(G\circ H)\geq \lambda_3(H)+\lambda_3(G)|V(H)|.
$$
Moreover, the lower bound is sharp.
\end{thm}

From Theorems \ref{th4-4} and \ref{th4-7}, and Yang and Xu' s result,
Li and Mao \cite{LM3} derived a upper bound of $\kappa_3(G\circ H)$.

\begin{thm}{\upshape \cite{LM3}} \label{th7-11}
Let $G$ and $H$ be two connected graphs. If $G$ is non-trivial and
non-complete, then $\kappa_3(G\circ H)\leq \lfloor
\frac{4}{3}\kappa_3(G)+r-\frac{4}{3}\lceil\frac{r}{2}\rceil\rfloor
|V(H)|$, where $r\equiv \kappa(G) \ (mod~4)$. Moreover, the bound is
sharp.
\end{thm}

The graph $P_n\circ P_3 \ (n\geq 4)$ indicates that both the lower
bound of Theorem \ref{th7-8} and the upper bound of Theorem
\ref{th7-11} are sharp.

In the same paper, they also derived the following upper bound of
$\kappa_3(G\Box H)$ from Theorems \ref{th4-4} and \ref{th4-7}, and
\u{S}pacapan's result.

\begin{thm}{\upshape \cite{LM3}} \label{th7-12}
Let $G$ and $H$ be two connected graphs. Then $\kappa_3(G\Box H)\leq
\min\{\lfloor
\frac{4}{3}\kappa_3(G)+r_1-\frac{4}{3}\lceil\frac{r_1}{2}\rceil\rfloor|V(H)|,
\lfloor\frac{4}{3}\kappa_3(H)+r_2-\frac{4}{3}\lceil\frac{r_2}{2}\rceil\rfloor
|V(G)|, \delta(G)+\delta(H)\}$, where $r_1\equiv \kappa(G) \ (mod \
4)$ and $r_2\equiv \kappa(H) \ (mod \ 4)$. Moreover, the bound is
sharp.
\end{thm}

The graph $P_n\circ P_m \ (n\geq 4, \ m\geq 4)$ is a sharp example
for the above theorem.

In \cite{LYZ}, Li, Yue and Zhao also obtained an upper bound of
$\lambda_3(G \circ H)$.

\begin{thm}{\upshape \cite{LYZ}} \label{th7-13}
Let $G$ be a connected graph, and $H$ be a non-trivial graph. Then
$\lambda_3(G\circ H)\leq \min \{\lfloor \frac{4\lambda_3(G)+2}{3}
\rfloor |V(H)|^2,\delta(H)+\delta(G)|V(H)|\}$. Moreover, the upper
bound is sharp.
\end{thm}

The graph $P_{t}\circ P_{n-t}$ is a sharp example for both Theorem
\ref{th7-10} and Theorem \ref{th7-13}.

\subsubsection{The case $k=n$}

Like that in \cite{Ku} for Cartesian product, Li, Li, Mao and Yue
\cite{LLMY} investigated the spanning tree packing number $\sigma$ of
lexicographic product graphs and hoped to obtain a lower bound of
$\sigma(G \circ H)$. Usually, in order to give such a lower bound,
one must find out as many spanning trees in $G\circ H$ as possible.
The following two procedures are given in their paper:

$\bullet$ \textbf{Graph decomposition}: Decompose the graph $G \circ
H$ into desired small graphs, such as parallel forests, good cycles,
and trees in $G \circ H$ corresponding to the spanning tree of $G$
or $H$ (see \cite{LLMY}).

$\bullet$ \textbf{Graph combination}: The above small graphs are
divided into groups each of which contains different kind of small
graphs. Then, combine the small graphs in each group to obtain a
spanning tree of $G \circ H$.

After the second procedure, they obtained the maximum number of
edge-disjoint spanning trees in $G \circ H$, which is a lower bound
of $\kappa_n(G \circ H)$.

\begin{thm}{\upshape \cite{LLMY}} \label{th7-9}
Let $G$ and $H$ be two connected graphs. $\kappa_n(G)=k$,
$\kappa_n(H)=\ell$, $|V(G)|=n_1$, and $|V(H)|=n_2$. Then

$(1)$ if $k n_2=\ell n_1$, then $\kappa_n(G \circ H)\geq k n_2(=\ell
n_1)$;

$(2)$ if $\ell n_1>k n_2$, then $\kappa_n(G \circ H)\geq
kn_2-\lceil\frac{k n_2-1}{n_1}\rceil+\ell-1$;

$(3)$ if $\ell n_1<k n_2$, then $\kappa_n(G \circ H)\geq
kn_2-\lceil\frac{2kn_2-4}{n_1+1}\rceil+\ell-1$.

Moreover, the lower bounds are sharp.
\end{thm}

To show the sharpness of the above lower bounds of Theorem
\ref{th7-9}, they considered the following example.

\noindent\textbf{Example 7.1}. $(1)$ Let $G$ and $H$ be two
connected graphs with $|V(G)|=n_1$ and $|V(H)|=n_2$ which can be
decomposed into exactly $k$ and $\ell$ edge-disjoint spanning trees
of $G$ and $H$, respectively, satisfying $k n_2=\ell n_1$. Then
$\kappa_n(G \circ H)=kn_2=\ell n_1$.

$(2)$ Let $G=P_3$ and $H=K_4$. Clearly, $\kappa_n(G)=k=1$,
$\kappa_n(H)=\ell=2$, $|V(G)|=n_1=3$, $|V(H)|=n_2=4$. Therefore,
$\kappa_n(G \circ H)=4=kn_2-\lceil\frac{k n_2-1}{n_1}\rceil+\ell-1$.

$(3)$ Let $G=P_2$ and $H=P_3$. Clearly, $\kappa_n(G)=k=1$,
$\kappa_n(H)=\ell=1$, $|V(G)|=n_1=2$, $|V(H)|=n_2=3$, $|E(G)|=1$.
Then $\kappa_n(G \circ H)=2=
kn_2-\lceil\frac{2kn_2-4}{n_1+1}\rceil+\ell-1$.

\section{Extremal problems}

In this section, we survey the results on the extremal problems of
generalized connectivity and generalized edge-connectivity.

\subsection{The minimal size of a graph with given generalized
$k$-(edge-)connectivity}

Li, Li and Shi \cite{LLShi} determined the minimal number of edges
among graphs with $\kappa_3(G)= 2$, i.e., graphs $G$ of order $n$
and size $e(G)$ with $\kappa_3(G)=2$, that is,

\begin{thm}{\upshape \cite{LLShi}}\label{th8-1}
If $G$ is a graph of order $n$ with $\kappa_3(G)=2$, then $e(G)\geq
\lceil\frac{6}{5}n\rceil$. Moreover, the lower bound is sharp for
all $n\geq 4$ but $n=9,10$, whereas the best lower bound for
$n=9,10$ is $\lceil\frac{6}{5}n\rceil+1$.
\end{thm}

They constructed a graph class to show that the bound of Theorem
\ref{th8-1} is sharp.

\noindent \textbf{Example 8.1.} For a positive integer $t\neq 2$,
let $C=x_1y_1x_2y_2\cdots x_{2t}y_{2t}x_1$ be a cycle of length
$4t$. Add $t$ new vertices $z_1,z_2,\cdots,z_t$ to $C$, and join
$z_i$ to $x_i$ and $x_{i+t}$, for $1\leq i\leq t$. The resulting
graph is denoted by $H$. Then $\kappa_3(H)=2$; see Figure $8.1$.

\begin{figure}[!hbpt]
\begin{center}
\includegraphics[scale=0.75]{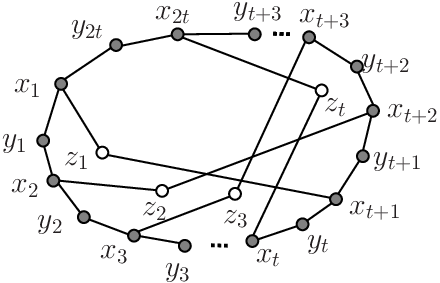}
\end{center}
\begin{center}
\caption{The graph $H$ with $\kappa_3(H)=2$}
\end{center}\label{fig6}
\end{figure}

Later, Li and Mao \cite{LM4} considered a generalization of the
above problem. Let $s(n,k,\ell)$ and $t(n,k,\ell)$ denote the
minimal number of edges of a graph $G$ of order $n$ with
$\kappa_k(G)=\ell \ (1\leq \ell\leq n-\lceil\frac{k}{2}\rceil)$ and
$\lambda_k(G)=\ell \ (1\leq \ell\leq n-\lceil\frac{k}{2}\rceil)$,
respectively.

From Theorem \ref{th8-1}, one can see that
$s(n,3,2)=\lceil\frac{6}{5}n\rceil$ for all $n\geq 4$ but $n=9,10$.
From Theorems \ref{th5-3} and \ref{th5-4}, we know that
$$
s(n,k,n-\lceil k/2\rceil)=t(n,k,n-\lceil k/2\rceil)=\left\{
\begin{array}{ll}
{n\choose 2},&for \ k \ even;\\
{n\choose 2}-\frac{k-1}{2},&for \ k \ odd.
\end{array}
\right.
$$
From Theorems \ref{th5-5} and \ref{th5-6}, we know that for $k$ even
$$
s(n,k,n-\lceil k/2\rceil-1)={n\choose 2}-k+1
$$
and
$$
t(n,k,n-\lceil k/2\rceil-1)={n\choose 2}-\lfloor n/2\rfloor.
$$

Li and Mao \cite{LM4} investigated $t(n,3,\ell)$ and derived the
following result.

\begin{thm}{\upshape \cite{LM4}} \label{th8-2}
Let $n$ be an integer with $n\geq 3$. Then

$(1)$ $t(n,3,n-2)={n\choose{2}}-1$;

$(2)$ $t(n,3,n-3)={n\choose{2}}-\lfloor\frac{n+3}{2}\rfloor$;

$(3)$ $t(n,3,1)=n-1$;

$(4)$ $t(n,3,\ell)\geq \big\lceil \frac{\ell(\ell+1)}{2\ell+1}n
\big\rceil$ for $n\geq 11$ and $2\leq \ell\leq n-4$. Moreover, the
bound is sharp.
\end{thm}

The complete bipartite graph $G=K_{\ell,\ell+1}$ is a sharp example
for the bound of Theorem \ref{th8-2}.

In \cite{SLi}, Li focused on the following problem: Given any
positive integer $n\geq 4$, is there a smallest integer $f(n)$ such
that every graph of order $n$ and size $e(G)\geq f(n)$ has
$\kappa_3(G)\geq 2$ ? She proved that every graph $G$ of order $n$
and size $e(G)=\frac{n^2}{2}-\frac{3n}{2}+3$ can be regarded as a
graph obtained from $K_n$ by deleting $n-3$ edges. Since
$\kappa_3(G)\geq 2$, $f(n)\leq \frac{n^2}{2}-\frac{3n}{2}+3$. On the
other hand, let $G$ be a graph obtained from $K_{n-1}$ by adding a
vertex $v$ and joining $v$ to one vertex of $K_{n-1}$. Clearly, the
order is $n$ and the size is $\frac{n^2}{2}-\frac{3n}{2}+2$. But
$\kappa_3(G)\leq \delta(G)=1$. So
$f(n)>\frac{n^2}{2}-\frac{3n}{2}+2$. Thus, the following result
is easily seen.

\begin{pro}{\upshape \cite{SLi}} \label{pro8-3}
Given any positive integer $n\geq 4$, there exists a smallest
integer $f(n)=\frac{n^2}{2}-\frac{3n}{2}+3$ such that every graph
$G$ of order $n$ and size $e(G)\geq f(n)$ has $\kappa_3(G)\geq 2$.
\end{pro}

Recall that a graph $G$ is \emph{minimal for $\kappa_k(G)=t$} if the
generalized $k$-connectivity of $G$ is $t$ but the generalized
$k$-connectivity of $G-e$ is less than $t$ for any edge $e$ of $G$.
Though it is easy to find the sharp lower bound of $e(G)$, very little
progress has been made on the sharp upper bound. So, Li phrased an
open problem as follows.

\noindent \textbf{Open Problem:} Let $G$ be a graph of order $n$ and
size $e(G)$ such that $G$ is minimal for $\kappa_3=2$. Find the
sharp upper bounds $g(n)$ of $e(G)$.

She proved that $2n-4\leq g(n)\leq 3n-10$, but the exact value of
$g(n)$ is still unknown.

\subsection{Maximum generalized local connectivity}

Recall that $\kappa(G)=\min\{\kappa_{G}(x,y)\,|\,x,y\in V(G), \
x\neq y\}$ is usually the connectivity of $G$. In contrast to this
parameter, $\overline{\kappa}(G)=\max\{\kappa_{G}(x,y)\,|\,x,y\in
V(G), \ x\neq y\}$, introduced by Bollob\'{a}s, is called the
\emph{maximum local connectivity} of $G$. The problem of determining
the smallest number of edges, $h(n;\overline{\kappa}\geq r)$, which
guarantees that any graph with $n$ vertices and
$h(n;\overline{\kappa}\geq r)$ edges will contain a pair of vertices
joined by $r$ internally disjoint paths was posed by Erd\"{o}s and
Gallai; see \cite{Bartfai} for details. Bollob\'{a}s
\cite{Bollobas1} considered the problem of determining the largest
number of edges, $f(n;\overline{\kappa}\leq \ell)$, for graphs with
$n$ vertices and local connectivity at most $\ell$, that is,
$f(n;\overline{\kappa}\leq \ell)=\max\{e(G)\,|\,|V(G)|=n\ and\
\overline{\kappa}(G)\leq \ell\}$. One can see that
$h_1(n;\overline{\kappa}\geq \ell+1)=f(n;\overline{\kappa}\leq
\ell)+1$. Similarly, let $\lambda_{G}(x,y)$ denote the local
edge-connectivity connecting $x$ and $y$ in $G$. Then
$\lambda(G)=\min\{\lambda_{G}(x,y)\,|\,x,y\in V(G), \ x\neq y\}$ and
$\overline{\lambda}(G)=\max\{\lambda_{G}(x,y)\,|\,x,y\in V(G), \
x\neq y\}$ are the edge-connectivity and maximum local
edge-connectivity, respectively. So the edge version of the above
problems can be given similarly. Set $g(n;\overline{\lambda}\leq
\ell)=\max\{e(G)\,|\, |V(G)|=n\ and\ \overline{\lambda}(G)\leq
\ell\}$. Let $h_2(n;\overline{\lambda}\geq r)$ denote the smallest
number of edges which guarantees that any graph with $n$ vertices
and $h_2(n;\overline{\kappa}\geq r)$ edges will contain a pair of
vertices joined by $r$ edge-disjoint paths. Similarly,
$h_2(n;\overline{\lambda}\geq \ell+1)= g(n;\overline{\lambda}\leq
\ell)+1$. The problem of determining the precise value of the
parameters $f(n;\overline{\kappa}\leq \ell)$,
$g(n;\overline{\lambda}\leq \ell)$, $h_1(n;\overline{\kappa}\geq
r)$, or $h_2(n;\overline{\kappa}\geq r)$ has obtained wide attention
and many results have been worked out; see \cite{Bollobas1,
Bollobas2, Bollobas3, Leonard1, Leonard2, Leonard3, Mader1, Mader2,
Thomassen}.

Similar to the classical maximum local connectivity, Li, Li and Mao
\cite{LLM} introduced the parameter
$\overline{\kappa}_k(G)=\max\{\kappa(S)\,|\,S\subseteq V(G), \
|S|=k\}$, which is called the \emph{maximum generalized local
connectivity} of $G$. There they considered the problem of determining
the largest number of edges, $f(n;\overline{\kappa}_k\leq \ell)$,
for graphs with $n$ vertices and maximal generalized local
connectivity at most $\ell$, that is, $f(n;\overline{\kappa}_k\leq
\ell)=\max\{e(G)\,|\,|V(G)|=n\ and\ \overline{\kappa}_k(G)\leq
\ell\}$. They also considered the smallest number of edges,
$h_1(n;\overline{\kappa}_k\geq r)$, which guarantees that any graph
with $n$ vertices and $h_1(n;\overline{\kappa}_k\geq r)$ edges will
contain a set $S$ of $k$ vertices such that there are $r$ internally
disjoint $S$-trees. It is easy to check that
$h_1(n;\overline{\kappa}_k\geq \ell+1)=f(n;\overline{\kappa}_k\leq
\ell)+1$ for $0\leq \ell \leq n-\lceil k/2\rceil-1$.

The edge version of these problems are also introduced and
investigated by Li and Mao in \cite{LM2}. Similarly, $g(n;\overline{\lambda}_k\leq
\ell)=\max\{e(G)\,|\,|V(G)|=n\ and\ \overline{\lambda}_k(G)\leq
\ell\}$, and $h_2(n;\overline{\lambda}_k\geq r)$ is the smallest
number of edges, $h_2(n;\overline{\lambda}_k\geq r)$, which
guarantees that any graph with $n$ vertices and
$h_2(n;\overline{\lambda}_k\geq r)$ edges will contain a set $S$ of
$k$ vertices such that there are $r$ edge-disjoint $S$-trees, and also
similarly, $h_2(n;\overline{\lambda}_k\geq \ell+1)=
g(n;\overline{\lambda}_k\leq \ell)+1$ for $0\leq \ell \leq n-\lceil
k/2\rceil-1$.

In order to make the parameter $f(n;\overline{\lambda}_k\leq \ell)$
to be meaningful, we need to determine the range of $\ell$. In fact,
with the help of the definitions of $\overline{\kappa}_k(G)$,
$\kappa_k(G)$, $\overline{\lambda}_k(G)$, $\lambda_k(G)$ and
Theorems \ref{th2-3} and \ref{th2-4}, Li and Mao got the
following observation, which implies that $1\leq \ell \leq n-\lceil
k/2 \rceil$.
\begin{obs}{\upshape \cite{LM2}}\label{obs8-3}
Let $k,n$ be two integers with $3\leq k\leq n$. Then for a connected
graph $G$ of order $n$, $1\leq \overline{\kappa}_k(G)\leq
\overline{\lambda}_k(G)\leq n-\lceil k/2 \rceil$. Moreover, the
upper and lower bounds are sharp.
\end{obs}

Let us now introduce a graph class $\mathcal{G}_n^*$ by a few steps.
For $r\geq 5$,
$\mathcal{G}_{n}=\{H_r^1,H_r^2,H_r^3,H_r^4,\linebreak[2]
H_r^5,H_r^6,H_r^7\}$ is a class of graphs of order $r$ (see Figure
$8.2$ for details).

\begin{figure}[!hbpt]
\begin{center}
\includegraphics[scale=0.75]{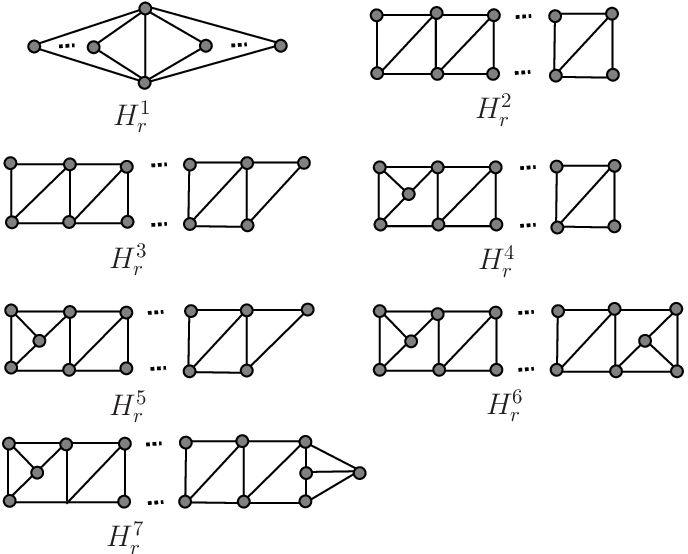}
\end{center}
\begin{center}
\caption{The graph class $\mathcal{G}_n$}
\end{center}\label{fig7}
\end{figure}

Li, Li and Mao introduced a graph operation. Let $H$ be a connected
graph, and $u$ a vertex of $H$. They defined the {\it attaching
operation at the vertex $u$} on $H$ as follows:

$\bullet$ identifying $u$ and a vertex of a $K_4$;

$\bullet$ $u$ is attached with only one $K_4$.

The vertex $u$ is called \emph{an attaching vertex}.

Let $\mathcal{H}_n^{i} \ (1\leq i\leq 7)$ be the class of graphs,
each of them is obtained from a graph $H_r^{i}$ by the attaching
operation at some vertices of degree $2$ on $H_r^{i}$, where $3\leq
r\leq n$ and $1\leq i\leq 7$ (note that $H_n^{i}\in
\mathcal{H}_n^{i}$). $\mathcal{G}_{n}^*$ is another class of graphs
that contains $\mathcal{G}_{n}$, given as follows:
$\mathcal{G}_{3}^*=\{K_3\}$, $\mathcal{G}_{4}^*=\{K_4\}$,
$\mathcal{G}_{5}^*=\{G_1\}\cup (\bigcup_{i=1}^7 \mathcal{H}_5^{i})$,
$\mathcal{G}_{6}^*=\{G_3,G_4\}\cup (\bigcup_{i=1}^7
\mathcal{H}_6^{i})$, $\mathcal{G}_{7}^*=\bigcup_{i=1}^7
\mathcal{H}_7^{i}$, $\mathcal{G}_{8}^*=\{G_2\}\cup (\bigcup_{i=1}^7
\mathcal{H}_8^{i})$, $\mathcal{G}_{n}^*=\bigcup_{i=1}^7
\mathcal{H}_n^{i}$ for $n\geq 9$ (see Figure $8.3$ for details).

\begin{figure}[!hbpt]
\begin{center}
\includegraphics[scale=0.8]{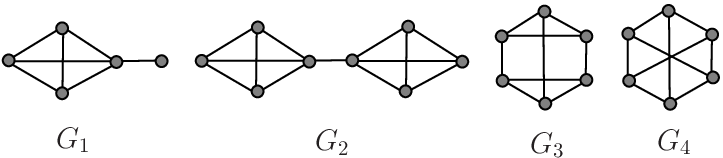}
\end{center}
\begin{center}
\caption{The graphs for Theorem \ref{th8-3}}
\end{center}\label{fig8}
\end{figure}

They obtained the following theorem.

\begin{thm}{\upshape \cite{LLM}}\label{th8-3}
Let $G$ be a connected graph of order $n$ such that
$\overline{\kappa}_3(G)\leq 2$. Then
$$
e(G)\leq \left\{
\begin{array}{ll}
2n-2,&if~n=4;\\
2n-3,&if~n\geq 3,~n\neq 4.
\end{array}
\right.
$$
with equality if and only if $G\in \mathcal{G}_n^*$.
\end{thm}

By the definition of $f(n;\overline{\kappa}_k\leq \ell)$, the
following corollary is immediate.

\begin{cor}{\upshape \cite{LLM}}\label{cro8-4}
$$ f(n;\overline{\kappa}_3\leq 2)=\left\{
\begin{array}{ll}
2n-2&if~n=4;\\
2n-3&if~n\geq 3,~n\neq 4.
\end{array}
\right.
$$
\end{cor}

For a general $\ell$, they constructed a graph class to give a lower
bound of $f(n;\overline{\kappa}_3\leq \ell)$.

\noindent \textbf{Example 8.2.} Let $n,\ell$ be odd, and $G'$ be a
graph obtained from an $(\ell-3)$-regular graph of order $n-2$ by
adding a maximum matching, and $G=G'\vee K_2$. Then
$\delta(G)=\ell-1$, $\overline{\kappa}_3(G)\leq \ell$ and
$e(G)=\frac{\ell+2}{2}(n-2)+\frac{1}{2}$.

Otherwise, let $G'$ be an $(\ell-2)$-regular graph of order $n-2$ and
$G=G'\vee K_2$. Then $\delta(G)=\ell$, $\overline{\kappa}_3(G)\leq
\ell$ and $e(G)=\frac{\ell+2}{2}(n-2)+1$.

Therefore,
$$
f(n;\overline{\kappa}_3\leq \ell)\geq \left\{
\begin{array}{cc}
\frac{\ell+2}{2}(n-2)+\frac{1}{2}&for~n,\ell~odd,\\
\frac{\ell+2}{2}(n-2)+1&~otherwise.
\end{array}
\right.
$$

One can see that for $\ell=2$ this bound is the best possible
($f(n;\overline{\kappa}_3\leq 2)=2n-3$). Actually, the graph
constructed for this bound is $K_2\vee (n-2)K_1$, which belongs to
$\mathcal{G}_{n}^*$.

Li and Zhao \cite{LZ} investigated the exact value of
$f(n;\overline{\kappa}_k=1)$. They introduced the following
operation and graph class: Let $H_1$ and $H_2$ be two connected
graphs. We obtain a graph $H_1+H_2$ from $H_1$ and $H_2$ by joining
an edge $uv$ between $H_1$ and $H_2$ where $u\in H_1$, $v\in H_2$.
We call this operation the \emph{adding operation}.
$\{C_3\}^i+\{C_4\}^j+\{C_5\}^k+\{K_1\}^l$ is a class of connected
graphs obtained from $i$ copies of $C_3$, $j$ copies of $C_4$, $k$
copies of $C_5$ and $\ell$ copies of $K_1$ by some adding operations
such that $0\leq i\leq \lfloor\frac{n}{3}\rfloor$, $0\leq j\leq 2$,
$0\leq k\leq 1$, $0\leq \ell \leq 2$ and $3i+4j+5k+\ell=n$. Note
that these operations are taken in an arbitrary order.

The following graphs shown in Figure $8.4$ will be used later.

\begin{figure}[!hbpt]
\begin{center}\label{fig8-4}
\includegraphics[scale=0.7]{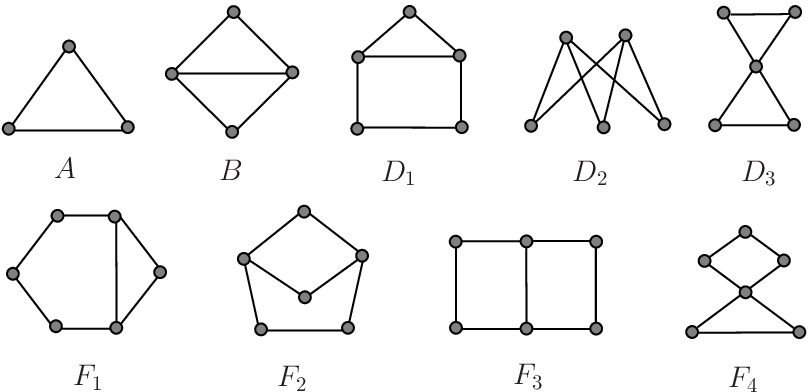}
\end{center}
\begin{center}
\caption{Graphs for $f(n;\overline{\kappa}_k=1)$.}
\end{center}
\end{figure}

At first, they studied the exact value of
$f(n;\overline{\kappa}_3=1)$ and characterized the graphs attaining
this value.

\begin{thm}{\upshape \cite{LZ}}\label{th8-6}
Let $n=3r+q \ (0\leq q \leq 2)$, and let $G$ be a connected graph of
order $n$ such that $\overline{\kappa}_3(G)=1$. Then

$$
e(G)\leq \frac{4n-3-q}{3}
$$
with equality if and only if $G\in \mathcal{G}_n^q$.
\end{thm}

The graph class $\mathcal{G}_n^q$ is defined as follows: Let
$n=3r+q$, $0\leq q \leq 2$. If $q=0$,
$\mathcal{G}_n^{0}=\{C_3\}^{r}$. If $q=1$,
$\mathcal{G}_n^{1}=\{C_3\}^{r}+K_1$ or $\{C_3\}^{r-1}+C_4$. If
$q=2$, $\mathcal{G}_n^{2}=\{C_3\}^{r}+\{K_1\}^2$ or
$\{C_3\}^{r-1}+C_4+K_1$ or $\{C_3\}^{r-1}+C_5$ or
$\{C_3\}^{r-2}+\{C_4\}^2$.

Next, they investigated the exact value of
$f(n;\overline{\kappa}_4=1)$ and characterized the graphs attaining
this value.

\begin{thm}{\upshape \cite{LZ}}\label{th8-7}
Let $n=4r+q$ where $0\leq q \leq 3$, and let $G$ be a connected
graph of order $n$ such that $\overline{\kappa}_4(G)=1$. Then
$$
e(G)\leq \left\{
\begin{array}{cc}
\frac{3n-2}{2}&if~q=0,\\
\frac{3n-3}{2}&if~q=1,\\
\frac{3n-4}{2}&if~q=2,\\
\frac{3n-3}{2}&if~q=3.
\end{array}
\right.
$$
with equality if and only if $G\in \mathcal{H}_n^q$.
\end{thm}

The graph class $\{A\}^{i_0} +\{B\}^{i_1}
+\{D_1\}^{i_2}+\{D_2\}^{i_3}+\{D_3\}^{i_4}
+\{F_1\}^{i_5}+\{F_2\}^{i_6}+\{F_3\}^{i_7}+\{F_4\}^{i_8}+\{K_1\}^{i_9}$
is composed of another connected graph class by some adding
operations satisfying the following conditions:

$\bullet$ $0\leq i_0 \leq 2$, $0\leq i_1 \leq
\lfloor\frac{n}{4}\rfloor$, $0\leq i_2+i_3+i_4 \leq 2$, $0\leq
i_5+i_6+i_7+i_8 \leq 1$, $0\leq i_9 \leq 2$;

$\bullet$ $D_i$ and $F_j$ are not simultaneously in a graph
belonging to this graph class where $1\leq i \leq 3$, $1\leq j \leq
4$;

$\bullet$ $3i_0+4i_1+5(i_2+ i_3+i_4)+6(i_5+i_6+i_7+i_8)+i_9=n$.

The graph class $\mathcal{H}_n^q$ is defined as follows: Let
$n=4r+q$, $0\leq q\leq 3$. If $q=0$, $\mathcal{H}_n^0=\{B\}^r$; If
$q=1$, $\mathcal{H}_n^1=\{B\}^r+K_1$ or $\{B\}^{r-1}+D_i \ (1\leq
i\leq 3$); If $q=2$, $\mathcal{H}_n^2=\{B\}^r+\{K_1\}^2$ or
$\{B\}^{r-1}+\{A\}^2$ or $\{B\}^{r-1}+D_i+K_1$ or
$\{B\}^{r-2}+D_i+D_j \ (1\leq i,j\leq 3)$ or $\{B\}^{r-1}+F_i \
(1\leq i\leq 4)$; If $q=3$, $\mathcal{H}_n^3=\{B\}^r+A$.

For a graph $G$, we say that a path $P=u_1u_2\cdots u_q$ is an {\it
ear} of G if $V(G)\cap V(P)=\{u_1, u_q\}$. If $u_1\neq u_q$, $P$ is
an {\it open ear}; otherwise $P$ is a {\it closed ear}. In their
proofs of Theorems \ref{th8-6} and \ref{th8-7}, they got necessary
and sufficient conditions for $\overline{\kappa}_k(G)=1$ with
$k=3,4$ by means of the number of ears of cycles. Naturally, one
might think that this method can always be applied for $k=5$, i.e.,
every cycle in $G$ has at most two ears, but unfortunately they
found a counterexample.

\noindent \textbf{Example 8.3.} Let $G$ be a graph which contains a
cycle with three independent closed ears. Set $C=u_1u_2u_3$,
$P_1=u_1v_1w_1u_1$, $P_2=u_2v_2w_2u_2$, and $P_3=u_3v_3w_3u_3$.
Then, $\overline{\kappa}_5(G)=1$. In fact, let $S$ be the set of
chosen five vertices. Obviously, for each $i$, if $v_i$ and $w_i$
are in $S$, then $\overline{\kappa}_k(S)=1$. So, only one vertex in
$P_i\setminus u_i$ can be chosen. Suppose that $v_1,v_2,v_3$ have
been chosen. By symmetry, $u_1,u_2$ are chosen. It is easy to check
that there is only one tree connecting $\{u_1,u_2,v_1,v_2,v_3\}$.
The remaining case is that all $u_1$, $u_2$ and $u_3$ are chosen.
Then, no matter which are the another two vertices, only one tree
can be found.

For a general $k$ with $5\leq k\leq n-1$, they obtained the
following lower bound of $f(n;\overline{\kappa}_k(G)=1)$ by
constructing a graph class as follows: If $q=0$, let
$G=\{K_{k-1}\}^r$, then $e(G)=r{k-1 \choose 2}+r-1$. If $1\leq q\leq
k$, let $G=\{K_{k-1}\}^r+K_q$, then $e(G)=r{k-1 \choose 2}+{q
\choose 2}+r$. So the following proposition is immediate.

\begin{pro}{\upshape \cite{LZ}}\label{pro8-8}
For $n=r(k-1)+q \ (0\leq q\leq k-2)$,
$$
f(n;\overline{\kappa}_k=1)\geq \left\{
\begin{array}{ll}
r{k-1 \choose 2}+r-1, & \hbox{if~$q=0$;} \\
r{k-1 \choose 2}+{q \choose 2}+r, & \hbox{if~$1\leq q\leq k-2$.}
\end{array}
\right.
$$
\end{pro}

Actually, Li and Zhao also got the exact value of
$f(n;\overline{\kappa}_k=1)$ for $k=n$.

\begin{thm}{\upshape \cite{LZ}}\label{th8-9}
Let $G$ be a connected graph of order $n$ such that
$\overline{\kappa}_n(G)=1$ where $n\geq 5$. Then
$$
e(G)\leq{n-1 \choose 2}+1
$$
with equality if and only if $G\in \mathcal{K}_n$.
\end{thm}

The graph class $\mathcal{K}_n$ is defined as follows: for $n=5$,
$\mathcal{K}_5=\{G:|V(G)|=5,e(G)=7\}$; for $n\geq 6$,
$\mathcal{K}_n=K_{n-1}+K_1$.

The following three corollaries are immediate from Theorems
\ref{th8-6}, \ref{th8-7} and \ref{th8-9}.

\begin{cor}{\upshape \cite{LZ}}\label{cor8-10}
For $n=3r+q \ (0\leq q \leq 2)$,
$$
f(n;\overline{\kappa}_3=1)=\frac{4n-3-q}{3}
$$.
\end{cor}

\begin{cor}{\upshape \cite{LZ}}\label{cor8-11}
For $n=4r+q \ (0\leq q \leq 3)$,
$$
f(n;\overline{\kappa}_4=1)= \left\{
\begin{array}{cc}
\frac{3n-2}{2}&if~q=0,\\
\frac{3n-3}{2}&if~q=1,\\
\frac{3n-4}{2}&if~q=2,\\
\frac{3n-3}{2}&if~q=3.
\end{array}
\right.
$$
\end{cor}

\begin{cor}{\upshape \cite{LZ}}\label{cor8-12}
For $n\geq 5$, $f(n;\overline{\kappa}_n=1)={n-1 \choose 2}+1$.
\end{cor}

Later, Li and Mao continued to study the above problems. Note that
for $k=n$ we have $1\leq \ell\leq \lfloor\frac{n}{2}\rfloor$ by
Observation \ref{obs8-3}. With the help of Theorem \ref{th1-1} (due
to Nash-Williams and Tutte), they determined the exact value of
$f(n;\overline{\kappa}_k\leq \ell)$ for $k=n$.

\begin{thm}{\upshape \cite{LM2}}\label{th2}
Let $G$ be a connected graph of order $n \ (n\geq 6)$. If
$\overline{\lambda}_n(G)\leq \ell \ (1\leq \ell\leq
\lfloor\frac{n}{2}\rfloor)$, then
$$
e(G)\leq \left\{
\begin{array}{ll}
{{n-1}\choose{2}}+\ell, &if~1\leq \ell\leq
\lfloor\frac{n-4}{2}\rfloor;\\
{{n-1}\choose{2}}+n-2,&if~\ell=
\lfloor\frac{n-2}{2}\rfloor~and~$n$~is~even;\\
{{n-1}\choose{2}}+\frac{n-3}{2},&if~\ell=
\lfloor\frac{n-2}{2}\rfloor~and~$n$~is~odd;\\
{{n}\choose{2}},&if~\ell=\lfloor\frac{n}{2}\rfloor.
\end{array}
\right.
$$
with equality if and only if $G\in \mathcal{G}_n$ for $1\leq
\ell\leq \lfloor\frac{n-4}{2}\rfloor$ where $\mathcal {G}_n$ is a
graph class obtained from a complete graph $K_{n-1}$ by adding a
vertex $v$ and joining $v$ to $\ell$ vertices of $K_{n-1}$;
$G=K_n\setminus e$ where $e\in E(K_n)$ for
$\ell=\lfloor\frac{n-2}{2}\rfloor$ and $n$ even; $G=K_n\setminus M$
where $M\subseteq E(K_n)$ and $|M|=\frac{n+1}{2}$ for
$\ell=\lfloor\frac{n-2}{2}\rfloor$ and $n$ odd; $G=K_n$ for
$\ell=\lfloor\frac{n}{2}\rfloor$.
\end{thm}

From the definition of $f(n;\overline{\kappa}_n\leq \ell)$ and
$g(n;\overline{\lambda}_n\leq \ell)$, the following corollary is
immediate.

\begin{cor}{\upshape \cite{LM2}}
For $1\leq \ell\leq \lfloor\frac{n}{2}\rfloor$ and $n\geq 6$,

$$
f(n;\overline{\kappa}_n\leq \ell)=g(n;\overline{\lambda}_n\leq
\ell)=\left\{
\begin{array}{ll}
{{n-1}\choose{2}}+\ell, &if~1\leq \ell\leq
\lfloor\frac{n-4}{2}\rfloor~or~\ell=
\lfloor\frac{n-2}{2}\rfloor~and~$n$~is~odd;\\
{{n-1}\choose{2}}+2\ell,&if~\ell=
\lfloor\frac{n-2}{2}\rfloor~and~$n$~is~even;\\
{{n}\choose{2}},&if~\ell=\lfloor\frac{n}{2}\rfloor.
\end{array}
\right.
$$
\end{cor}

For $k=n-1$, $1\leq \ell\leq \lfloor\frac{n+1}{2}\rfloor$ by
Observation \ref{obs8-3}. In order to determine the exact value of
$f(n;\overline{\kappa}_{n-1}\leq \ell)$ for a general $\ell \ (1\leq
\ell\leq \lfloor\frac{n+1}{2}\rfloor)$, Li and Mao first focused on
the cases $\ell=\lfloor\frac{n+1}{2}\rfloor$ and
$\lfloor\frac{n-1}{2}\rfloor$. This is also because by
characterizing the graphs with
$\overline{\kappa}_{n-1}(G)=\lfloor\frac{n+1}{2}\rfloor$ and $
\lfloor\frac{n-1}{2}\rfloor$, the difficult case
$\ell=\lfloor\frac{n-3}{2}\rfloor$ can be dealt with. Next, they
considered the case $1\leq \ell\leq \lfloor\frac{n-5}{2}\rfloor$ and
summarized the results for a general $\ell$.

\begin{thm}{\upshape \cite{LM2}}\label{th8-16}
Let $G$ be a connected graph of order $n \ (n\geq 12)$. If
$\overline{\kappa}_{n-1}(G)\leq \ell \ (1\leq \ell\leq
\lfloor\frac{n+1}{2}\rfloor)$, then
$$
e(G)\leq \left\{
\begin{array}{ll}
{{n-2}\choose{2}}+2\ell, &if~1\leq \ell\leq
\lfloor\frac{n-5}{2}\rfloor;\\
{{n-2}\choose{2}}+n-2,&if~\ell=
\lfloor\frac{n-3}{2}\rfloor~and~$n$~is~odd;\\
{{n-2}\choose{2}}+n-4,&if~\ell=
\lfloor\frac{n-3}{2}\rfloor~and~$n$~is~even;\\
{{n-1}\choose{2}}+n-2,&if~\ell=
\lfloor\frac{n-1}{2}\rfloor~and~$n$~is~odd;\\
{{n-1}\choose{2}}+\frac{n-2}{2},&if~\ell=
\lfloor\frac{n-1}{2}\rfloor~and~$n$~is~even;\\
{{n}\choose{2}},&if~\ell=\lfloor\frac{n+1}{2}\rfloor.
\end{array}
\right.
$$
with equality if and only if $G\in \mathcal{H}_n$ for $1\leq
\ell\leq \lfloor\frac{n-5}{2}\rfloor$ where $\mathcal {H}_n$ is a
graph class obtained from the complete graph of order $n-2$ by
adding two nonadjacent vertices and joining each of them to any
$\ell$ vertices of $K_{n-2}$; $G=K_n\setminus M$ where $|M|=n-1$ for
$\ell=\lfloor\frac{n-3}{2}\rfloor$ and $n$ odd; $G\in \mathcal{H}_n$
for $\ell=\lfloor\frac{n-3}{2}\rfloor$ and $n$ even; $G=K_n\setminus
e$ where $e\in E(K_n)$ for $\ell=\lfloor\frac{n-1}{2}\rfloor$ and
$n$ odd; $G=K_n\setminus M$ where $|M|=\frac{n}{2}$ for
$\ell=\lfloor\frac{n-1}{2}\rfloor$ and $n$ even; $G=K_n$ for
$\ell=\lfloor\frac{n+1}{2}\rfloor$.
\end{thm}

The following corollary is immediate from Theorem \ref{th8-16}.

\begin{cor}{\upshape \cite{LM2}}\label{cor8-17}
For $1\leq \ell\leq \lfloor\frac{n+1}{2}\rfloor$ and $n\geq 12$,
$$
f(n;\overline{\kappa}_{n-1}\leq \ell)=\left\{
\begin{array}{ll}
{{n-2}\choose{2}}+2\ell, &if~1\leq \ell\leq
\lfloor\frac{n-5}{2}\rfloor,~or~\ell=
\lfloor\frac{n-3}{2}\rfloor~and~$n$~is~even;\\
{{n-2}\choose{2}}+2\ell+1,&if~\ell=
\lfloor\frac{n-3}{2}\rfloor~and~$n$~is~odd;\\
{{n-1}\choose{2}}+\ell,&if~\ell=
\lfloor\frac{n-1}{2}\rfloor~and~$n$~is~even;\\
{{n-1}\choose{2}}+2\ell-1,&if~\ell=
\lfloor\frac{n-1}{2}\rfloor~and~$n$~is~odd;\\
{{n}\choose{2}},&if~\ell=\lfloor\frac{n+1}{2}\rfloor.
\end{array}
\right.
$$
\end{cor}

Applying Theorem \ref{th8-16} and the relation between
$\overline{\kappa}_{k}$ and $\overline{\lambda}_{k}$, they
investigated the edge case and derived the following result.

\begin{thm}{\upshape \cite{LM2}}\label{th8-18}
Let $G$ be a connected graph of order $n\ (n\geq 12)$. If
$\overline{\lambda}_{n-1}(G)\leq \ell \ (1\leq \ell\leq
\lfloor\frac{n+1}{2}\rfloor)$, then
$$
e(G)\leq \left\{
\begin{array}{ll}
{{n-2}\choose{2}}+2\ell, &if~1\leq \ell\leq
\lfloor\frac{n-5}{2}\rfloor;\\
{{n-2}\choose{2}}+n-2,&if~\ell=
\lfloor\frac{n-3}{2}\rfloor~and~$n$~is~odd;\\
{{n-2}\choose{2}}+n-4,&if~\ell=
\lfloor\frac{n-3}{2}\rfloor~and~$n$~is~even;\\
{{n-1}\choose{2}}+n-2,&if~\ell=
\lfloor\frac{n-1}{2}\rfloor~and~$n$~is~odd;\\
{{n-1}\choose{2}}+\frac{n-2}{2},&if~\ell=
\lfloor\frac{n-1}{2}\rfloor~and~$n$~is~even;\\
{{n}\choose{2}},&if~\ell=\lfloor\frac{n+1}{2}\rfloor.
\end{array}
\right.
$$
with equality if and only if $G\in \mathcal{H}_n$ for $1\leq
\ell\leq \lfloor\frac{n-5}{2}\rfloor$ where $\mathcal {H}_n$ is a
graph class obtained from the complete graph of order $n-2$ by
adding two nonadjacent vertices and joining each of them to any
$\ell$ vertices of $K_{n-2}$; $G=K_n\setminus M$ where $|M|=n-1$ for
$\ell=\lfloor\frac{n-3}{2}\rfloor$ and $n$ odd; $G\in \mathcal{H}_n$
for $\ell=\lfloor\frac{n-3}{2}\rfloor$ and $n$ even; $G=K_n\setminus
e$ where $e\in E(K_n)$ for $\ell=\lfloor\frac{n-1}{2}\rfloor$ and
$n$ odd; $G=K_n\setminus M$ where $|M|=\frac{n}{2}$ for
$\ell=\lfloor\frac{n-1}{2}\rfloor$ and $n$ even; $G=K_n$ for
$\ell=\lfloor\frac{n+1}{2}\rfloor$.
\end{thm}

\begin{cor}{\upshape \cite{LM2}}\label{cor8-19}
For $1\leq \ell\leq \lfloor\frac{n+1}{2}\rfloor$ and $n\geq 12$,
$$
g(n;\overline{\lambda}_{n-1}\leq \ell)=\left\{
\begin{array}{ll}
{{n-2}\choose{2}}+2\ell, &if~1\leq \ell\leq
\lfloor\frac{n-5}{2}\rfloor,~or~\ell=
\lfloor\frac{n-3}{2}\rfloor~and~$n$~is~even;\\
{{n-2}\choose{2}}+2\ell+1,&if~\ell=
\lfloor\frac{n-3}{2}\rfloor~and~$n$~is~odd;\\
{{n-1}\choose{2}}+\ell,&if~\ell=
\lfloor\frac{n-1}{2}\rfloor~and~$n$~is~even;\\
{{n-1}\choose{2}}+2\ell-1,&if~\ell=
\lfloor\frac{n-1}{2}\rfloor~and~$n$~is~odd;\\
{{n}\choose{2}},&if~\ell=\lfloor\frac{n+1}{2}\rfloor.
\end{array}
\right.
$$
\end{cor}

\begin{rem} It is not easy to determine the exact value of
$f(n;\overline{\kappa}_k\leq \ell)$ and
$g(n;\overline{\lambda}_k\leq \ell)$ for a general $k$. So they hoped
to give a sharp lower bound of them. They construct a graph $G$ of
order $n$ as follows: Choose a complete graph $K_{k-1} \ (1\leq
\ell\leq \lfloor\frac{k-1}{2}\rfloor)$. For the remaining $n-k+1$
vertices, join each of them to any $\ell$ vertices of $K_{k-1}$.
Clearly, $\overline{\kappa}_{n-1}(G)\leq
\overline{\lambda}_{n-1}(G)\leq \ell$ and $e(G)=
{{k-1}\choose{2}}+(n-k+1)\ell$. So $f(n;\overline{\kappa}_k\leq
\ell)\geq {{k-1}\choose{2}}+(n-k+1)\ell$ and
$g(n;\overline{\lambda}_k\leq \ell)\geq
{{k-1}\choose{2}}+(n-k+1)\ell$. From Theorems \ref{th8-16} and
\ref{th8-18}, one knows that these two bounds are sharp for
$k=n,n-1$.
\end{rem}

\section{For random graphs}

In this section, we survey the results for random graphs. The two
most frequently occurring probability models of random graphs are
$G(n,M)$ and $G(n,p)$. The first consists of all graphs with $n$
vertices having $M$ edges, in which the graphs have the same
probability. The model $G(n,p)$ consists of all graphs with $n$
vertices in which the edges are chosen independently and with
probability $p$. Given sequences $a_n$ and $b_n$ of real numbers
(possibly taking negative values), we write $a_n=O(b_n)$ if there is
a constant $C>0$ such that $|a_n|\leq C|b_n|$ for all $n$; write
$a_n=o(b_n)$ if $lim_{n\rightarrow \infty} a_n/b_n=0$. Write
$a_n=\Omega(b_n)$ if $a_n\geq 0$ and $b_n=O(a_n)$; $a_n
=\omega(b_n)$ if $a_n\geq 0$ and $b_n=o(a_n)$; $a_n=\Theta(b_n)$ if
$a_n\geq 0$, $a_n=\Omega(b_n)$ and $a_n=\Theta(b_n)$. We say that an
event $\mathcal{A}$ happens \textit{almost surely} if the
probability that it happens approaches $1$ as $n\rightarrow \infty
$, i.e., $Pr[\mathcal{A}]=1-o_n(1)$. Sometimes, we say \textit{a.s.}
for short. We will always assume that $n$ is the variable that tends
to infinity. Given a sequence of events $(E_n)_{n\in N}$, we say
that $E_n$ happens \emph{asymptotically almost surely} $(a.a.s.)$ if
$Pr(E_n)\rightarrow 1$ as $n\rightarrow \infty $.

For a graph property $P$, a function $p(n)$ is called \emph{a
threshold function} of $P$ if:
\begin{itemize}
\item for every $r(n)=O(p(n))$, $G(n, r(n))$ almost surely satisfies $P$ ; and

\item for every $r'(n)=o(p(n))$, $G(n, r'(n))$ almost surely
does not satisfy $P$.
  \end{itemize}

Furthermore, $p(n)$ is called \emph{a sharp threshold function} of
$P$ if there exist two positive constants $c$ and $C$ such that:
\begin{itemize}
\item for every $r(n)\geq C\cdot p(n)$, $G(n,r(n))$ almost surely satisfies $P$ ; and
\item for every $r'(n)\leq c\cdot p(n)$, $G(n,r'(n))$ almost surely
does not satisfy $P$.
\end{itemize}

\subsection{Results for $k=n$}

The spanning tree packing problem has long been one of the main
motives in Graph Theory. Frieze and Luczak \cite{Frieze} firstly
considered the maximum number of edge-disjoint
spanning trees contained in the random graph $G_{k\text{-out}}$, and
studied the random graph $G_k=G_{k\text{-out}}$. This random graph has vertex set
$V=\{1,2,\cdots,n\}$. Each $v\in V_n$ independently chooses a set
$out(v)$ of distinct vertices as neighbours, where each $k$-subset
of $V_n-\{v\}$ is equally likely to be chosen. This produces a
random $k$ out-regular diagraph $D_k$, which has been selected
uniformly from $(\frac{n-1}{k})^n$ distinct possibilities, where
$G_k$ is obtained by ignoring orientation but without coalescing
edges; see \cite{FexnerF, Frieze0} for properties of this model.
They obtained that for a fixed integer $k\geq 2$ the random graph
$G_{k\text{-out}}$ almost surely has $k$ edge-disjoint spanning
trees.

Moreover, Palmer and Spencer \cite{Palmer1} proved that in almost
every random graph process, the hitting time for having $k$
edge-disjoint spanning trees equals the hitting time for having
minimum degree $k$, for any fixed positive integer $k$. In other
words, considering the random graph $G(n,p)$, for any fixed positive
integer $k$, if $p(n)\leq \frac{\log n+k\log\log n-\omega(1)}{n}$
(which is the maximal $p$ for which $\delta(G(n,p))\leq k $ a.s.),
the probability that the spanning tree packing number equals the
minimum degree approaches to $1$ as $n\rightarrow \infty$. On the
other hand, in Catlin's paper \cite{Catlin1} it was found that if
the edge probability was rather large, then almost surely the random
graph $G(n,p)$ has $\lambda_n(G)=\lfloor|E(G)|/(n-1)\rfloor$, which
is less than the minimum degree of $G$. We refer papers
\cite{Catlin1} and \cite{Palmer} to the reader for more details.

A natural question is whether there exists a largest $p(n)$ such
that for every $r'(n)\leq p(n)$, almost surely the random graph
$G(n,p)$ satisfies that the spanning tree packing number equals the
minimum degree.

In \cite{CLL}, Chen, Li and Lian partly answered this question by
establishing the following two theorems for multigraphs. The first
theorem establishes a lower bound of $q(n)$ with $q(n)\geq (1.1\log
n)/n$. Note that this bound for $p$ will allow the minimum degree to
be a function of $n$, and in this sense they improved the result of
Palmer and Spencer.

\begin{thm}{\upshape \cite{CLL}}\label{th9-1}
For any $p$ such that $(\log n+\omega(1))/n\leq p\leq (1.1\log
n)/n$, almost surely the random graph $G\sim G(n,p)$ satisfies that
the spanning tree packing number is equal to the minimum degree,
i.e.
\begin{equation*}
\lim_{n\rightarrow \infty}\mathbf{Pr}(\lambda_n(G)=\delta(G))=1.
\end{equation*}
\end{thm}

The second theorem gives an upper bound of $q(n)$ with $q(n)\leq
(51\log n)/n$.
\begin{thm}{\upshape \cite{CLL}}\label{th9-2}
For any $p$ such that $p\geq (51\log n)/n$, almost surely the random
graph  $G\sim G(n,p)$ satisfies that the spanning tree packing
number is less than the minimum degree, i.e.,
\begin{equation*}
\lim_{n\rightarrow \infty}\mathbf{Pr}(\lambda_n(G)<\delta(G))=1.
\end{equation*}
\end{thm}

\begin{rem} From Theorems \ref{th9-1} and \ref{th9-2}, one can see that
$\log n/n$ is a sharp threshold function for the graph property that
the spanning tree packing number is equal to the minimum degree.
\end{rem}

Later, Gao, P\'{e}rez-Gims\'{e}nez and Sato strengthened the
previous results. In order to introduce their work, we first need
more notations and concepts. Let $\bar{d}(G)=2m(G)/(|V(G)|-1)$. Note
that $\bar{d}(G)$ differs from the average degree of $G$ by a small
factor of $|V(G)|/(|V(G)|-1)$. In particular, in their paper, all
constants involved in these notations do not depend on $p$ under
discussion. For instance, if we have $a_n=\Omega(b_n)$, where $b_n$
may be an expression involving $p=p(n)$, then it means that there
are constants $C>0$ and $n_0$ (both independent with $p$), such that
$a_n\geq C|b_n|$ uniformly for all $n\geq n_0$ and for all $p$ in
the range under discussion. For any graph $G$, let $T(G)$ and $A(G)$
denote the maximum number of edge-disjoint spanning trees in $G$
(possibly $0$ if G is disconnected) and the minimum number of
subforests of G which cover the whole edge set of $G$, respectively.
This number $A(G)$ is known as the \emph{arboricity} of $G$.

They proved that for all $p\in [0, 1]$, the STP number is $a.a.s.$
the minimum between $\delta$ and $m/(n-1)$, where $\delta$ and $m$
respectively denote the minimum degree and the number of edges of
$G(n, p)$.

\begin{thm}{\upshape \cite{Gao}}\label{th9-3}
For every $p = p(n)\in [0,1]$, we have that $a.a.s.$
$$
T(G(n,p))=\min\left\{\delta(G(n, p)),\Big \lfloor \frac{\bar{d}(G(n,
p))}{2} \Big \rfloor \right\}.
$$
\end{thm}

Note that the quantities $\delta$ and $m/(n-1)$ above correspond to
the two trivial upper bounds observed earlier for arbitrary graphs,
so this implies that we can $a.a.s.$ find a best-possible number of
edge-disjoint spanning trees in $G(n, p)$. Their argument uses
several properties of $G(n, p)$ in order to bound the number of
crossing edges between subsets of vertices with
certain restrictions, and then applies
the characterization of the STP number by Tutte and Nash-Williams
stated in Theorem \ref{th1-1}. Moreover, they determined the ranges
of $p$ for which the STP number takes each of these two values:
$\delta$ and $m/(n-1)$. In spite of the fact that the property
$\{\delta\leq m/(n-1)\}$ is not necessarily monotonic with respect
to $p$, they showed that it has a sharp threshold at $p\sim \beta
\log n/n$, where $\beta \approx 6.51778$ is a constant defined in
the following theorem.

\begin{thm}{\upshape \cite{Gao}}\label{th9-4}
Let $\beta=2/\log(e/2)\approx 6.51778$. Then

$(1)$ if $p=\frac{\beta(\log n-\log \log n/2-\omega(1))}{n-1}$, then
$a.a.s$ $\delta(G(n, p))\leq \big \lfloor \frac{\bar{d}(G(n, p))}{2}
\big \rfloor$ and so $T(G(n,p))=\delta(G(n, p))$;

$(2)$ if $p=\frac{\beta(\log n-\log \log n/2+\omega(1))}{n-1}$, then
$a.a.s$ $\delta(G(n, p))> \big \lfloor \frac{\bar{d}(G(n, p))}{2}
\big \rfloor$ and so $T(G(n,p))=\delta(G(n, p))$.
\end{thm}

Below this threshold, the STP number of $G(n, p)$ is $a.a.s.$ equal
to $\delta$; and above the threshold it is $a.a.s.$ $m/(n-1)$. In
particular, this settles the question raised by Chen, Li and Lian
\cite{CLL}.

They further considered the random graph process
$G_0,G_1,\cdots,G_{{n}\choose{2}}$ defined as follows: for each $m
=0,\cdots,{{n}\choose{2}}$, $G_m$ is a graph with vertex set $[n]$;
the graph $G_0$ has no edges; and, for each $1\leq m\leq
{{n}\choose{2}}$, the graph $G_m$ is obtained by adding one new edge
to $G_{m-1}$ chosen uniformly at random among the edges not present
in $G_{m-1}$. Equivalently, we can choose uniformly at random a
permutation $(e_1,\cdots, e_{{{n}\choose{2}}})$ of the edges of the
complete graph with vertex set $[n]$, and define each $G_m$ to be
the graph on vertex set $[n]$ and edges $e_1,\cdots,e_m$.

They also included a stronger version of these results in the
context of the random graph process in which $p$ gradually grows
from $0$ to $1$ (or, similarly, the edges are added one by one).
This provides a full characterization of the STP number that holds
$a.a.s.$ simultaneously during the whole random graph process.

\begin{thm}{\upshape \cite{Gao}}\label{th9-5}
Let $\beta=2/\log(e/2)\approx 6.51778$. The following holds in the
random graph process $G_0,G_1,\cdots,G_{{n}\choose{2}}$.

$(1)$ $a.a.s$ $T(G_m)=\min \{\delta(G_m),\lfloor m/(n-1)\rfloor\}$
for every $0\leq m\leq {{n}\choose{2}}$.

$(2)$ Moreover, for any constant $\epsilon>0$, $a.a.s$

$\bullet$ $\delta(G_m)\leq \lfloor m/(n-1)\rfloor$ for every $0\leq
m\leq \frac{(1-\epsilon)\beta}{2}n \log n$, and

$\bullet$ $\delta(G_m)>\lfloor m/(n-1)\rfloor$ for every
$\frac{(1-\epsilon)\beta}{2}n \log n\leq m\leq {{n}\choose{2}}$.
\end{thm}

The argument combines a more accurate version of the same ideas used
in the analysis of the STP number of $G(n, p)$ together with
multiple couplings of $G(n, p)$ at different values of $p$. In
addition, the article contains several results about the arboricity
of $G(n, p)$. As an almost direct application of their result on the
STP number, for $p$ above the threshold $\beta log n/n$, they
determined the arboricity of $G(n, p)$ to be $a.a.s.$
equal to $m/(n-1)$. This significantly extends the range of $p$ in
the result by Catlin, Chen and Palmer \cite{Catlin1}. They further
proved that for all other values of $p$, the arboricity of $G(n, p)$
is concentrated on at most two values.

\begin{thm}{\upshape \cite{Gao}}\label{th9-6}
Let $\beta=2/\log(e/2)\approx 6.51778$.

$(1)$ For all $p=\frac{\beta(\log n-\log \log n/2-\omega(1))}{n-1}$,
$a.a.s$. $A(G(n, p))=\big \lceil \frac{\bar{d}(G(n, p))}{2} \big
\rceil$; for all $p=\omega(1/n)$, $a.a.s$ $A(G(n, p))\in \left\{
\big \lceil \frac{\bar{d}(G(n, p))}{2} \big \rceil, \big \lceil
\frac{\bar{d}(G(n, p))}{2} \big \rceil+1 \right\}$;

$(2)$ For all $p=\Theta(1/n)$, $a.a.s$ $A(G(n,
p))=(1+\Theta(1))pn/2$. Moreover, there exists a $k>0$ (depending on
$p$), such that $a.a.s.$ $A(G(n, p))\in \{k,k+1\}$.

$(3)$ If $p=o(1/n)$, then $a.a.s.$ $A(G(n, p))\leq 1$.
\end{thm}

In order to prove this for the case $pn\rightarrow \infty$, they
added $o(n)$ edges to $G(n, p)$ in a convenient way that guarantees
a full decomposition of the resulting graph into edge-disjoint
spanning trees. This construction builds upon some of the ideas
previously used to study the STP number. The case $pn=O(1)$ uses
different proof techniques which rely on the structure of the
$k$-core of $G(n, p)$ together with the
Nash-Williams characterization of arboricity stated in Theorem
\ref{th1-2}.

Finally, some of the aforementioned results on the arboricity are
also given below in the more precise context of the random graph
process, similarly as they did for the STP number.

\begin{thm}{\upshape \cite{Gao}}\label{th9-7}
Let $\beta=2/\log(e/2)\approx 6.51778$. The following holds in the
random graph process $G_0,G_1,\cdots,G_{{n}\choose{2}}$.

$(1)$ Let $m_0$ be any function of n such that $m_{0}/n\rightarrow
\infty$ and let $\epsilon>0$ be any constant. Then, $a.a.s.$
simultaneously for all $m\geq m_0$ such that $\delta(G_m)\leq
\bar{d}(G_m)/2$,

$$
\Big \lceil\frac{m+\phi_1}{n-1} \Big\rceil \leq A(G_m)\leq \Big
\lceil\frac{m+\phi_2}{n-1} \Big\rceil,
$$
where $\phi_1=n/\exp(\frac{(1+\epsilon)}{\beta}\frac{2m}{n})=o(n)$
and $\phi_1=n/\exp(\frac{1-\epsilon}{\beta}\frac{2m}{n})=o(n)$. In
particular, $a.a.s.$ $A(G_m)\in
\{\lceil\frac{m}{n-1}\rceil,\lceil\frac{m}{n-1}\rceil+1\}$ for all
$m$ in that range.

$(2)$ Moreover, a.a.s. simultaneously for every m such that
$\delta(G_m)\geq \bar{d}(G_m)/2$ we have

$$
A(G_m)\leq \lfloor m/(n-1)\rfloor.
$$
\end{thm}

\begin{cor}{\upshape \cite{Gao}}\label{cor9-8}
Let $m_{A=i}$ denote the minimum $m$ such that $A(G_m)$ becomes $i$
in the random graph process $G_0,G_1,\cdots,G_{{n}\choose{2}}$. Let
$i_0$ be any function of $n$ such that $i_0\rightarrow \infty$ and
$\epsilon >0$ be a constant. Then $a.a.s.$

$(1)$ for every $i_0\leq i\leq (1-\epsilon)\beta \log n/2$,
$$
(i-1)(n-1)-\phi_2< m_{A=i} < (i-1)(n-1)-\phi_1,
$$
where $\phi_1=n/\exp(\frac{2(1+\epsilon)}{\beta}i)=o(n)$ and
$\phi_2=n/\exp(\frac{2(1-\epsilon)}{\beta})=o(n)$; and

$(2)$ for every $(1+\epsilon)\beta \log n/2\leq i\leq n/2$,
$$
m_{A=i}= (i-1)(n-1)+1
$$
\end{cor}

\subsection{Results for $k=3$}

As well-known, for the vertex connectivity, Bollob\'{a}s and
Thomason \cite{Bollobas5} gave the following result.

\begin{thm}{\upshape \cite{Bollobas5}}\label{th9-9}
If $\ell \in\mathbb{N}$ and $y\in \mathbb{R}$ are fixed,  and
$M=\frac{n}{2}(\log n + \ell\log \log n+y+o(1))\in\mathbb{N}$, then
$$\Pr \left[{\kappa \left({G\left({n,M}\right)}\right)=\ell}
\right] \to 1-{e^{-{e^{-y/\ell!}}}}$$

and $$\Pr \left[{\kappa \left( {G\left( {n,M} \right)}
\right)=\ell+1} \right] \to {e^{-{e^{-y/\ell!}}}}.$$
\end{thm}

Gu, Li and Shi \cite{GLS} focused their attention on the generalized
$3$-connectivity of random graphs for simple graphs. They got the
following theorem, which could be seen as a generalization of
Theorem \ref{th9-9}. At first, they proved that there exists a
constant $c$ such that if $p'<c\frac{{\log n+(\ell+1)\log \log
n-\log \log \log n}}{n}$ then $\kappa_3(G(n, p))<\ell$ almost surely
holds. Then, they showed that for any  three vertices in $G(n, p)$,
where $p=\frac{{\log n + (\ell+1)\log \log n-\log \log \log n}}{n}$,
there almost surely exist three trees of some typical depths rooted
at these three vertices, respectively. Combining some branches of
these trees, $\ell$ internally disjoint trees connecting any three
vertices can be constructed, which implies that $\kappa_3(G(n,
p))\geq \ell$. Hence, they derived the following result.

\begin{thm} {\upshape \cite{GLS}} \label{th9-10}
Let $\ell\geq 1$ be a fixed integer. Then $p=\frac{{\log n +
(\ell+1)\log \log n-\log \log \log n}}{n}$ is a sharp threshold
function for the property $\kappa_3(G(n, p)) \geq \ell$.
\end{thm}

\section{An application problem}

For a network, we usually want to search for a minimum network such
that some local parts have the connectivity we want and the other
parts only need to be connected.

Li, Li and Mao \cite{LLM2} noticed an interesting problem: What is
the smallest number of edges $f(n,k,\ell)$ for a connected graph $G$
of order $n$ that contains $\ell$ edge-disjoint $S$-Steiner trees
for given $k$ vertices of $G$. They determined the
exact value of the parameter $f(n,k,\ell)$ and characterized all the
graphs attaining this value.

This problem has its strong application backgrounds. Suppose that
$G$ is a secure information-gathering network. We denote by $S$ the
set of core departments and each department wants to exchange
important information with others. So we need some Steiner trees to
connect them. But for a vertex not belonging to a Steiner tree, we
let it be an agent. Usually an agent only needs to connect to its
superior leader. So the global network should be connected.

A graph $G$ is called a {\it $(k,\ell)$-minimum connected graph} if
$|V(G)|=n$, $e(G)=f(n,k,\ell)$ and there exist $\ell$ edge-disjoint
$S$-trees for some $S\subseteq V(G)$ with $|S|=k$.

They obtained the exact value of $f(n,k,\ell)$.

\begin{thm}{\upshape \cite{LLM2}}\label{th10-1}
For $3\leq k\leq n$, $f(n,k,\ell)= (k-1)\ell+n-k$.
\end{thm}

Motivated by characterizing the $(k,\ell)$-minimum connected graphs,
they first introduced the notion of an {\it initial graph} and three
graph operations.

\textbf{The initial graph}. Let $k$ and $\ell$ be two integers, and
$K_{1,k-1}$ be a star. Assume that $u_1$ is the center of the star
$K_{1,k-1}$, and $u_2,u_3,\ldots,u_{k-1}$ are the leaves of the
star. A {\it $(k,\ell)$-initial graph} is a graph obtained from the
star $K_{1,k-1}$ by replacing each edge of $K_{1,k-1}$ by $\ell$
multiple edges. Clearly, the graph $G$ contains $\ell$ edge-disjoint
spanning trees $T_1,T_2,\cdots, T_{\ell}$ such that each $T_i$ is a
star $K_{1,k-1}$, that is, $u_1u_2\cup u_1u_3\cup\cup u_1u_{k-1}$.
The edge-disjoint spanning trees $T_1,T_2,\cdots,T_{\ell}$ are
called {\it $(k,\ell)$-initial trees}. Note that
$G=\bigcup_{i=1}^{\ell}T_i$ and $G$ have $k$ vertices and $(k-1)\ell$
edges.

\textbf{Operation I}. For a tree $T_i$, we add an edge $e$ to $T_i$
such that $e\not\in E(T_i)$ and $e$ joins two vertices of $T_i$.
Thus $T_i+e$ contains a unique cycle, say $C$. Pick up an edge
$e'\in E(C) \ (e'\neq e)$, then we obtain a new tree $T'_i$ by deleting
$e'$ from $T_i+e$. Set $T_i:=T_i'$ and $G:=\bigcup_{i=1}^{\ell}T_i$.

\textbf{Operation II}. For a tree $T_i$, pick up a vertex $v\in
V(T_i)$. Let $N=N_{T_i}(v)$. Divide $N$ into three subsets $N_1,
N_2$ and $N_3$ such that $N=N_1\cup N_2\cup N_3$ and $N_3=\{u\}$
(Note that $N_1$ and $N_2$ could be empty sets). Replace the vertex
$v$ by two vertices $v'$ and $v''$, and join $v'$ to each vertex of
$N_1$, $v''$ to each vertex of $N_2$, and $u$ to $v'$ and $v''$. See
Figure $10.1$ for details. If $v\in S$, then set $v:=v'$; otherwise,
we do nothing. Denote the new tree by $T_i'$. Set $T_i:=T_i'$ and
$G:=\bigcup_{i=1}^{\ell}T_i$.

\begin{figure}[!hbpt]
\begin{center}
\includegraphics[scale=0.7]{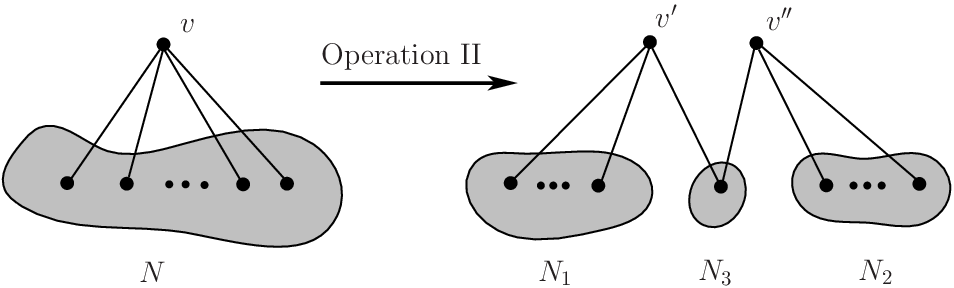}
\end{center}
\begin{center}
\caption{Operation II}
\end{center}\label{fig10-1}
\end{figure}

\textbf{Operation III}. For a tree $T_i$, pick up a vertex $v\in
V(T_i)$. Let $N=N_{T_i}(v)$. Divide $N$ into two subsets $N_1$ and
$N_2$ such that $N=N_1\cup N_2$ (Note that $N_1$ and $N_2$ could be
empty sets). We replace vertex $v$ by two new vertices $v'$ and
$v''$, and join $v'$ to each vertex of $N_1$, $v''$ to each vertex
of $N_2$, and join $v'$ to $v''$. See Figure $10.2$ for details. If
$v\in S$, then set $v:=v'$; otherwise, we do nothing. Set
$T_i:=T_i'$ and $G:=\bigcup_{i=1}^{\ell}T_i$.

\begin{figure}[!hbpt]
\begin{center}
\includegraphics[scale=0.7]{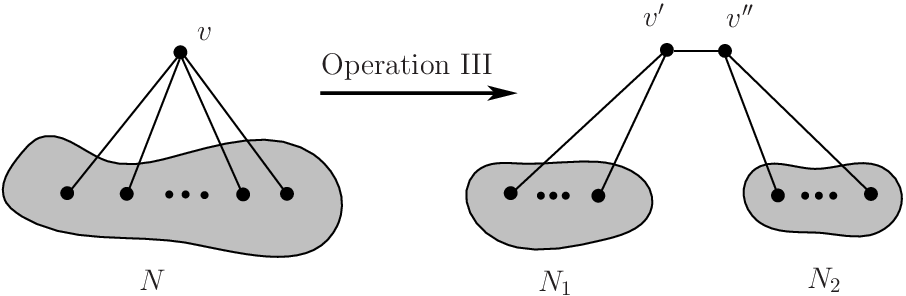}
\end{center}
\begin{center}
\caption{Operation III}
\end{center}\label{fig9-2}
\end{figure}

Next, they proved that any $(k,\ell)$-minimum connected graph can be
obtained from a $(k,\ell)$-initial graph.

\begin{thm}{\upshape \cite{LLM2}}\label{th10-2}
Any $(k,\ell)$-minimum connected graph can be obtained from a
$(k,\ell)$-initial graph by doing a sequence of Operations I, II and
III on the initial graph.
\end{thm}

\section{Concluding remarks}

As we mentioned before, the pendant tree-connectivity and
path-connectivity are also nice and natural generalizations of the
classical connectivity, and they have very close relation with the generalized
(edge-)connectivity. However, in order to make things clear (not to cause confusion),
in this survey we only concentrate on one kind of generalization,
i.e., the generalized (edge-)connectivity. Actually, for the pendant tree-connectivity and
path-connectivity, there have been many results published recently.
For these results we refer the reader to \cite{Hager,
Hager2, Mao1, Mao2, Mao3}.


\begin{thebibliography}{1}

\bibitem{Achuthan}
N. Achuthan, N.R. Achuthan, L. Caccetta, \emph{On the
Nordhaus-Gaddum problems}, Australasian J. Combin. 2(1990), 5--27.

\bibitem{AkersK}
S.B. Akers, B. Krishnamurthy, \emph{A
group-theoretic model for symmetric interconnection networks}, IEEE
Trans. Comput. 38(4)(1989), 555--566.


\bibitem{Alavi}
Y. Alavi, J. Mitchem, \emph{The connectivity and edge-connectivity
of complementary graphs}, Lecture Notes Math. 186(1971), 1--3.

\bibitem{Aouchiche}
M. Aouchiche, P. Hansen, \emph{A survey of Nordhaus-Gaddum type
relations}, Discrete Appl. Math. 161(4-5)(2013), 466--546.

\bibitem{Babai}
L. Babai, \emph{Automorphism groups, isomorphism,
reconstruction, in: R.L. Graham et al. (Eds.)}, Handbook of
Combinatorics, Elsevier Science, Amsterdam, 1995, 1449--1540.

\bibitem{Bao}
F. Bao, Y. Igarashi, S.R. \"{O}hring, \emph{Reliable broadcasting in
product networks}, Discrete Appl. Math. 83(1998), 3--20.

\bibitem{Barden}
B. Barden, R. Libeskind-Hadas, J. Davis, W. Williams, \emph{On
edge-disjoint spanning trees in hypercubes}, Inform. Proc. Lett.
70(1999), 13--16.

\bibitem{Bartfai}
P. B\'{a}rtfai, \emph{Solution of a problem proposed by P.
Erd\"{o}s (in Hungarian)}, Mat. Lapok. (1960), 175--140.

\bibitem{Bauer}
D. Bauer, H. Broersma, E. Schmeichel, \emph{Toughness in graphs: A
survey}, Graphs \& Combin. 22(2006), 1--35.

\bibitem{Beineke1}
L.W. Beineke, R.J. Wilson, \emph{Topics in Structural Graph
Theory}, Cambrige University Press, 2013.

\bibitem{Beineke2}
L.W. Beineke, O.R. Oellermann, R.E. Pippert, \emph{The average
connectivity of a graph}, Discrete Math. 252(2002), 31--45.

\bibitem{Biggs}
N. Biggs, Algebraic Graph Theory, Cambridge
University Press, New York, 1992.

\bibitem{Blasiak}
A. Blasiak, R. Kleinberg, E. Lubetzky, \emph{Lexicographic products
and the power of non-linear network coding}, FOCS (2011), 609--618.

\bibitem{Boesch}
F.T. Boesch, S. Chen, \emph{A generalization of line connectivity
and optimally invulnerable graphs}, SIAM J. Appl. Math. 34(1078),
657--665.

\bibitem{Bollobas1}
B. Bollob\'{a}s, \emph{Extremal Graph Theory},
Acdemic press, 1978.

\bibitem{Bollobas2}
B. Bollob\'{a}s, \emph{On graphs with at most three independent
paths connecting any two vertices}, Studia Sci. Math. Hungar.
1(1966), 137--140.

\bibitem{Bollobas3}
B. Bollob\'{a}s, \emph{Cycles and semi-topological configurations,
in: ``Theory and Applications of graphs'' (Y. Alavi and D.R. Lick,
eds)} Lecture Notes Math. 642, Springer, 1978, 66--74.

\bibitem{Bollobas4}
B. Bollob\'{a}s, \emph{Random Graphs}, Cambridge
University Press, 2001.

\bibitem{Bollobas5}
B. Bollob\'{a}s, A. Thomason, \emph{Random graphs of small order},
in Random Graphs, Ann. Discrete. Math. (1985), 47--97.

\bibitem{bondy}
J.A. Bondy, U.S.R. Murty, \emph{Graph Theory}, GTM 244,
Springer, 2008.

\bibitem{Catlin1}
P. Catlin, Z. Chen, E. Palmer, \emph{On the edge-arboricity of a
random graph}, Ars Combin. 35(A)(1993), 129--134.

\bibitem{Catlin}
P. Catlin, J. Grossman, A. Hobbs, H. Lai, \emph{Fractional
arboricity, strength, and principal partitions in graphs and
matroids}, Discrete Appl. Math. 40(1992), 285--302.

\bibitem{Chartrand1}
G. Chartrand, S. Kappor, L. Lesniak, D. Lick, \emph{Generalized
connectivity in graphs}, Bull. Bombay Math. Colloq. 2(1984), 1--6.

\bibitem{Chartrand2}
G. Chartrand, F. Okamoto, P. Zhang, \emph{Rainbow trees in graphs
and generalized connectivity}, Networks 55(4)(2010), 360--367.

\bibitem{Steeart}
G. Chartrand, M. Steeart, \emph{The connectivity of line graphs},
Math. Ann. 182(1969), 170--174.

\bibitem{CLLM}
L. Chen, X. Li, M. Liu, Y. Mao, \emph{A solution to
a conjecture on the generalized connectivity of graphs},
arXiv:1304.6153 [math.CO] 2013.

\bibitem{CLL}
X. Chen, X. Li, H. Lian, \emph{Note on packing of
edge-disjoint spanning trees in sparse random graphs},
arXiv:1301.1097 [math.CO] 2013.

\bibitem{Cheng}
X. Cheng, D. Du, \emph{Steiner Trees in Industry}, Kluwer Academic
Publisher, Dordrecht, 2001.

\bibitem{Chiue}
W. Chiue, B. Shieh, \emph{On connectivity of Cartesian product of
two graphs}, Appl. Math. Comput. 102(1999), 129--137.

\bibitem{Cunningham}
W. Cunningham, \emph{Optimal attack and reinforcement of a netwok},
J. $ACM$ 32(1985), 549--561.

\bibitem{DayA}
K. Day, A.E. Al-Ayyoub, \emph{The cross product of interconnection
networks}, IEEE Trans. Parallel \& Distributed Systems 8(2)(1997),
109--118.

\bibitem{Day}
D.P. Day, O.R. Oellermann, H.C. Swart, \emph{The $\ell$-connectivity
function of trees and complete multipartite graphs}, J. Combin.
Math. Combin. Comput. 10(1991), 183--192.

\bibitem{Dirac}
G.A. Dirac, \emph{In abstrakten graphen vorhandene
vollst\"{a}ndige 4-graphen und ihre unterteilungen}, Math. Nach
22(1960), 61--85.

\bibitem{Du}
D. Du, X. Hu, \emph{Steiner Tree Problems in Computer Communication
Networks}, World Scientific, 2008.

\bibitem{Esf} A.H. Esfahanian, \emph{Generalized measures of fault tolerance
with application to $N$-cube networks}, IEEE Trans. Computers 38(1989),
1586--1591.

\bibitem{Feng}
M. Feng, M. Xu, K. Wang, \emph{Identifying codes of lexicographic
product of graphs}, Electron. J. Combin. 19(4)(2012), 56--63.

\bibitem{FexnerF}
T.I. Fexner, A.M. Frieze, On the connectivity of
random $m$-orientable graphs and digraphs, Combinatorica 2(1982),
347--369.

\bibitem{Fragopoulou}
P. Fragopoulou, S.G. Akl, \emph{Edge-disjoint spanning trees on the
star network with applications to fault tolerance}, IEEE Trans.
Computers 45(2)(1996), 174--185.

\bibitem{Frank}
A. Frank, T. Kir\'{a}by, M. Kriesell, \emph{On
decomposing a hypergraph into $k$ connected sub-hypergraphs},
Technical Report published by the Egrevdry Research Group, Budapest,
Hungary, ISSN 1587-4451, 2001.

\bibitem{Frieze0}
A.M. Frieze, \emph{Maximum matchings in a class of
random graphs}, J. Comb. Theory B 40(1986), 196--212.

\bibitem{Frieze}
A.M. Frieze, T. Luczak, \emph{Edge-disjoint spanning trees in random
graphs}, Per. Math. Hung. 21(1990), 35--37.

\bibitem{Gao}
P. Gao, X. P\'{e}rez-Gim\'{e}nez, C.M. Sato, \emph{Arboricity and
spanning-tree packing of random graphs}, arXiv:1303.517[math.CO]
2013.

\bibitem{Goldsmith1}
D.L. Goldsmith, \emph{On the second-order edge-connectivity of a
graph}, Congressus Numerantium 29(1980), 479--484.

\bibitem{Goldsmith2}
D.L. Goldsmith, \emph{On the $n$-order connectivity
of a graph}, Congressus Numerantium 32(1981), 375--382.

\bibitem{Goldsmith3}
D.L. Goldsmith, B. Manval, V. Faber, \emph{Seperation of graphs into
three components by removal of edges}, J. Graph Theory 4(1980),
213--218.

\bibitem{Grotschel1}
M. Gr\"{o}tschel, \emph{The Steiner tree packing problem in $VLSI$
design}, Math. Program. 78(1997), 265--281.

\bibitem{Grotschel2}
M. Gr\"{o}tschel, A. Martin, R. Weismantel, \emph{Packing Steiner
trees: A cutting plane algorithm and commputational results}, Math.
Program. 72(1996), 125--145.

\bibitem{GLS}
R. Gu, X. Li, Y. Shi,
\emph{The generalized 3-connectivity of random graphs},
Acta Math. Sin. 57(2)(2014), 321--330.

\bibitem{Gusfield}
D. Gusfield, \emph{Connectivity and edge-disjoint spanning trees},
Infor. Process. Lett. 16(1983), 87--89.

\bibitem{Hager}
M. Hager, \emph{Pendant tree-connectivity}, J. Comb. Theory
38(1985), 179--189.

\bibitem{Hager2}
M. Hager, \emph{Path-connectivity in graphs}, Discrete Math.
59(1986), 53--59.

\bibitem{Hammack}
R. Hammack, W. Imrich, S. Klav\u{z}r,
\emph{Handbook of Product Graphs}, Secend Edition, CRC Press, 2011.

\bibitem{Harary}
F. Harary, \emph{The maximum connectivity of a graph}, Proc. Nat.
Acad. Sci. USA, 1142--1146.

\bibitem{Hellwig}
A. Hellwig, L. Volkmann, \emph{The connectivity of a graph and its
complement}, Discrete Appl. Math. 156(2008), 3325--3328.

\bibitem{Heydemann}
M.C. Heydemann, \emph{Cayley graphs and
interconnection networks}, in: G. Hahn and G. Sabidussi eds., Graph
Symmetry, Kluwer Academic Publishing, Dordrecht, 1997, 167--224.

\bibitem{Hind}
H.R. Hind, O.R. Oellermann, \emph{Menger-type results for three or
more vertices}, Congressus Numerantium 113(1996), 179--204.

\bibitem{Itai}
A. Itai, M. Rodeh, \emph{The multi-tree approach to reliability in
distributed networks}, Infor. \& Comput. 79(1988), 43--59.

\bibitem{Ivchenko}
G. Ivchenko, \emph{The strength of connectivity of a random graph},
Theory Probab. Appl. 18(2009), 396--403.

\bibitem{Jain}
K. Jain, M. Mahdian, M. Salavatipour, \emph{Packing Steiner trees},
in: Proc. 14th $ACM$-$SIAM$ symposium on Discterte
Algorithms, Baltimore, 2003, 266--274.

\bibitem{Klavzar}
S. Klav\v{z}ar, S. \v{S}pacapan, \emph{On the edge-connectivity of
Cartesian product graphs}, Asian-Eur. J. Math. 1(2008), 93--98.

\bibitem{Kriesell1}
M. Kriesell, \emph{Edge-disjoint trees containing some given
vertices in a graph}, J. Combin. Theory Ser. B 88(2003), 53--65.

\bibitem{Kriesell2}
M. Kriesell, \emph{Edge-disjoint Steiner trees in graphs without
large bridges}, J. Graph Theory 62(2009), 188--198.

\bibitem{Kriesell3}
M. Kriesell, \emph{Local spanning trees in graphs and hypergraph
decomposition with respect to edge-connectivity}, Technical Report
257, University of Hannover, 1999.

\bibitem{Ku}
S. Ku, B. Wang, T. Hung, \emph{Constructing edge-disjoint spanning
trees in product networks}, Parallel and Distributed Systems, IEEE
Transactions on parallel and disjoited systems 14(3)(2003),
213--221.

\bibitem{Lakshmivarahan}
S. Lakshmivarahan, J.S. Jwo, S.K. Dhall,
\emph{Symmetry in interconnection networks based on Cayley graphs of
permutation groups: a survey}, Parallel Comput. 19(4)(1993),
361--407.

\bibitem{LA}
R. Laskar, B. Auerbach, \emph{On descomposition of $r$-partite
graphs into edge-disjoint Hamilton circuits}, Discrete Math.
14(1976), 265--258.

\bibitem{Lau}
L. Lau, \emph{An approximate max-Steiner-tree-packing
min-Steiner-cut theorem}, Combinatorica 27(2007), 71--90.

\bibitem{Leonard1}
J. Leonard, \emph{On a conjecture of Bollob\'{a}s and Edr\"{o}s},
Period. Math. Hungar. 3(1973), 281--284.

\bibitem{Leonard2}
J. Leonard, \emph{On graphs with at most four edge-disjoint paths
connecting any two vertices}, J. Comb. Theory Ser. B 13(1972),
242--250.

\bibitem{Leonard3}
J. Leonard, \emph{Graphs with $6$-ways}, Canad. J. Math. 25(1973),
687--692.

\bibitem{LXZW}
F. Li, Z. Xu, H. Zhao, W. Wang, \emph{On the number of spanning
trees of the lexicographic product of networks}, Sci. China Ser. F
42(2012), 949--959.

\bibitem{LLM}
H. Li, X. Li, Y. Mao, \emph{On extremal graphs with at
most two internally disjoint Steiner trees connecting any three
vertices}, Bull. Malays. Math. Sci. Soc. 37(3)(2014), 747--756.

\bibitem{LLM2}
H. Li, X. Li, Y. Mao, \emph{The minimally connected graphs
containing $\ell$ disjoint Steiner trees for given vertices},
Preprint 2012.

\bibitem{LLMS}
H. Li, X. Li, Y. Mao, Y. Sun, \emph{Note on the generalized
connectivity}, Ars Combin. 114(2014), 193--202.

\bibitem{LLMY}
H. Li, X. Li, Y. Mao, J. Yue, \emph{Note on the spanning-tree
packing number of lexicographic product graphs},
Discrete Math. 338(5-6)(2015), 669--673.

\bibitem{LLSun}
H. Li, X. Li, Y. Sun, \emph{The generalied $3$-connectivity of
Cartesian product graphs}, Discrete Math. Theor. Comput. Sci.
14(1)(2012), 43--54.

\bibitem{SLi}
S. Li, \emph{Some topics on generalized connectivity of graphs},
Thesis for Doctor Degree, Nankai University, 2012.

\bibitem{LLL1}
S. Li, W. Li, X. Li, \emph{The generalized connectivity of complete
bipartite graphs}, Ars Combin. 104(2012), 65--79.

\bibitem{LLL2}
S. Li, W. Li, X. Li, \emph{The generalized connectivity of complete
equipartition $3$-partite graphs}, Bull. Malays. Math. Sci. Soc.
(2)37(1)(2014), 103--121.

\bibitem{LL}
S. Li, X. Li, \emph{Note on the hardness of generalized
connectivity}, J. Comb. Optim. 24(2012), 389--396.

\bibitem{LLShi}
S. Li, X. Li, Y. Shi, \emph{The minimal size of a graph with
generalized connectivity $\kappa_3(G)=2$}, Australasian J. Combin.
51(2011), 209--220.

\bibitem{LLZ}
S. Li, X. Li, W. Zhou, \emph{Sharp bounds for the generalized
connectivity $\kappa_3(G)$}, Discrete Math. 310(2010), 2147--2163.

\bibitem{WLi}
W. Li, \emph{On the generalized connectivity of
complete multipartite graphs}, Thesis for Doctor Degree, Nankai
University, 2012.

\bibitem{LM}
X. Li, Y. Mao, \emph{Nordhaus-Gaddum-type results for the generalized
edge-connectivity of graphs}, Discrete Appl. Math.
185(2015), 102--112.

\bibitem{LM2}
X. Li, Y. Mao, \emph{On extremal graphs with at most $\ell$ internally
disjoint Steiner trees connecting any $n-1$ vertices},
Graphs \& Combin., in press.

\bibitem{LM3}
X. Li, Y. Mao, \emph{The generalized $3$-connectivity of
lexicographic product graphs}, Discrete Math. Theor.
Comput. Sci. 16(1)(2014), 339--354.

\bibitem{LM4}
X. Li, Y. Mao, \emph{The minimal size of a graph
with given generalized $3$-edge-connectivity}, Ars
Combin. 118(2015), 63--72.

\bibitem{LMS}
X. Li, Y. Mao, Y. Sun, \emph{On the generalized (edge-)connectivity
of graphs}, Australasian J. Comb. 58(2)(2014), 304--319.

\bibitem{LMW2}
X. Li, Y. Mao, \emph{Graphs with large
generalized (edge-)connectivity}, arXiv: 1305.1089 [math.CO] 2013.

\bibitem{LYZ}
X. Li, J. Yue, Y. Zhao, \emph{The generalized 3-edge-connectivity
of lexicographic product graphs}, in Comb. Optim. Appl., LNCS 8881 (Proc.
COCOA2014, Maui, HI, USA), pp.412--425.

\bibitem{LZ}
X. Li, Y. Zhao, \emph{On graphs with only one Steiner
tree connecting any $k$ vertices}, arXiv:1301.4623[math.CO] 2013.

\bibitem{Libeskind}
R. Libeskind-Hadas, D. Mazzoni, R. Rajagopalan,
\emph{Tree-based multicasting in wormhole-routed irregular
topologies}, Proc. Merged 12th Int'l Parallel Processing Symp. and
the Ninth Symp. Parallel and Distributed Processing, 244-249, Apr.
1998.

\bibitem{LCXu}
M. L\"{u}, G.L. Chen, J. Xu, \emph{On super
edge-connectivity of Cartesian product graphs}, Networks
49(2)(2007), 135--157.

\bibitem{Mader1}
W. Mader, \emph{Ein extremalproblem des zusammenhangin endlichen
graphen}, Math. Z. 131(1973), 223--231.

\bibitem{Mader2}
W. Mader, \emph{Grad und lokaler zusammenhangs von graphen}, Math.
Ann. 205(1973), 9--11.

\bibitem{Mader3}
W. Mader, \emph{\"{U}ber die maximalzahl kantendisjunkter A-Wege},
Arch. Math. 30(1978), 325--336.

\bibitem{Mader4}
W. Mader, \emph{\"{U}ber die maximalzahl kreuzungsfreier H-wege},
Arch. Math. 31(1978), 387-402.

\bibitem{Mao1}
Y. Mao, \emph{Constructing internally disjoint
pendant Steiner trees in Cartesian product networks}, submitted.

\bibitem{Mao2}
Y. Mao, \emph{On the pendant tree-connectivity of
graphs}, submitted.

\bibitem{Mao3}
Y. Mao, \emph{Path connectivity of lexicographical
product graphs}, Int. J. Comput. Math., in press.

\bibitem{Matula}
D. Matula, \emph{Determining edge-connectivity in $O(mn)$},
Proceeding of $28th$ Symp. Foundation Computer Science (1987),
249--251.

\bibitem{Nash}
C.St.J.A. Nash-Williams, \emph{Edge-disjonint spanning trees of
finite graphs}, J. London Math. Soc. 36(1961), 445--450.

\bibitem{Nash2}
C.St.J.A. Nash-Williams, \emph{Decomposition of finite graphs into
forests}, J. London Math. Soc. 39(1964), 12.

\bibitem{Oellermann1}
O.R. Oellermann, \emph{Connectivity and edge-connectivity in graphs:
A survey}, Congessus Numerantium 116(1996), 231--252.

\bibitem{Oellermann2}
O.R. Oellermann, \emph{On the $\ell$-connectivity of a graph}.
Graphs \& Combin. 3(1987), 285--299.

\bibitem{Oellermann3}
O.R. Oellermann, \emph{A note on the $\ell$-connectivity function of
a graph}, Congessus Numerantium 60(1987), 181--188.


\bibitem{Okamoto}
F. Okamoto, P. Zhang, \emph{The tree connectivity of regular
complete bipartite graphs}, J. Combin. Math. Combin. Comput. 74
(2010), 279--293.

\bibitem{OY}
K. Ozeki, T. Yamashita, \emph{Spanning trees: A survey}, Graphs \&
Combin. 27(1)(2011), 1--26.

\bibitem{Palmer}
E. Palmer, \emph{On the spanning tree packing number of a graph: a
survey}, Discrete Math. 230(2001), 13--21.

\bibitem{Palmer1}
E. Palmer, J. Spencer, \emph{Hitting time for $k$ edge-disjoint
spanning trees in a random graph}, Period. Math. Hungar. 31(1995),
151--156.

\bibitem{Peng}
Y. Peng, C. Chen, K. Koh, \emph{On edge-toughness of a
complete $n$-partite graph}, Research Report No. 304, Lee Kong Chian
Centre for Mathematical Reserch, National University of Singapore,
1987.

\bibitem{Peng2}
Y. Peng, T. Tay, \emph{On edge-toughness of a graph}, J. Graph
Theory 17(1993), 233--246.

\bibitem{Petingi}
L. Petingi, J. Rodriguez, \emph{Bounds on the maximum number of
edge-disjoint Steiner trees of a graph}, Congressus Numerantium
145(2000), 43--52.

\bibitem{Rao}
S. Rao, \emph{Graph and its complement}, Proc. Indian Nat. Sci.
Acad. Part A 41(1975), 297--304.

\bibitem{Robertson}
N. Robertson, P. Seymour, \emph{Graph minors $XIII$. The disjoint
path problems}, J. Comb. Theory Ser. B 63(1995), 65--110.

\bibitem{Roskind}
J. Roskind, R. Tarjan, \emph{A Note on Finding Maximum-Cost
Edge-Disjoint Spanning Trees}, Math. Operations Research,
10(2)(1985), 701--708.

\bibitem{Sabidussi}
Sabidussi, \emph{Graphs with given group and given graph theoretical
properties}, Canadian J. Math. 9(1957), 515--525.

\bibitem{Schrijver}
A. Schrijver, \emph{Combinatorial optimization: Polyhedra and efficiency.
Vol. B, Volume 24 of Algorithms and Combinatorics}, Springer-Verlag,
Berlin, 2003.

\bibitem{Sherwani}
N. Sherwani, \emph{Algorithms for $VLSI$
Physical Design Automation}, 3rd Edition, Kluwer Acad. Pub., London,
1999.

\bibitem{Spacapan}
S. \u{S}pacapan, \emph{Connectivity of Cartesian products of
graphs}, Appl. Math. Lett. 21(2008), 682--685.

\bibitem{SL}
Y. Sun, X. Li, \emph{On the difference of two
generalized connectivities of a graph}, accepted by J. Comb.
Optim.

\bibitem{SZ}
Y. Sun, S. Zhou, \emph{Tree connectivities of Caylay
graphs on Abelian groups with small degrees}, accepted by Bull.
Malays. Math. Sci. Soc.

\bibitem{Thomassen}
B. S{\o}rensen, C. Thomassen, \emph{On $k$-rails in graphs}, J.
Comb. Theory 17(1974), 143--159.

\bibitem{Tutte}
W. Tutte, \emph{On the problem of decomposing a graph into $n$
connected factors}, J. London Math. Soc. 36(1961), 221--230.

\bibitem{Wang}
H. Wang and D. Blough, \emph{Construction of edge-disjoint spanning
trees in the torus and application to multicast in wormhole-routed
networks}, Proc. Int'l Conf. Parallel and Distributed Computing
Systems, 1999.

\bibitem{Welsh}
D. Welsh,
\emph{Matroid Theorey}, Academic Press, London, 1976.

\bibitem{West}
D. West, H. Wu, \emph{Packing Steiner trees and $S$-connectors in
graphs}, J. Comb. Theory Ser. B 102(2012), 186--205.

\bibitem{Wilson}
E.L. Wilson, R.L. Hemminger, M.D. Plimmer, \emph{A family
of path properties for graphs}, Math. Ann. 197(1972), 107--122.

\bibitem{Yang}
C. Yang, J. Xu, \emph{Connectivity of lexicographic product and
direct product of graphs}, Ars Combin. 111(2013), 3--12.

\bibitem{Zhou}
S. Zhou, \emph{A class of arc-transitive Cayley graphs as models for
interconnection networks}, SIAM J. Discrete Math. 23(2009),
694--714.



\end{thebibliography}
\end{document}